\documentclass[preprint,10pt]{elsarticle}

\usepackage{mathrsfs,amsmath,amssymb}

\usepackage{mathrsfs,amsmath,amssymb}
\usepackage{mathrsfs,amsmath,amssymb}
\usepackage{amscd} 
\usepackage{relsize,amsmath} 
\newcommand{\lsum}{\mathlarger{\sum}}

\newcommand{\lint}{\mathlarger{\int}}

\DeclareMathOperator{\esssup}{ess \,sup}

\makeatletter
\newcommand{\biggg}[1]{{\hbox{$\left#1\vbox to 20.5pt{}\right.\n@space$}}}
\newcommand{\Biggg}[1]{{\hbox{$\left#1\vbox to 23.5pt{}\right.\n@space$}}}
\newcommand{\bigggg}[1]{{\hbox{$\left#1\vbox to 26.5pt{}\right.\n@space$}}}
\newcommand{\Bigggg}[1]{{\hbox{$\left#1\vbox to 29.5pt{}\right.\n@space$}}}
\newcommand{\biggggg}[1]{{\hbox{$\left#1\vbox to 32.5pt{}\right.\n@space$}}}
\newcommand{\Biggggg}[1]{{\hbox{$\left#1\vbox to 35.5pt{}\right.\n@space$}}}
\newcommand{\bigggggg}[1]{{\hbox{$\left#1\vbox to 38.5pt{}\right.\n@space$}}}
\newcommand{\Bigggggg}[1]{{\hbox{$\left#1\vbox to 41.5pt{}\right.\n@space$}}}
\makeatother

\makeatletter
\@addtoreset{equation}{section}
\makeatother

\usepackage{braket}

\begin{document}

\newtheorem{thm}{Theorem}
\newtheorem{lem}[thm]{Lemma}
\newdefinition{rmk}{Remark}
\newproof{pf}{Proof}
\newproof{pot}{Proof of Theorem \ref{thm2}}

\begin{frontmatter}



\title{Global structure of 
solutions toward \\
the rarefaction waves 
for the Cauchy problem \\
of the scalar conservation law 
with nonlinear viscosity}


\author[label1,label2]{Natsumi Yoshida}
\ead{14v00067@gst.ritsumei.ac.jp/jt-bnk68@mail.doshisha.ac.jp}

\address[label1]{OIC Research Organization, 
Ritsumeikan University, Ibaraki, Osaka 567-8570, Japan.}
\address[label2]{Faculty of Culture and Information Science, 
Doshisha University, Kyotanabe, Kyoto 610-0394, Japan.}
\address{}

\begin{abstract}
In this paper, we investigate the global structure of solutions 
to the Cauchy problem 
for the scalar viscous conservation law 
where the far field states are prescribed. 
Especially, we deal with the case 
when 
the viscous/diffusive flux $\sigma(v) \sim |\, v \,|^p$
is of non-Newtonian type (i.e., $p>0$), including a pseudo-plastic case (i.e., $p<0$). 
When the corresponding Riemann problem 
for the hyperbolic part 
admits a Riemann solution which 
consists of single rarefaction wave, 
under a condition on nonlinearity of the viscosity,
it has been recently proved by Matsumura-Yoshida \cite{matsumura-yoshida'} that 
the solution of the Cauchy problem tends toward the rarefaction wave 
as time goes to infinity for the case $p>3/7$ without any smallness conditions. 
The new ingredients we obtained 
are the extension to the stability results in \cite{matsumura-yoshida'} to 
the case $p > 1/3$ (also without any smallness conditions), 
and furthermore their precise time-decay estimates.
\end{abstract}


\end{frontmatter}

%



\pagestyle{myheadings}
\thispagestyle{plain}
\markboth{N. YOSHIDA}{GLOBAL STRUCTURE FOR RAREFACTION WAVE}

\section{Introduction and main theorems}
In this paper, 
we consider the asymptotic behavior of solutions to the Cauchy problem 
for a one-dimensional scalar conservation law with nonlinear viscosity 
\begin{eqnarray}
 \left\{\begin{array}{ll}
  \partial_tu +\partial_x \big( \, f(u) - \sigma(\partial_xu )\, \big) =0
  \qquad &\big( \, t>0, \: x\in \mathbb{R} \, \big), \\[5pt]
  u(0, x) = u_0(x) \rightarrow  u_{\pm}  \qquad &( x \rightarrow  \pm \infty ).   
 \end{array}
 \right.\,
\end{eqnarray}
Here, $u=u(t, \, x)$ is the unknown function of $t>0$ and $x\in \mathbb{R}$, 
the so-called conserved quantity, 
the functions $f$ and $-\sigma$ stand for the convective flux and 
viscous/diffusive one, respectively, 
$u_0$ is the initial data, 
and $u_{\pm } \in \mathbb{R}$ 
are the prescribed far field states. 
We suppose that $f$ is a smooth function, 
and $\sigma$ is a smooth 
function
satisfying 
\begin{equation}
\sigma(0)=0,\quad \sigma'(v)>0\quad (v\in \mathbb{R}),
\end{equation}
and 
for some $p>0$ 
\begin{equation}
|\,\sigma(v)\,| \sim | \, v \,|^p,
\quad |\,\sigma'(v)\,| \sim | \, v \,|^{p-1}\quad (\,| \, v \,| \to \infty \,).
\end{equation}
A typical example of $\sigma$ in the field of viscous fluid is
\begin{equation}
\sigma(\partial_xu) = \mu \, 
    \left( \, 1 + | \, \partial_xu \, |^2 \, \right)^{\frac{p-1}{2}} \partial_xu. 
\end{equation}
Here $\mu$ is a positive constant and the so-called viscous coefficient, 
which describes a nonlinear
relation between the internal stress $\sigma$ 
and the deformation velocity $\partial_xu$ ($u$ is the fluid velocity). 
It should be noted that
the cases $p>1$, $p=1$ and $p<1$ physically correspond to 
where the fluid is of dilatant type, Newtonian and pseudo-plastic type,
respectively (see \cite{chh1}, \cite{chh2}, 
\cite{chh-ric}, \cite{jah-str-mul}, \cite{liep-rosh}, 
\cite{ma}, \cite{ma-pr-st}, \cite{soc} and so on).
We are interested in the global structure 
(existence, uniqueness, smoothness and time-decay properties) 
for the solution of 
(1.1), in particular, for the pseudo-plastic case $p<1$. 
In this case, recently, Matsumura-Yoshida \cite{matsumura-yoshida'} 
first showed in the case $p>3/7$ by a technical energy method 
that the following: when $u_- = u_+\,(=:\tilde{u})$, then 
the solution of (1.1) globally tends toward the constant 
state $\tilde{u}$ as time goes to infinity, 
and when $u_- < u_+$ and $f''(u)>0\; (u \in \mathbb{R})$, then 
the solution of (1.1) a rarefaction wave $u^r$ (see (1.7)) 
as time goes to infinity. 
However, the global stabilities were open in the case $0<p<3/7$, 
and their time-decay rates were even open in case $0<p<1$. 
Under such a situation, we have succeeded in 
proving these global stabilities in the case $p>1/3$, 
and moreover, obtaining their precise time-decay estimates in this case. 

In fact, for the case $u_- = u_+ =\tilde{u}$, we can show the following. 

\medskip

\noindent
{\bf Theorem 1.1.}\quad{\it
Assume the far field states satisfy $u_- = u_+ =\tilde{u}$, 
the convective flux $f\in C^3(\mathbb{R})$, 
the viscous flux 
$\sigma\in C^2(\mathbb{R})$, {\rm(1.2)}, {\rm(1.3)} 
and $p>1/3$. 
Further assume the initial data satisfy
$u_0-\tilde{u} \in H^2$.
Then the Cauchy problem {\rm(1.1)} has a 
unique global in time 
solution $u$ 
satisfying 
\begin{eqnarray*}
\left\{\begin{array}{ll}
u-\tilde{u} \in C^0\cap L^{\infty}
( [\, 0, \, \infty \, ) \, ; H^2 \,),\\[5pt]
\partial _x u \in L^2\bigl( \, 0,\, \infty \, ; H^2 \,\bigr),\\[5pt]
\partial _t u \in C^0
                  \cap L^{\infty}\bigl( \, [\, 0,\, \infty \, ) \, ; L^2 \,\bigr)
                  \cap L^2\bigl( \, 0,\, \infty \, ; H^1 \,\bigr),
\end{array} 
\right.\,
\end{eqnarray*}
and the asymptotic behavior 
\begin{eqnarray*}
\left\{\begin{array}{ll}
\displaystyle{
\lim _{t \to \infty}\sup_{x\in \mathbb{R}} \, 
\left| \, u(t,x) - \tilde{u} \, \right| = 0, 
\quad 
\lim _{t \to \infty}\sup_{x\in \mathbb{R}} \, 
\left| \, \partial_{x}u(t,x) \, \right| = 0
},\\[10pt]
\displaystyle{
\lim _{t \to \infty} \, 
\| \, \partial_{x}u(t) \, \|_{H^1} = 0, 
\quad 
\lim _{t \to \infty} \, 
\| \, \partial_{t}u(t) \, \|_{L^2} = 0. 
}
\end{array} 
\right.\,
\end{eqnarray*}
}

\medskip 

The decay results correspond to Theorem 1.1 are the following. 

\medskip

\noindent
{\bf Theorem 1.2.} \quad{\it
Under the same assumptions as in Theorem 1.1, 
the unique global in time 
solution $u$ 
of the Cauchy problem {\rm(1.1)} 
has the following time-decay estimates 
\begin{eqnarray*}
\left\{\begin{array}{ll}
\| \, u(t)-\tilde{u} \, \|_{L^q}
\le C_{u_{0}} \, 
(1+t)^{-\frac{1}{4}\left( 1-\frac{2}{q} \right)}
\quad \big( \, t\ge0 \, \big),\\[5pt]
\| \, u(t)-\tilde{u} \, \|_{L^{\infty}}
\le C_{u_{0}, \epsilon} \, (1+t)^{-\frac{1}{4}+\epsilon}
\quad \big( \, t\ge0 \, \big),
\end{array} 
\right.\,
\end{eqnarray*}
for $q \in [\, 2, \, \infty)$ and any $\epsilon>0$. 
}

\medskip

\noindent
{\bf Theorem 1.3.} \quad{\it
Under the same assumptions as in Theorem 1.1, 
if the initial data further satisfies $u_0-\tilde{u}\in L^1$, 
then the unique global in time 
solution $u$ 
of the Cauchy problem {\rm(1.1)} 
has the following time-growth and time-decay estimates 
\begin{eqnarray*}
\left\{\begin{array}{ll}
\| \, u(t)-\tilde{u} \, \|_{L^1}
\le C_{u_{0}} \, 
\max \big\{ \, 1, \, \ln (1+t) \, \big\}
\quad \big( \, t\ge0 \, \big),\\[5pt]
\| \, u(t)-\tilde{u} \, \|_{L^q}
\le C_{u_{0}} \, 
(1+t)^{-\frac{1}{2}\left( 1-\frac{1}{q} \right)} \, 
\max \big\{ \, 1, \, \ln (1+t) \, \big\}
\quad \big( \, t\ge0 \, \big),\\[5pt]
\| \, u(t)-\tilde{u} \, \|_{L^{\infty}}
\le C_{u_{0}, \epsilon} \, (1+t)^{-\frac{1}{2}+\epsilon}
\quad \big( \, t\ge0 \, \big),
\end{array} 
\right.\,
\end{eqnarray*}
for $q \in ( 1, \, \infty)$ and any $\epsilon>0$. 
}

\medskip

\noindent
{\bf Theorem 1.4.} \quad{\it
Under the same assumptions as in Theorem 1.1, 
the unique global in time 
solution $u$ 
of the Cauchy problem {\rm(1.1)} 
has the following time-decay estimates 
for the derivatives 
\begin{eqnarray*}
\left\{\begin{array}{ll}
\| \, \partial_{x}u(t) \, \|_{L^{q+1}}
\le C_{u_{0}, \epsilon} \, 
(1+t)^{-\frac{2q+1}{2q+2}+\epsilon}
\quad \big( \, t\ge0 \, \big),\\[5pt]
\| \, \partial_{x}u(t) \, \|_{L^{\infty}}
\le C_{u_{0}, \epsilon} \, 
(1+t)^{-1+\epsilon} 
\quad \big( \, t\ge0 \, \big),\\[5pt]
\| \, \partial_{t}u(t) \, \|_{L^{2}}
\le C_{u_{0}, \epsilon} \, (1+t)^{-\frac{3}{4}+\epsilon}
\quad \big( \, t\ge0 \, \big),\\[5pt]
\| \, \partial_{x}^2u(t) \, \|_{L^{2}}
\le C_{u_{0}, \epsilon} \, (1+t)^{-\frac{3}{4}+\epsilon}
\quad \big( \, t\ge0 \, \big),
\end{array} 
\right.\,
\end{eqnarray*}
for $q \in ( 1, \, \infty)$ and any $\epsilon>0$. 
}

\medskip

Next, we consider the case 
where the convective flux function $f$ is fully convex, 
that is, 
\begin{equation}
f''(u)>0\quad (u \in \mathbb{R}), 
\end{equation}
and $u_-<u_+$.
Then, since the corresponding Riemann problem (cf. \cite{lax}, \cite{smoller}) 
\begin{eqnarray}
 \left\{\begin{array} {ll}
 \partial _t u + \partial _x \bigl( \, f(u) \, \bigr)=0 
\quad \big( \, t>0,\: x\in \mathbb{R} \, \big),\\[5pt]
u(0, x)=u_0 ^{\rm{R}} (x)
      := \left\{\begin{array} {ll}
         u_-  & \; (x < 0),\\[5pt]
         u_+  & \; (x > 0),
         \end{array}\right.
 \end{array}
  \right.\,
\end{eqnarray}
turns out to admit a single
rarefaction wave solution, 
we expect that 
the solution of the Cauchy problem (1.1) globally tends toward
the rarefaction wave 
as time goes to infinity. 
Here, the rarefaction wave connecting $u_-$ to $u_+$ 
is given by 
\begin{equation}
u^r \left( \, \frac{x}{t}\: ;\:  u_- ,\:  u_+ \,\right)
= \left\{
\begin{array}{ll}
  u_-  & \; \bigl(\, x \leq f'(u_-)\,t \, \bigr),\\[7pt]
  \displaystyle{ (f')^{-1}\left( \frac{x}{t}\right) } 
  & \; \bigl(\, f'(u_-)\,t \leq x \leq f'(u_+)\,t\,  \bigr),\\[7pt]
   u_+ & \; \bigl(\, x \geq f'(u_+)\,t \, \bigr).
\end{array}
\right. 
\end{equation} 
Then we can show the following 

\medskip

\noindent
{\bf Theorem 1.5.} \quad{\it
Assume the far field states satisfy $u_- < u_+$, 
the convective flux $f\in C^4(\mathbb{R})$, {\rm(1.5)}, 
the viscous flux $\sigma\in C^2(\mathbb{R})$, {\rm(1.2)}, {\rm(1.3)}, 
and $p>1/3$. 
Further assume the initial data satisfy
$u_0-u_0 ^{\rm{R}} \in L^2$ and
$\partial _xu_0 \in H^1$. 
Then the Cauchy problem {\rm(1.1)} has a 
unique global in time 
solution $u$ 
satisfying 
\begin{eqnarray*}
\left\{\begin{array}{ll}
u-u_0 ^{\rm{R}} \in C^0\cap L^{\infty}
\bigl( \, [\, 0, \, \infty \, ) \, ; L^2 \, \bigr),\\[5pt]
\partial _x u \in C^0
                  \cap L^{\infty}\bigl( \, [\, 0,\, \infty \, ) \, ; H^1 \,\bigr)
                  \cap L^2_{\rm{loc}}\bigl( \, 0,\, \infty \, ; H^2 \,\bigr),\\[5pt]
\partial _t u \in C^0
                  \cap L^{\infty}\bigl( \, [\, 0,\, \infty \, ) \, ; L^2 \,\bigr)
                  \cap L^2_{\rm{loc}}\bigl( \, 0,\, \infty \, ; H^1 \,\bigr),
\end{array} 
\right.\,
\end{eqnarray*}
and the asymptotic behavior 
\begin{eqnarray*}
\left\{\begin{array}{ll}
\displaystyle{
\lim _{t \to \infty}\sup_{x\in \mathbb{R}} \, 
\left| \, 
u(t,x) - u^r\left(\, \frac{x}{t}\: ;\: u_-, \, u_+ \, \right) 
\, \right|} = 0,\\[10pt]
\displaystyle{
\lim _{t \to \infty}\sup_{x\in \mathbb{R}} \, 
| \, 
\partial_{x}u(t,x) 
- \partial_{x}u^r\left(\, t, \, x\: ;\: u_-, \, u_+ \, \right)
\, |} = 0, \\[10pt]
\displaystyle{
\lim _{t \to \infty} \, 
\| \, 
\partial_{x}u(t) 
- \partial_{x}u^r\left(\, t, \, \cdot \: ;\: u_-, \, u_+ \, \right) 
\, \|_{H^{1}}
} = 0,\\[5pt]
\displaystyle{
\lim _{t \to \infty} \, 
\| \, 
\partial_{t}u(t) 
- \partial_{t}u^r\left(\, t, \, \cdot \: ;\: u_-, \, u_+ \, \right) 
\, \|_{L^{2}}
} = 0.
\end{array} 
\right.\,
\end{eqnarray*}
where $\partial_{x}u^r$ and $\partial_{t}u^r$ are given by 
\begin{equation*}
\partial_{x}u^r\left(\, t, \, x \: ;\: u_-, \, u_+ \, \right)
= \left\{
\begin{array}{ll}
  0  & \; \bigl(\, x \leq f'(u_-)\,t \, \bigr),\\[7pt]
  \displaystyle{ 
  \frac{1}
  {f''\bigg( \, (f')^{-1}\left( \displaystyle{\frac{x}{t}} \right) \, \bigg)} 
  \, \frac{1}{t}
  } 
  & \; \bigl(\, f'(u_-)\,t \leq x \leq f'(u_+)\,t\,  \bigr),\\[15pt]
  0 & \; \bigl(\, x \geq f'(u_+)\,t \, \bigr),
\end{array}
\right. 
\end{equation*} 
and 
\begin{equation*}
\partial_{t}u^r\left(\, t, \, x \: ;\: u_-, \, u_+ \, \right)
= \left\{
\begin{array}{ll}
  0  & \; \bigl(\, x \leq f'(u_-)\,t \, \bigr),\\[7pt]
  \displaystyle{ 
  \frac{-1}
  {f''\bigg( \, (f')^{-1}\left( \displaystyle{\frac{x}{t}} \right) \, \bigg)} 
  \, \frac{x}{t^2}
  } 
  & \; \bigl(\, f'(u_-)\,t \leq x \leq f'(u_+)\,t\,  \bigr),\\[15pt]
  0 & \; \bigl(\, x \geq f'(u_+)\,t \, \bigr),
\end{array}
\right. 
\end{equation*} 
respectively. 
}

\medskip 

The decay results correspond to Theorem 1.5 are the following 
(see also Remarks 2.3 and 3.9). 

\medskip

\noindent
{\bf Theorem 1.6.} \quad{\it
Under the same assumptions as in Theorem 1.5, 
the unique global in time 
solution $u$ 
of the Cauchy problem {\rm(1.1)} 
has the following time-decay estimates 
\begin{eqnarray*}
\left\{\begin{array}{ll}
\displaystyle{
\Big\| \, 
u(t) - u^r\left(\, \frac{\cdot }{t}\: ;\: u_-, \, u_+ \, \right) 
\, \Big\|_{L^q}
}
\le C_{u_{0}} \, 
(1+t)^{-\frac{1}{4}\left( 1-\frac{2}{q} \right)}
\quad \big( \, t\ge0 \, \big),\\[10pt]
\displaystyle{
\Big\| \, 
u(t) - u^r\left(\, \frac{\cdot }{t}\: ;\: u_-, \, u_+ \, \right) 
\, \Big\|_{L^{\infty}}
}
\le C_{u_{0}, \epsilon} \, (1+t)^{-\frac{1}{4}+\epsilon}
\quad \big( \, t\ge0 \, \big),
\end{array} 
\right.\,
\end{eqnarray*}
for $q \in [\, 2, \, \infty)$ and any $\epsilon>0$. 
}

\medskip

\noindent
{\bf Theorem 1.7.} \quad{\it
Under the same assumptions as in Theorem 1.5, 
if the initial data further satisfies $u_0-u_0 ^{\rm{R}}\in L^1$, 
then the unique global in time 
solution $u$ 
of the Cauchy problem {\rm(1.1)} 
has the following time-growth and time-decay estimates 
\begin{eqnarray*}
\left\{\begin{array}{ll}
\displaystyle{
\Big\| \, 
u(t) - u^r\left(\, \frac{\cdot }{t}\: ;\: u_-, \, u_+ \, \right) 
\, \Big\|_{L^1}
}
\le C_{u_{0}, \epsilon} \, 
(1+t)^{\epsilon}
\quad \big( \, t\ge0 \, \big),\\[10pt]
\displaystyle{
\Big\| \, 
u(t) - u^r\left(\, \frac{\cdot }{t}\: ;\: u_-, \, u_+ \, \right) 
\, \Big\|_{L^q}
}\\[10pt]
\quad 
\le C_{u_{0}} \, 
(1+t)^{-\frac{1}{2}\left( 1-\frac{1}{q} \right)} \,
\max \big\{ \, 1, \, \ln (1+t) \, \big\}
\quad \big( \, t\ge0 \, \big),\\[10pt]
\displaystyle{
\Big\| \, 
u(t) - u^r\left(\, \frac{\cdot }{t}\: ;\: u_-, \, u_+ \, \right) 
\, \Big\|_{L^{\infty}}
}
\le C_{u_{0}, \epsilon} \, (1+t)^{-\frac{1}{2}+\epsilon}
\quad \big( \, t\ge0 \, \big),
\end{array} 
\right.\,
\end{eqnarray*}
for $q \in ( 1, \, \infty)$ and any $\epsilon>0$. 
}

\medskip

\noindent
{\bf Theorem 1.8.} \quad{\it
Under the same assumptions as in Theorem 1.5, 
the unique global in time 
solution $u$ 
of the Cauchy problem {\rm(1.1)} 
has the following time-decay estimates 
for the derivatives 
\begin{eqnarray*}
\left\{\begin{array}{ll}
\| \, 
\partial_{x}u(t) 
- \partial_{x}u^r\left(\, 1+t, \, \cdot \: ;\: u_-, \, u_+ \, \right) 
\, \|_{L^{q+1}}
\le C_{u_{0}} \, 
(1+t)^{-\frac{q}{q+1}}
\quad \big( \, t\ge0 \, \big),\\[5pt]
\| \, 
\partial_{x}u(t) 
- \partial_{x}u^r\left(\, 1+t, \, \cdot \: ;\: u_-, \, u_+ \, \right) 
\, \|_{L^{\infty}}
\le C_{u_{0}, \epsilon} \, 
(1+t)^{-1+\epsilon} 
\quad \big( \, t\ge0 \, \big),\\[5pt]
\| \, 
\partial_{t}u(t) 
- \partial_{t}u^r\left(\, 1+t, \, \cdot \: ;\: u_-, \, u_+ \, \right) 
\, \|_{L^{2}}
\le C_{u_{0}} \, (1+t)^{-\frac{1}{2}}
\quad \big( \, t\ge0 \, \big),\\[5pt]
\| \, 
\partial_{x}^2u(t) 
- \partial_{x}^2u^r\left(\, 1+t, \, \cdot \: ;\: u_-, \, u_+ \, \right) 
\, \|_{L^{2}}
\le C_{u_{0}, \epsilon} \, (1+t)^{-\frac{3}{4}+\epsilon}
\quad \big( \, t\ge0 \, \big),\\[5pt]
\| \, 
\partial_{x}u(t) 
- \partial_{x}U^r\left(\, t, \, \cdot \: ;\: u_-, \, u_+ \, \right) 
\, \|_{L^{q+1}}
\le C_{u_{0}, \epsilon} \, 
(1+t)^{-\frac{2q+1}{2q+2}+\epsilon}
\quad \big( \, t\ge0 \, \big),\\[5pt]
\| \, 
\partial_{x}u(t) 
- \partial_{x}U^r\left(\, t, \, \cdot \: ;\: u_-, \, u_+ \, \right) 
\, \|_{L^{\infty}}
\le C_{u_{0}, \epsilon} \, 
(1+t)^{-1+\epsilon} 
\quad \big( \, t\ge0 \, \big),\\[5pt]
\| \, 
\partial_{t}u(t) 
- \partial_{t}U^r\left(\, t, \, \cdot \: ;\: u_-, \, u_+ \, \right) 
\, \|_{L^{2}}
\le C_{u_{0}, \epsilon} \, (1+t)^{-\frac{3}{4}+\epsilon}
\quad \big( \, t\ge0 \, \big),\\[5pt]
\| \, 
\partial_{x}^2u(t) 
- \partial_{x}^2U^r\left(\, t, \, \cdot \: ;\: u_-, \, u_+ \, \right) 
\, \|_{L^{2}}
\le C_{u_{0}, \epsilon} \, (1+t)^{-\frac{3}{4}+\epsilon}
\quad \big( \, t\ge0 \, \big),
\end{array} 
\right.\,
\end{eqnarray*}
for $q \in [\, 1, \, \infty)$ and any $\epsilon>0$, 
where 
$\partial_{x}^2u^r$ is given by 
\begin{equation*}
\partial_{x}^2u^r\left(\, t, \, x \: ;\: u_-, \, u_+ \, \right)
= \left\{
\begin{array}{ll}
  0  & \; \bigl(\, x \leq f'(u_-)\,t \, \bigr),\\[7pt]
  \displaystyle{ 
  \frac{-f'''\bigg( \, (f')^{-1}\left( \displaystyle{\frac{x}{t}} \right) \, \bigg)}
  {\bigg( \, 
  f''\bigg( \, (f')^{-1}\left( \displaystyle{\frac{x}{t}} \right) \, \bigg)
  \, \bigg)^3} 
  \, \frac{1}{t^2}
  } 
  & \; \bigl(\, f'(u_-)\,t \leq x \leq f'(u_+)\,t\,  \bigr),\\[25pt]
  0 & \; \bigl(\, x \geq f'(u_+)\,t \, \bigr),
\end{array}
\right. 
\end{equation*} 
and $U^r$ is a smooth approximation for $u^r$ 
which is defined by {\rm(3.1)}. 
}

\medskip

Now let us review the known results concerning the nonlinear stabilities of 
various of nonlinear waves to the Cauchy problem (1.1). 
For the case $p=1$ (Newtonian type viscosity), global 
stability of both rarefaction wave and viscous shock wave 
were first obtained by
Il'in-Ole{\u\i}nik \cite{ilin-oleinik}. 
Harabetian \cite{har} further obtained 
in the case $p=1$ (and more general case (1.8)) that 
the precise time-decay estimates of global 
stability of single rarefaction wave 
if the far field states satisfy $0\le u_{-}<u_{+}$  or $0<u_{-}<u_{+}$, 
with the aid of the arguments 
on monotone semigroups 
associated with the following type quasilinear parabolic equations 
\begin{equation}
\partial_{t}u
+ \partial_{x} \big( \, f(u)-A'(u)\, \partial_{x}u \, \big)=0, 
\end{equation}
where $A'(u)\ge0 \; (u\in \mathbb{R})$, 
by Osher-Ralston \cite{osh-ral} (see Remark 1.10, see also \cite{cra-tar}). 
In particular, for the Burgers case, 
that is, $f(u)=u^2/2\; (u\in \mathbb{R})$ and $A(u)=u \; (u\in \mathbb{R})$, 
under the condition $u_{-}<u_{+}$, 
Hattori-Nishihara \cite{hattori-nishihara} obtained 
the precise poinwise and time-decay estimates 
of the difference $|u-u^r|$ (see Remark 1.10). 
For the case $p>1$ 
(dilatant type viscosity, 
where the viscosity also originated from 
the Lady{\v{z}}enskaja model, see \cite{du-gu}, \cite{lad}), 
when the convective flux satisfies (1.5) 
and viscous flux 
is even the Ostwald-de Waele type 
($p$-Laplacian type, see \cite{de waele}, \cite{ost}), 
that is, 
\begin{equation}
\sigma(v)= \mu \, 
    \left| \, v \, \right|^{p-1} v, 
\end{equation}
Matsumura-Nishihara \cite{matsu-nishi2} 
proved that if $u_- = u_+ =\tilde{u}$, 
then the solution globally tends toward the constant state $\tilde{u}$, 
and if $u_- < u_+$, then toward the single rarefaction wave (cf. \cite{yoshida8}). 
Yoshida \cite{yoshida3} also obtained 
the precise time-decay estimates of the solution
toward the constant state and the single rarefaction wave. 
For $p>3/7$, Matsumura-Yoshida \cite{matsumura-yoshida} very recently 
showed that if $u_- = u_+ =\tilde{u}$, 
then the solution globally tends toward the constant state $\tilde{u}$, 
and if $u_- < u_+$, then toward the rarefaction wave (see Remark 5.6). 
For $p \ge 1$,
it is further considered a case where the flux function $f$ is smooth
and convex  
on the whole $\mathbb{R}$ except 
a finite interval $I := (a,\, b) \subset \mathbb{R}$, and 
linearly degenerate on $I$, that is, 
\begin{equation}
\left\{
\begin{array}{ll}
  f''(u) >0 & \; \bigl( \, u \in (-\infty ,\, a\, ]\cup [\, b,\, \infty ) \, \bigr),\\[5pt]
  f''(u) =0 & \; \bigl( \, u \in (a,\, b) \, \bigr).
\end{array}\right.
\end{equation} 
Under the conditions $p = 1$, $u_{-}<u_{+}$, (1.8) and (1.9), 
it has been proved by Matsumura-Yoshida \cite{matsumura-yoshida} 
that the unique global in time solution 
to the Cauchy problem (1.1) globally tends toward 
the multiwave pattern of the combination of the rarefaction waves 
and the viscous contact wave 
as time goes to infinity, where the viscous contact wave is given by 
\begin{align}
\begin{aligned}
&U\left(\frac{x-\tilde{\lambda} \, t}{\sqrt{t}}\: ;\: u_- ,\: u_+ \right)\\
& \quad 
 :=u_- +\frac{u_+ - u_-}{\sqrt{\pi}}
   \lint ^{\frac{\mathlarger{x-\tilde{\lambda} t}}{\mathlarger{\sqrt{4\mu t}}}}_{-\infty} 
 \mathrm{e}^{-\xi^2}\, \mathrm{d}\xi 
\quad \, \left( \, \tilde{\lambda}:= \frac{f(b)-f(a)}{b-a} \, \right), 
\end{aligned}
\end{align}
which is constructed by the linear heat kernel 
(in detail, see Section 10, see also \cite{yoshida7}). 
Also the precise time-decay estimates of above stability 
was obtained by Yoshida \cite{yoshida1}. 
Yoshida further proved 
under the conditions $p > 1$, $u_{-}<u_{+}$, (1.8) and (1.9) 
that the global stability of the multiwave pattern 
to the Cauchy problem (1.1) (see \cite{yoshida3}, \cite{yoshida4}), 
and obtained its precise time-decay estimates (see \cite{yoshida5}). 
The multiwave pattern is the combination of the rarefaction waves 
and the viscous contact wave which is given by 
\begin{align}
U\left(\, \frac{x-\tilde{\lambda} \, t}
          {t^{\frac{1}{p+1}}}\: ;\: u_- ,\: u_+ \right) 
:= u_- + \displaystyle{
        \lint^{\frac{\mathlarger{x-\tilde{\lambda} t}}
        { \mathlarger{t}^{{\scriptscriptstyle \frac{1}{p+1}}} }}_{-\infty}
             \Bigl(  \, 
                     \max \big\{ \,  A-B \, \xi^2, \, 0 \, \big\} 
                     \, \Bigr)^{\frac{1}{p-1}}
                \, \mathrm{d}\xi }, 
\end{align}
where 
$$
\int^{\infty}_{-\infty}
                     \Bigl(  \, 
                     \max \big\{ \,  A-B \, \xi^2, \, 0 \, \big\} 
                     \, \Bigr)^{\frac{1}{p-1}}
                     \, \mathrm{d}\xi=| \, u_+ - u_- \, |, 
$$
which is constructed by the Barenblatt-Kompaneec-Zel'dovi{\v{c}} solution 
(see also 
\cite{barenblatt}, \cite{carillo-toscani}, \cite{huang-pan-wang}, \cite{kamin}, 
\cite{pattle}, \cite{vaz1}, \cite{vaz2}, \cite{zel-kom}) 
of the porous medium equation. 
On the other hand, 
under the Rankine-Hugoniot condition 
\begin{equation}
-s\, (\, u_{+} - u_{-}) + f(u_{+}) - f(u_{-}) = 0, 
\end{equation}
and Ole{\u\i}nik's shock condition 
\begin{eqnarray}
- s\, (\, u - u_{\pm}\, ) + f(u) - f(u_{\pm}) 
\left\{\begin{array}{ll}
< 0 
\quad \big( \, u \in (u_{+}, \: u_{-}) \, \big), \\[5pt]
> 0 
\quad \big( \, u \in (u_{-}, \: u_{+}) \, \big),
 \end{array}
 \right.\,
\end{eqnarray}
the local asymptotic stability 
of single viscous shock wave is proved for $p=1$ by 
Matsumura-Nishihara \cite{matsu-nishi3}
(and moreover, they obtained its precise time-decay estimates), 
and recently for 
any $p>0$, more generally, for the case where smooth $\sigma$ satisfies 
\begin{equation}
\sigma(0) = 0,\qquad
\sigma'(v) > 0 \quad ( v \in \mathbb{R} ),\qquad
\displaystyle{\lim_{v\to \pm \infty}} \sigma(v) =\pm \infty, 
\end{equation}
by Yoshida \cite{yoshida6}. 
Also in the case $p=1$, 
under (1.13) and Lax' shock condition 
\begin{equation}
f'(u_{+}) < s < f'(u_{-}) \quad (u_{-}>u_{+}), 
\end{equation}
Freist{\"{u}}hler-Serre \cite{fre-ser} obtained 
the $L^1$-global asymptotic stability of single viscous shock wave, 
and its pointwise decay estimates were also obtained by Deng-Wang \cite{den-wan}. 
Furthermore, for $p<0$ (in particular, for $p=-1$), 
Kurganov-Levy-Rosenau \cite{kur-lev-ros} 
investigated 
the existence, uniqueness and smoothness of a solution 
to the conservation law with nonmonotonic viscous flux, 
that is, 
\begin{equation}
\displaystyle{
\sigma (\partial_{x}u) 
= \pm \, \mu \, 
  \frac{\partial_xu}
  {1 + | \, \partial_xu \, |^2}. 
}
\end{equation}

\medskip

The proofs of Theorems 1.1 and 1.5 are given 
by the technical energy method. 
Also the proofs of Theorems 1.2-1.4 and 1.6-1.8 
are given by the technical time-weighted energy method 
(see Hashimoto-Ueda-Kawashima \cite{hashimoto-kawashima-ueda}) 
and the bootstrap argument 
by Yoshida \cite{yoshida2}, \cite{yoshida4}. 
Because the proofs of Theorems 1.1-1.4 and 
the proofs of Theorem 1.2-1.8 with $p \geq 1$ are much 
easier than that of Theorems 1.2-1.8 
with $0<p<1$, we only show Theorems 1.2-1.8 under the assumption $0<p<1$ in 
the present paper. 

\medskip

\noindent
{\it Remark 1.9.}\quad 
The global structure 
(the existence, uniqueness, smoothness and time-decay properties) 
of the solution to (1.1) 
in all above Theorems 1.1-1.8 coincides with that in the case where 
the viscous flux is linear, 
that is $\sigma(\partial_{x}u)= \mu \, \partial_{x}u$ in (1.1). 

\medskip

\noindent
{\it Remark 1.10.}\quad 
Under the condition $0\le u_{-}<u_{+}$  or $0<u_{-}<u_{+}$, 
Harabetian \cite{har} obtained that the solution to (1.8) 
has the following time-decay properties: 
for $f''(u)>0\; (u\in \mathbb{R})$ and $A'(u)\ge0 \; (u\in \mathbb{R})$, 
\begin{align}
\begin{aligned}
&\Big\| \, 
u(t) - u^r\left(\, \frac{\cdot }{t}\: ;\: u_-, \, u_+ \, \right) 
\, \Big\|_{L^q}\\
&\le C \, 
(1+t)^{-\frac{1}{2}\left( 1-\frac{1}{q} \right)} \,
\big( \, \ln (1+t) \, \big)^{{\frac{1}{2}\left( 1+\frac{1}{q} \right)}}
\quad \big( \, t\ge0 \, \big), 
\end{aligned}
\end{align}
and for $f(u)=u^2/2\; (u\in \mathbb{R})$ and $A(u)=u \; (u\in \mathbb{R})$ 
(Burgers case), 
\begin{equation}
\Big\| \, 
u(t) - u^r\left(\, \frac{\cdot }{t}\: ;\: u_-, \, u_+ \, \right) 
\, \Big\|_{L^q}
\sim 
(1+t)^{-\frac{1}{2}\left( 1-\frac{1}{q} \right)} 
\quad \big( \, t\ge0 \, \big), 
\end{equation}
for $q \in \ [\, 1, \, \infty\,]$. 
Also in such Burgers case, under the conditions 
\begin{equation}
u_{-}<u_{+}, \quad 
u_{0}-u_{0}^{\mathrm{R}} \in L^1, \quad 
\int_{-\infty}^{\infty}
\big( \, u_{0}(x)-u_{0}^{\mathrm{R}}(x) \, \big) 
\, \mathrm{d}x=0, 
\end{equation}
the pointwise and time-decay estimates 
obtained by Hattori-Nishihara \cite{hattori-nishihara} 
are 
\begin{align}
\begin{aligned}
&\left| \, 
 u(t,x)-u^r \left( \, \frac{x}{t}\: ;\:  u_- ,\:  u_+ \,\right) 
 \, \right|\\
&\le C \, 
\left\{
\begin{array}{ll}
  &\big( \, |\, x-u_{-}\, t \, |^{-1} + |\, x-u_{+}\, t \, |^{-1} \, \big) \\[7pt]
& \quad 
  \big(\, 
    u_{-}\, t + K \, \sqrt{\, 2 \, t} \le x
    \le u_{+}\, t - K \, \sqrt{\, 2 \, t} \, \bigr),\\[7pt]
  &t^{-\frac{1}{2}} \quad 
  \big(\, 
    u_{+}\, t - K \, \sqrt{\, 2 \, t} \le x
    \le u_{+}\, t + K \, \sqrt{\, 2 \, t} \, \bigr),\\[7pt]
  &t^{-\frac{1}{2}} \quad 
  \big(\, 
    u_{-}\, t - K \, \sqrt{\, 2 \, t} \le x
    \le u_{-}\, t + K \, \sqrt{\, 2 \, t} \, \bigr),\\[7pt]
  &\min \bigg\{ \, t^{-\frac{1}{2}}, \, 
  t^{-\frac{1}{2}} \, \mathrm{e}^{-(1-\theta)\frac{|x-u_{+}t |^2}{2t}}
  + \displaystyle{\sup_{2 |\eta| \ge \theta |x-u_{+}t |}} \, 
    |\, u_{0}(\eta)-u_{+} \, |
  \, \bigg\}\\[7pt]
& \quad 
  \bigl(\, x \geq u_{+}\, t + K \, \sqrt{\, 2 \, t} \, \bigr),\\[7pt]
  &\min \bigg\{ \, t^{-\frac{1}{2}}, \, 
  t^{-\frac{1}{2}} \, \mathrm{e}^{-(1-\theta)\frac{|x-u_{-}t |^2}{2t}}
  + \displaystyle{\sup_{2 |\eta| \ge \theta |x-u_{-}t |}} \, 
    |\, u_{0}(\eta)-u_{-} \, |
  \, \bigg\}\\[7pt]
& \quad 
  \bigl(\, x \leq u_{-}\, t - K \, \sqrt{\, 2 \, t} \, \bigr)
\end{array}
\right.\\[7pt]
& \big( \, t>0, \: x \in \mathbb{R} \, \big), 
\end{aligned}
\end{align}
where $K\gg1$ and $\theta \in (0, \, 1)$, and 
\begin{equation}
\Big\| \, 
u(t) - u^r\left(\, \frac{\cdot }{t}\: ;\: u_-, \, u_+ \, \right) 
\, \Big\|_{L^q}
\le C \, 
t^{-\frac{1}{2}\left( 1-\frac{1}{q} \right)} 
\quad ( t>0 ), 
\end{equation}
for $q \in \ [\, 1, \, \infty\,]$, respectively. 
From these decay results, 
we note that the time-decay rates in (1.18), (1.19) and (1.22) are 
almost the same as the rates in both Theorems 1.3 and 1.7.

\medskip

This paper is organized as follows. 
In Section 2, we prepare the basic
properties of the rarefaction wave. 
In Section 3, we reformulate the problem 
in terms of the deviation from 
the asymptotic state. 
Also,
in order to show the global existence 
and asymptotic behavior of solution to the reformulated problem, 
we show the strategy how the local existence and
the {\it a priori} estimates are combined.
In Sections 4 and 5,
we give the proof of the {\it a priori} estimates
step by step 
by using the technical energy method.
In Sections 6-9, 
we further give the proof of time-decay estimates 
step by step 
by using the technical time-weighted energy method
with the help of all previous sections. 
Finally in Section 10, we considere the global structure of 
the conservation law without the convective term, 
and compare the structure with that in Section 1 
and some known properties on the viscous contact wave. 

\medskip

{\bf Some Notation.}\quad 
We denote by $C$ generic positive constants unless 
they need to be distinguished. 
In particular, use 
$C_{\alpha,\beta,\cdots }$ 
when we emphasize the dependency on $\alpha,\: \beta,\: \cdots $.
For function spaces, 
${L}^p = {L}^p(\mathbb{R})$ and ${H}^k = {H}^k(\mathbb{R})$ 
denote the usual Lebesgue space and 
$k$-th order Sobolev space on the whole space $\mathbb{R}$ 
with norms $||\cdot||_{{L}^p}$ and $||\cdot||_{{H}^k}$, 
respectively. 

\bigskip 

\noindent
\section{Preliminaries} 
In this section, 
we 
prepare a couple of lemmas concerning with 
the basic properties of 
the rarefaction wave.
Since the rarefaction wave $u^r$ is not smooth enough, 
we need some smooth approximated one. 
We start with the rarefaction wave solution $w^r$ 
to the Riemann problem 
for the non-viscous Burgers equation:
\begin{equation}
\label{riemann-burgers}
  \left\{\begin{array}{l}
  \partial _t w + 
  \displaystyle{ \partial _x \Big( \, \frac{1}{2} \, w^2 \Big) } = 0 
  \quad \quad \qquad \qquad \big( \, t > 0,\: x\in \mathbb{R} \, \big),
\\[7pt]
  w(0, x) = w_0 ^{\rm{R}} ( \, x\: ;\: w_- ,\: w_+\,)
:= \left\{\begin{array}{ll}
w_+ &\quad (x>0),\\[5pt]
w_- &\quad (x<0),
\end{array}
\right.
\end{array}
  \right.
\end{equation}
where $w_\pm \in \mathbb{R}$ are 
the prescribed far field states
satisfying $w_-<w_+$. 
The unique global weak solution 
$w=w^r\left(\,{x}/{t}\: ;\: w_-,\: w_+\,\right)$ 
of (\ref{riemann-burgers}) is explicitly given by 
\begin{equation}
\label{rarefaction-burgers}
w^r \Big(\, \frac{x}{t}\: ;\: w_-,\: w_+ \, \Big) 
= 
  \left\{\begin{array}{ll}
  w_{-} & \bigl(\, x \leq w_{-} \, t \, \bigr),\\[5pt]
  \displaystyle{ \frac{x}{t} } & \bigl(\, w_{-}\, t \leq x \leq w_{+}\, t \, \bigr),\\[5pt]
  w_+ & \bigl(\, x\geq w_{+}\, t \, \bigr).
  \end{array}\right.
\end{equation} 
Next, under the condition 
$f\in C^4(\mathbb{R})$ with $f''(u)>0\ (u\in \mathbb{R})$ and $u_-<u_+$, 
the rarefaction wave solution 
$u=u^r\left( {x}/{t}\: ;\: u_-,\: u_+ \, \right)$ 
of the Riemann problem (1.6) 
for hyperbolic conservation law 
is exactly given by 
\begin{equation}
u^r\left( \, \frac{x}{t} \: ; \:  u_-,\: u_+ \, \right) 
= (\lambda)^{-1}\Big(\, 
  w^r\left( \, \frac{x}{t} \: ; \:  \lambda_-,\: \lambda_+ \, \right)
  \,\Big)
\end{equation}
which is nothing but (1.7), 
where $\lambda(u):=f'(u)$ and $\lambda_\pm := \lambda(u_\pm) = f'(u_\pm)$. 
We define a smooth approximation of $w^r( {x}/{t}\: ;\: w_-,\: w_+ \,)$ 
by the unique classical solution 
$$
w=w(\, t, \, x\: ;\: q,\: w_-,\: w_+ \,)
\in C^{\infty }\big( \, [\, 0,\: \infty \,)\times \mathbb{R} \,\big)
$$
to the Cauchy problem for the following 
non-viscous Burgers equation
\begin{eqnarray}
\label{smoothappm}
\left\{\begin{array}{l}
 \partial _t w 
 + \displaystyle{ \partial _x \Big( \, \frac{1}{2} \, w^2 \, \Big) } =0 
 \, \, \; \; \quad \qquad \qquad \qquad \qquad  \: \,
 \big(\, t>0,\: x\in \mathbb{R} \, \big),\\[7pt]
 w(0, x) 
 = w_0(x) \\[7pt]
 \quad \qquad 
 := \displaystyle{ \frac{w_-+w_+}{2} 
    + \frac{w_+-w_-}{2}\,K_{q} \, \int_{0}^{x} \frac{\mathrm{d}y}{(1+y^2)^q} }
 \quad (x\in \mathbb{R}), 
\end{array}
\right.
\end{eqnarray}   
where $K_{q}$ is a positive constant such that 
$$
K_{q} \, \int_{0}^{\infty} \frac{\mathrm{d}y}{(1+y^2)^q} =1 
\quad \bigg( \, q>\frac{1}{2} \, \bigg). 
$$
By applying the method of characteristics, 
we get the following formula 
\begin{eqnarray}
 \left\{\begin{array} {l}
 w(t, x)=w_0\bigl( \, x_0(t, x) \, \bigr)\\[7pt]
 \: \quad \qquad =
 \displaystyle{ \frac{\lambda_-+\lambda_+}{2} } 
+ \displaystyle{ \frac{\lambda_+-\lambda_-}{2}
\,K_{q} \, \int_{0}^{x_0(t, x) } \frac{\mathrm{d}y}{(1+y^2)^q}  } ,\\[7pt]
 x=x_0(t, x)+w_0\bigl( \, x_0(t, x) \, \bigr)\,t.
 \end{array}
  \right.\,
\end{eqnarray}
By making use of (2.5) similarly as in \cite{matsu-nishi1}, 
we obtain the properties of 
the smooth approximation $w=w(\,t, \, x \: ; \: q, \: w_-, \: w_+\,)$ 
in the next lemma.

\medskip

\noindent
{\bf Lemma 2.1.}\quad{\it
Assume $q>1/2$ and $w_-<w_+$. 
Then the classical solution $w=w(\,t, \, x \: ; \: q, \: w_-, \: w_+\,)$
given by {\rm(2.5)} 
satisfies the following properties. 

\noindent
{\rm (1)}\ \ $w_- < w(t, x) < w_+$ and\ \ $\partial_xw(t, x) > 0$  
\quad  $\big( \, t>0, \: x\in \mathbb{R} \, \big)$.

\smallskip

\noindent
{\rm (2)}\ For any $r \in [\, 1, \, \infty \,]$, there exists a positive 
constant $C_{w_{\pm}, q, r}$ such that
             \begin{eqnarray*}
                 \begin{array}{l}
                    \| \, \partial_x w(t) \, \|_{L^r} \leq 
                    C_{w_{\pm},q,r} \, (1+t)^{-1+\frac{1}{r}} 
                    \quad \bigl( \, t\ge 0 \, \bigr),\\[5pt]
                    \| \, \partial_t w(t) \, \|_{L^r} \leq 
                    C_{w_{\pm},q,r} \, (1+t)^{-1+\frac{1}{r}} 
                    \quad \bigl( \, t\ge 0 \, \bigr),\\[5pt]
                    \| \, \partial_x^2 w(t) \, \|_{L^r} \leq 
                    C_{w_{\pm},q,r} \, (1+t)^{-1-\frac{1}{2q} \, 
                    \left( 1 - \frac{1}{r} \right)} 
                    \quad \bigl( \, t\ge 0 \, \bigr),\\[5pt]
                    \| \, \partial_x^3 w(t) \, \|_{L^r} \leq 
                    C_{w_{\pm},q,r} \, (1+t)^{-1-\frac{1}{2q} \, 
                    \left( 2 - \frac{1}{r} \right)} 
                    \quad \bigl( \, t\ge 0 \, \bigr).
                    \end{array}       
              \end{eqnarray*}
              
\smallskip

\noindent
{\rm (3)}\; $\displaystyle{\lim_{t\to \infty} 
\sup_{x\in \mathbb{R}}
\left| \,w(t, x)- w^r \left( \frac{x}{t} \right) \, \right| = 0}.$

\smallskip

\noindent
{\rm (4)}\ 
There exists a positive 
constant $C_{w_{\pm}, q}$ such that
             \begin{align*}
             \begin{aligned}
                    &| \, w(t, x)-w_{-} \, | \\
                    &\leq 
                    C_{w_{\pm},q} \, 
                    \min \Big\{ \, 
                    \big( \, 1+|\,x - w_{-}\,t\,| \, \big)^{-2q+1}, \, 
                    ( 1+t )^{-\frac{2q-1}{2q}}
                    \, \Big\}
                    \quad \big( \, t\ge0, \: x\le w_{-}\,t \, \big),\\[5pt]
                    &| \, \partial_x w(t, x) \, | 
                    \sim | \, \partial_t w(t, x) \, | \\
                    & \qquad \qquad \; \;
                    \leq 
                    C_{w_{\pm},q} \, 
                    \big( \, 1+|\,x - w_{-}\,t\,|^{2q} +t \, \big)^{-1} 
                    \quad \big( \, t\ge0, \: x\le w_{-}\,t \, \big).
                    \end{aligned}       
              \end{align*}
\smallskip

\noindent
{\rm (5)}\ 
There exists a positive 
constant $C_{w_{\pm}, q}$ such that
             \begin{align*}
             \begin{aligned}
                    &| \, w(t, x)-w_{+} \, | \\
                    &\leq 
                    C_{w_{\pm},q} \, 
                    \min \Big\{ \, 
                    \big( \, 1+|\,x - w_{+}\,t\,| \, \big)^{-2q+1}, \, 
                    ( 1+t )^{-\frac{2q-1}{2q}}
                    \, \Big\}
                    \quad \big( \, t\ge0, \: x\ge w_{+}\,t \, \big),\\[5pt]
                    &| \, \partial_x w(t, x) \, | 
                    \sim | \, \partial_t w(t, x) \, | \\
                    & \qquad \qquad \; \;
                    \leq 
                    C_{w_{\pm},q} \, 
                    \big( \, 1+|\,x - w_{+}\,t\,|^{2q} +t \, \big)^{-1}
                    \quad \big( \, t\ge0, \: x\ge w_{+}\,t \, \big).
                    \end{aligned}       
              \end{align*}
\smallskip

\noindent
{\rm (6)}\ 
There exists a positive 
constant $C_{w_{\pm}, q}$ such that
             \begin{align*}
             \begin{aligned}
                    &\displaystyle{
                    \left| \, w(t, x)-\frac{x}{1+t} \, \right|}
                    \leq 
                    C_{w_{\pm},q} \, 
                    ( 1+t )^{-\frac{2q-1}{2q}} 
                    \quad 
                    \big( \, t\ge0, 
                    \: w_{-}\,t \le x \le w_{+}\,t \, \big),\\[5pt]
                    &\displaystyle{
                    \left| \, \partial_x w(t, x)-\frac{1}{1+t} \, \right|}
                    \sim 
                    \displaystyle{
                    \left| \, \partial_t w(t, x)+\frac{x}{(1+t)^2} \, \right|}\\
                    & \qquad \qquad \qquad \qquad \; 
                    \leq 
                    C_{w_{\pm},q} \, 
                    ( 1+t )^{-1} 
                    \quad 
                    \big( \, t\ge0, 
                    \: w_{-}\,t \le x \le w_{+}\,t \, \big).
                    \end{aligned}       
              \end{align*}
}

\medskip

{\bf Proof of Lemma 2.1.}
The proofs of (1)-(3) in Lemma 2.1 are well-known and given 
in \cite{hashimoto-matsumura}, \cite{hattori-nishihara}, \cite{liu-matsumura-nishihara}, \cite{matsu-nishi1}, 
\cite{matsumura-yoshida}, \cite{yoshida1}, and so on, 
instead of the decay estimate $\| \, \partial_t w \, \|_{L^r}$ in (2). 
However, this estimate is immediately given by 
$| \, \partial_t w(t,x) \, | \sim | \, \partial_x w(t,x) \, |$. 
We also note that the proofs of (5)-(6) are similarly given as that of (5). 
Therefore, we only give the proof of (5). 

First, direct calculation shows from (2.5) that
\begin{align}
\begin{aligned}
x_0(t, w_{+} \, t)
&=\Big( \, w_{+}-w_0\bigl( \, x_0(t, w_{+} \, t) \, \bigr) \, \Big) \, t\\
&=\frac{w_{+}-w_{-}}{2} \,K_{q} \, t 
  \, \int_{0}^{x_0(t, x) } \frac{\mathrm{d}y}{(1+y^2)^q}\\
& \sim \bigl( \, x_0(t, w_{+} \, t) \, \bigr)^{-2q+1} \, t
\quad \big( \, t\ge0, \: x\ge w_{+}\,t \, \big).
\end{aligned}
\end{align}
Therefore, 
\begin{equation}
x_0(t, w_{+} \, t)
\sim t^{\frac{1}{2q}}
\quad \big( \, t\ge0, \: x\ge w_{+}\,t \, \big).
\end{equation}
We also note from (2.5) that 
\begin{equation}
w_0\bigl( \, x_0(t, w_{+} \, t) \, \bigr)
\le w_0\bigl( \, x_0(t, x) \, \bigr)=w(t,x)
\quad \big( \, t\ge0, \: x\ge w_{+}\,t \, \big)
\end{equation}
since 
\begin{equation}
w_0'(x_0)
= \frac{w_{+}-w_{-}}{2} \,
  \frac{K_{q}}{(1+x_0^2)^q}>0
\quad \big( \, t\ge0, \: x\ge w_{+}\,t \, \big).
\end{equation}
Then, we can estimate by using (2.7) and (2.8) as 
\begin{align}
\begin{aligned}
| \, w(t,x)-w_{+} \, |
&\le w_{+}-w_0\bigl( \, x_0(t, w_{+} \, t) \, \bigr)\\
&=\frac{x_0(t, w_{+} \, t)}{t}\\
& \sim t^{-1+\frac{1}{2q}}
\quad \big( \, t>0, \: x\ge w_{+}\,t \, \big).
\end{aligned}
\end{align}
On the other hand, we can also estimate as 
\begin{align}
\begin{aligned}
| \, w(t,x)-w_{+} \, |
&= \frac{w_{+}-w_{-}}{2} \,K_{q} 
  \, \int_{0}^{x_0(t, x) } \frac{\mathrm{d}y}{(1+y^2)^q}\\
& \sim \big( \, 1+x_0(t, x) \, \big)^{-2q+1}\\
& \le \big( \, 1+|\,x - w_{+}\,t\,| \, \big)^{-2q+1}
\quad \big( \, t\ge0, \: x\ge w_{+}\,t \, \big).
\end{aligned}
\end{align}
From (2.10) and (2.11), we have the pointwise estimate of $| \, w(t,x)-w_{+} \, |$. 

We also have the pointwise estimate of $| \, \partial_{x}w(t,x) \, |$ from 
\begin{align*}
\begin{aligned}
| \, \partial_{x}w(t,x) \, |
&= \frac{w_0'\bigl( \, x_0(t, x) \, \bigr)}{1+w_0'\bigl( \, x_0(t, x) \, \bigr) \, t}\\
& \sim \frac{1}{1+\bigl( \, x_0(t, x) \, \bigr)^{2q} + t}
\quad \big( \, t\ge0, \: x\ge w_{+}\,t \, \big).
\end{aligned}
\end{align*}
Therefore, we have (5). 

Thus, the proof of Lemma 2.1 is completed. 

%
%

\bigskip

By using $w$, 
we now define the approximation for 
the rarefaction wave $u^r\left( \,{x}/{t}\: ;\: u_-,\: u_+ \, \right)$ by 
\begin{equation}
U^r(\, t, \, x\: ; \: q, \: u_-,\: u_+ \,) 
= (\lambda)^{-1} 
\bigl( \,w(\, t, \, x\: ;\: q, \: \lambda_-,\: \lambda_+\,) \, \bigr), 
\end{equation}
Noting the assumption of the convective flux $f$, 
we have the next lemma (see also Remark 3.9). 

\medskip

\noindent
{\bf Lemma 2.2.}\quad{\it 
Assume $q>1/2$, $u_-<u_+$ and 
$f\in C^4(\mathbb{R})$ with $f''(u)>0 \ (u \in \mathbb{R})$.  
Then the approximation $U^r_{1}$ defined by {\rm (2.6)} 
satisfies the following properties. 

\noindent
{\rm (1)}\ $U^r(t, x)$ defined by {\rm (2.6)} is 
the unique $C^3$-global solution to
the Cauchy problem
\begin{equation*}
\left\{
\begin{array}{l} 
\partial _t U^r +\partial _x \bigl( \, f(U^r ) \, \bigr) = 0 \quad
\bigl(\, t>0, \: x\in \mathbb{R} \, \bigr),\\[7pt]
U^r(0,x) 
= \displaystyle{ (\lambda)^{-1} \left( \, \frac{\lambda_- + \lambda_+}{2} 
+ \frac{\lambda_+ - \lambda_-}{2} 
  \,K_{q} \, \int_{0}^{x} \frac{\mathrm{d}y}{(1+y^2)^q} \, \right) } 
\quad( x\in \mathbb{R}),\\[10pt]
\displaystyle{\lim_{x\to \pm \infty}} U^r(t, x) =u_{\pm} \quad
\bigl(\, t\ge 0 \, \bigr).
\end{array}
\right.\,     
\end{equation*}
{\rm (2)}\ \ $u_- < U^r(t, x) < u_+$ and\ \ $\partial_xU^r(t, x) > 0$  
\quad  $\bigl(\, t>0, \: x\in \mathbb{R} \, \bigr)$.

\smallskip

\noindent
{\rm (3)}\ For any $r \in [\,1,\, \infty \,]$, there exists a positive 
constant $C_{q,r}$ such that
             \begin{eqnarray*}
                 \begin{array}{l}
                    \| \, \partial_x U^r(t) \, \|_{L^r} 
                    \leq C_{\lambda_{\pm}, q,r}\, (1+t)^{-1+\frac{1}{r}} 
                    \quad \bigl(\, t\ge 0 \, \bigr),\\[5pt]
                    \| \, \partial_t U^r(t) \, \|_{L^r} 
                    \leq C_{\lambda_{\pm}, q,r}\, (1+t)^{-1+\frac{1}{r}} 
                    \quad \bigl(\, t\ge 0 \, \bigr),\\[5pt]
                    \| \, \partial_x^2 U^r(t) \, \|_{L^r} 
                    \leq C_{\lambda_{\pm}, q,r}\, (1+t)^{-1-\frac{1}{2q} \, 
                    \left( 1 - \frac{1}{r} \right)}
                    \quad \bigl(\, t\ge 0 \, \bigr),\\[5pt]
                    \| \, \partial_x^3 U^r(t) \, \|_{L^r} 
                    \leq C_{\lambda_{\pm}, q,r}\, (1+t)^{-1-\frac{1}{2q} \, 
                    \left( 2 - \frac{1}{r} \right)}
                    \quad \bigl(\, t\ge 0 \, \bigr).
                    \end{array}       
              \end{eqnarray*}

\smallskip

\noindent
{\rm (4)}\; $\displaystyle{\lim_{t\to \infty} 
\sup_{x\in \mathbb{R}}
\left| \,U^r(t, x)- u^r \left( \frac{x}{t} \right) \, \right| = 0}.$
              
\smallskip

\noindent
{\rm (5)}\ 
There exists a positive 
constant $C_{\lambda_{\pm}, q}$ such that
\begin{align*}
             \begin{aligned}
                    &| \, U^r(t, x)-u_{-} \, |  \\
                    &\leq 
                    C_{\lambda_{\pm},q} \, 
                    \min \Big\{ \, 
                    \big( \, 1+|\,x - \lambda_{-}\,t\,| \, \big)^{-2q+1}, \, 
                    ( 1+t )^{-\frac{2q-1}{2q}}
                    \, \Big\}
                    \quad \big( \, t\ge0, \: x\le \lambda_{-}\,t \, \big),\\[5pt]
                    &| \, \partial_x U^r(t, x) \, | 
                    \sim | \, \partial_t U^r(t, x) \, | \\
                    & \qquad \qquad \; \; \; \;
                    \leq 
                    C_{\lambda_{\pm},q} \, 
                    \big( \, 1+|\,x - \lambda_{-}\,t\,|^{2q} +t \, \big)^{-1} 
                    \quad \big( \, t\ge0, \: x\le \lambda_{-}\,t \, \big).
                    \end{aligned}       
              \end{align*}
\smallskip

\noindent
{\rm (6)}\ 
There exists a positive 
constant $C_{\lambda_{\pm}, q}$ such that
\begin{align*}
             \begin{aligned}
                    &| \, U^r(t, x)-u_{+} \, | \\
                    &\leq 
                    C_{\lambda_{\pm},q} \, 
                    \min \Big\{ \, 
                    \big( \, 1+|\,x - \lambda_{+}\,t\,| \, \big)^{-2q+1}, \, 
                    ( 1+t )^{-\frac{2q-1}{2q}}
                    \, \Big\}
                    \quad \big( \, t\ge0, \: x\ge \lambda_{+}\,t \, \big),\\[5pt]
                    &| \, \partial_x U^r(t, x) \, | 
                    \sim | \, \partial_t U^r(t, x) \, | \\
                    & \qquad \qquad \; \; \; \;
                    \leq 
                    C_{\lambda_{\pm},q} \, 
                    \big( \, 1+|\,x - \lambda_{+}\,t\,|^{2q} +t \, \big)^{-1} 
                    \quad \big( \, t\ge0, \: x\ge \lambda_{+}\,t \, \big).
                    \end{aligned}       
              \end{align*}
\smallskip

\noindent
{\rm (7)}\ 
There exists a positive 
constant $C_{\lambda_{\pm}, q}$ such that
            \begin{align*}
             \begin{aligned}
                    &\displaystyle{
                    \left| \, U^r(t, x)
                    -\lambda^{-1} \left( \frac{x}{1+t} \right) \, \right|} 
                    \leq 
                    C_{\lambda_{\pm},q} \, 
                    ( 1+t )^{-\frac{2q-1}{2q}}
                    \quad 
                    \big( \, t\ge0, 
                    \: \lambda_{-}\,t \le x \le \lambda_{+}\,t \, \big),\\[5pt]
                    &\displaystyle{
                    \left| \, 
                    \partial_x U^r(t, x)
                    -\partial_x \left( \, \lambda^{-1} 
                    \left( \frac{x}{1+t} \right) \, \right) 
                    \, \right|}\\
                    &\sim \displaystyle{
                    \left| \, 
                    \partial_t U^r(t, x)
                    -\partial_t \left( \, \lambda^{-1} 
                    \left( \frac{x}{1+t} \right) \, \right)\, \right|}
                    \leq 
                    C_{\lambda_{\pm},q} \, 
                    ( 1+t )^{-1} 
                    \quad \big( \, t\ge0, 
                    \: \lambda_{-}\,t \le x \le \lambda_{+}\,t \, \big).
                    \end{aligned}       
              \end{align*}

\smallskip

 \noindent
 {\rm (8)}\ 
 There exists a positive 
 constant $C_{\lambda_{\pm}, q}$ such that 
\begin{align*}
\begin{aligned}
&\left\|
  \,U^r(t,\cdot \: ) - u^r\left( \frac{\cdot }{1+t}\right) \, 
 \right\|_{L^r} \\
 &\leq 
 C_{\lambda_{\pm}, q}\, 
 \left| \, \frac{2}{(2\, q-1) \, (1-\theta) \, r -1} \, \right|^{\frac{1}{r}}\, 
 ( 1+t )^{\left(-1+\frac{1}{2q} \right) \theta}\\
 &\quad + 
 C_{\lambda_{\pm}, q}\, 
 |\,\lambda_{+} - \lambda_{-}\,|^{\frac{1}{r}} \, 
 ( 1+t )^{-1+\frac{1}{r}+\frac{1}{2q}}
 \quad \bigl(t \ge 0 \bigr),
\end{aligned}       
\end{align*}
where $\theta \in [\, 0, \, 1 \, ]$ and $r \in [\, 1, \, \infty\, ]$ 
satisfying $(2 \, q-1) \, (1-\theta) \, r>1$. 

\smallskip

 \noindent
 {\rm (9)}\ 
 There exists a positive 
 constant $C_{\lambda_{\pm}, q}$ such that 
\begin{align*}
\begin{aligned}
&\left\|
  \,U^r(t,\cdot \: ) - u^r\left( \frac{\cdot }{1+t}\right) \, 
 \right\|_{L^r} \\ 
 &\leq C_{\lambda_{\pm}, q}\, 
 \max \bigg\{ \, 
 \left| \, \frac{2}{2 \, q -1} \, \right|^{\frac{1}{r}}, \, 
 |\,\lambda_{+} - \lambda_{-}\,|^{\frac{1}{r}} 
 \, \bigg\}\, 
 (1+t)^{-1+\frac{1}{r}+\frac{1}{2q}} 
 \quad \bigl(t \ge 0 \bigr),
\end{aligned}       
\end{align*}
where $2 \, q/(2 \, q -1)\le r \le \infty$. 

\smallskip

 \noindent
 {\rm (10)}\ 
 There exists a positive 
 constant $C_{\lambda_{\pm}, q}$ such that 
\begin{align*}
\begin{aligned}
&\|
  \,\partial_x U^r(t,\cdot \: ) - \partial_x u^r(1+t,\cdot \: ) \, 
 \|_{L^r} 
 \sim 
 \|
  \,\partial_t U^r(t,\cdot \: ) - \partial_t u^r(1+t,\cdot \: ) \, 
 \|_{L^r}
 \\ 
 &\leq C_{\lambda_{\pm}, q}\, 
 \max \bigg\{ \, 
 \left| \, \frac{2}{2 \, q -1} \, \right|^{\frac{1}{r}}, \, 
 |\,\lambda_{+} - \lambda_{-}\,|^{\frac{1}{r}} 
 \, \bigg\}\, 
 (1+t)^{-1+\frac{1}{r}} 
 \quad \bigl(t \ge 0 \bigr),
\end{aligned}       
\end{align*}
where $r \in [\, 1, \, \infty\, ]$. 

\smallskip

 \noindent
 {\rm (11)}\ 
 There exists a positive 
 constant $C_{\lambda_{\pm}, q}$ such that 
\begin{equation*}
\|
  \,\partial_x^2 U^r(t,\cdot \: ) - \partial_x^2 u^r(1+t,\cdot \: ) \, 
 \|_{L^2} 
 \leq C_{\lambda_{\pm}, q}\, 
 (1+t)^{-1-\frac{1}{4q}} 
 \quad \bigl(t \ge 0 \bigr).
\end{equation*}
}

\medskip

{\bf Proof of Lemma 2.2.}
Because the proofs of (1)-(7) in Lemma 2.2 are easily given 
by Lemma 2.1 
and the proof of (10) is similarly given as the proofs of (8) and (9), 
we only give the proofs of (8), (9) and (11). 

We first note by using (5) that 
for any $\theta \in [\, 0, \, 1 \, ]$, 
\begin{align}
\begin{aligned}
                    &| \, U^r(t, x)-u_{-} \, |  \\
                    &\leq 
                    C_{\lambda_{\pm},q} \, 
                    ( 1+t )^{\left(-1+\frac{1}{2q} \right) \theta}\, 
                    \big( \, 1+|\,x - \lambda_{-}\,t\,| \, \big)^{-(2q-1)(1-\theta)}
                    \quad \big( \, t\ge0, \: x\le \lambda_{-}\,t \, \big).
\end{aligned}       
\end{align}
For $1 \le r < \infty$, we then have 
\begin{align}
\begin{aligned}
&\int_{-\infty}^{\lambda_{-}t}| \, U^r(t, x)-u_{-} \, |^r
\, \mathrm{d}x\\
&\leq 
\frac{C_{\lambda_{\pm},q}^r}{(2\, q-1) \, (1-\theta) \, r -1}\, 
                    ( 1+t )^{\left(-1+\frac{1}{2q} \right) r \theta}
                    \quad \big( \, t \geq 0  \: ; \: \theta \in [\, 0, \, 1 \, ]\, \big) 
\end{aligned}       
\end{align}
when $(2 \, q-1) \, (1-\theta) \, r>1$. 
Similarly by using (5) and (6), 
we also have 
\begin{align}
\begin{aligned}
&\int_{\lambda_{+}t}^{\infty}| \, U^r(t, x)-u_{+} \, |^r
\, \mathrm{d}x\\
&\leq  
\frac{C_{\lambda_{\pm},q}^r}{(2\, q-1) \, (1-\theta) \, r -1}\, 
                    ( 1+t )^{\left(-1+\frac{1}{2q} \right) r \theta}
                    \quad \big( \, t \geq 0  \: ; \: \theta \in [\, 0, \, 1 \, ]\, \big) 
\end{aligned}       
\end{align}
when $(2 \, q-1) \, (1-\theta) \, r>1$, and 
\begin{equation}
\int_{\lambda_{-}t}^{\lambda_{+}t}
\left| \, U^r(t, x)-\lambda^{-1} \left( \frac{x}{1+t} \right) \, \right|^r
\, \mathrm{d}x
\leq C_{\lambda_{\pm},q}^r \, 
|\,\lambda_{+} - \lambda_{-}\,|\, 
                    ( 1+t )^{\left(-1+\frac{1}{2q} \right) r +1}
                    \quad \big( \, t \geq 0  \, \big), 
\end{equation}
respectively. 
Therefore, combining (2.14), (2.15) and (2.16), 
we immediately get the $L^r$-estimate in (8) 
with $q>1/2$, $\theta \in [\, 0, \, 1 \, ]$ and $r \in [\, 1, \, \infty)$ 
satisfying $(2 \, q-1) \, (1-\theta) \, r>1$. 
Also taking the limit $r\rightarrow \infty$ to this, we have (8). 

Next, we show (9). 
If we choose 
$$
\theta=1-\frac{2 \, q}{r \, (2\, q -1)}\in [\, 0, \, 1 \, ]
\quad \left( \, 0 < \frac{2 \, q}{2 \, q -1}\le r  \, \right)
$$
in order to satisfy 
$$
\left(-1+\frac{1}{2 \, q} \right) r \, \theta=\left(-1+\frac{1}{2 \, q} \right) r +1, 
$$
then (8) immediately becomes (9).

We finally show (11). 
By using (3) in Lemma 2.2, we have 
\begin{align*}
\begin{aligned}
&\|
 \,\partial_x^2 U^r(t,\cdot \: ) - \partial_x^2 u^r(1+t,\cdot \: ) \, 
 \|_{L^2}^2 \\
& \leq \| \,\partial_x^2 U^r(t ) \, \|_{L^2}^2 
  + \frac{1}{(1+t)^4} \, 
    \int_{\lambda_{-}t}^{\lambda_{+}t} 
    \displaystyle{ 
  \frac{\bigg| \, 
  \lambda''\bigg( \, \lambda^{-1}\left( \displaystyle{\frac{x}{1+t}} \right) \, \bigg)
  \, \bigg|^2}
  {\bigg| \, 
  \lambda'\bigg( \, \lambda^{-1}\left( \displaystyle{\frac{x}{1+t}} \right) \, \bigg)
  \, \bigg|^6} 
  }
  \, \mathrm{d}x\\
 & \leq C_{\lambda_{\pm}, q}\, (1+t)^{-2-\frac{1}{2q}}
   + C_{\lambda_{\pm}, q}\, (1+t)^{-3}
   \quad \big( \, t \geq 0  \, \big).
\end{aligned}       
\end{align*}
Therefore, we have (11).


Thus, the proof of Lemma 2.2 is completed.

\medskip

\noindent
{\it Remark 2.3.}\quad 
Now rewriting $q$ and $r$, in (8) and (9), as $s$ and $q$, respectively, 
we have 
\begin{equation}
\left\|
  \,U^r(t,\cdot \: ) - u^r\left( \frac{\cdot }{1+t} \right) \, 
 \right\|_{L^q} 
\leq C_{\lambda_{\pm}, q, \epsilon}\, 
 (1+t)^{-1+\frac{1}{q}+\epsilon} 
 \quad \bigl(t \ge 0 \bigr),  
\end{equation}
for $q \in [\, 1, \, \infty\, ]$ 
by choosing $s=1/(2 \, \epsilon)\gg1$ $(0<\epsilon \ll1)$, 
and 
\begin{equation}
\left\|
  \,U^r(t,\cdot \: ) - u^r\left( \frac{\cdot }{1+t} \right) \, 
 \right\|_{L^q} 
\leq C_{\lambda_{\pm}, q}\, 
 (1+t)^{-\frac{1}{2}+\frac{1}{q}} 
 \quad \bigl(t \ge 0 \bigr),  
\end{equation}
for $q \in [\, 2, \, \infty\, ]$ 
by choosing 
$$
s=\max_{q\ge2} \, \frac{2 \, q}{3 \, q-2}=1 
$$
in order to satisfy 
$$
-\frac{1}{4}\, \left( 1- \frac{2}{q} \right) 
\le -1 + \frac{1}{q} +\frac{1}{2 \, s} 
$$
(compare 
(2.17) with $L^q$-estimates of $\phi$ in Theorem 3.7, 
and (2.18) with that of $\phi$ in Theorem 3.5).

%
%

%
%

\bigskip 

\noindent
\section{Reformulation of the problem} 
In this section, we reformulate our problem (1.1) 
in terms of the deviation from the asymptotic state. 
To do that, using (2.12), 
and noting Remark 2.3 (in particular, (2.17) and (2.18)), we prepare 
two approximations $U^r_{1}$ and $U^r_{2}$ for 
the rarefaction wave $u^r\left( \,{x}/{t}\: ;\: u_-,\: u_+ \, \right)$, 
that is, 
\begin{equation}
U^r_{1}(\, t, \, x\: ; \: u_-,\: u_+ \,) 
:= (\lambda)^{-1} 
\bigl( \,w(\, t, \, x\: ;\: q=1, \: \lambda_-,\: \lambda_+\,) \, \bigr), 
\end{equation}
\begin{equation}
U^r_{2}(\, t, \, x\: ; \: u_-,\: u_+ \,) 
:= (\lambda)^{-1} 
\bigl( \,w(\, t, \, x\: ;\: q=r, \: \lambda_-,\: \lambda_+\,) \, \bigr) 
\quad (r \gg1).
\end{equation}
Now letting 
\begin{equation}
u(t, x) = U^r_{1}(t, x) + \phi(t, x), 
\end{equation}
we reformulate the problem (1.1) in terms of 
the deviation $\phi $ from $U^r_{1}$ as 
\begin{eqnarray}
 \left\{\begin{array}{ll}
  \partial _t\phi + \partial_x \big( \, f(\phi+U^r_{1}) - f(U^r_{1}) \, \big) \\[5pt]
- \partial_x 
    \Big( \, 
    \sigma ( \, \partial_x \phi + \partial_x U^r_{1} \, ) 
    - \sigma ( \, \partial_x U^r_{1}  \, )  \, 
    \Big)
= 
    \partial_x \Big( \, \sigma ( \, \partial_x U^r_{1} \,)  \, \Big) 
    \quad \big( \, t>0,\: x\in \mathbb{R} \, \big),\\[5pt]
  \phi(0, x) = \phi_0(x) 
  := u_0(x)-U^r_{1}(0, x) \rightarrow 0 \quad (x \rightarrow \pm \infty).
 \end{array}
 \right.\,
\end{eqnarray}
Then we look for 
the unique global in time solution 
$\phi $ which has the asymptotic behavior 
\begin{equation}
\displaystyle{
\sup_{x \in \mathbb{R}}
\big( \, \left|\,  \phi(t, x) \, \right|+
\left|\,  \partial_{x}\phi(t, x) \, \right|
\, \big)
\to  0 \qquad (t\to  \infty)}. 
\end{equation}
Here we note 
that $\phi_0 \in H^2$ by the assumptions on $u_0$, 
and Lemma 2.2. 
Then the corresponding theorem for $\phi$ to Theorem 1.2  we should prove is 
as follows. 

\medskip

\noindent
{\bf Theorem 3.1.} \quad{\it
Assume the far field states satisfy $u_- < u_+$, 
the convective flux $f \in C^4(\mathbb{R})$, {\rm(1.5)}, 
the viscous flux $\sigma \in C^2(\mathbb{R})$, {\rm(1.2)}, {\rm(1.3)}, 
and $p>1/3$. 
Further assume the initial data satisfy
$\phi_0 \in H^2$. 
Then the Cauchy problem {\rm(3.4)} has a 
unique global in time 
solution $\phi$ 
satisfying 
\begin{eqnarray*}
\left\{\begin{array}{ll}
\phi \in C^0\cap L^{\infty}
\bigl(\,[\, 0, \, \infty \, ) \, ; H^2 \,\bigr),\\[5pt]
\partial _x \phi  \in L^2\bigl( \, 0,\, \infty \, ; H^2 \,\bigr),\\[5pt]
\partial _t \phi  \in C^0\bigl(\,[\, 0, \, \infty \, ) \, ; L^2 \,\bigr)
                      \cap L^2\bigl( \, 0,\, \infty \, ; H^1 \,\bigr),
\end{array} 
\right.\,
\end{eqnarray*}
and the asymptotic behavior 
$$
\lim _{t \to \infty}\, \sup_{x\in \mathbb{R}}
\big( \, \left|\,  \phi(t, x) \, \right|+
\left|\,  \partial_{x}\phi(t, x) \, \right|
\, \big) = 0, 
$$
$$
\lim _{t \to \infty} \, 
\big( \, \|\,  \partial_{t}\phi(t) \, \|_{L^2}+
\|\,  \partial_{x}\phi(t) \, \|_{H^1}
\, \big) = 0. 
$$
}

Theorem 3.1 is shown by combining the 
local existence of the solution together with
the {\it a priori} estimates. 
To state the local existence precisely, 
the Cauchy problem at general 
initial time $\tau \ge 0$ with the given initial 
data $\phi_\tau \in H^2$ is formulated as follows. 
\begin{eqnarray}
 \left\{\begin{array}{ll}
  \partial _t\phi + \partial_x \big( \, f(\phi+U^r_{1}) - f(U^r_{1}) \, \big) \\[5pt]
- \partial_x 
    \Big( \, 
    \sigma \big( \, \partial_x \phi + \partial_x U^r_{1} \,\big) 
    - \sigma \big( \, \partial_x U^r_{1}  \,\big)  \, 
    \Big)
    = 
    \partial_x \Big( \, \sigma\big( \, \partial_x U^r_{1}  \,\big)  \, \Big) \quad
\big(\,t>\tau,\: x\in \mathbb{R}\, \big),
    \\[5pt]
  \phi(\tau, x) = \phi_\tau(x) 
  := u_\tau(x)-U^r_{1}(\tau, x)
  \rightarrow 0 \quad (x\rightarrow \pm \infty).
 \end{array}
 \right.\,
\end{eqnarray}

\medskip

\noindent
{\bf Theorem 3.2} (local existence){\bf .}\quad{\it
For any $ M > 0 $, there exists a positive constant 
$t_0=t_0(M)$ not depending on $\tau$ 
such that if $\phi_{\tau} \in H^2$ and 
$\displaystyle{ 
\| \, \phi_{\tau} \, \|_{{H}^{2}} \leq M}$, then
the Cauchy problem {\rm (3.6)} has a unique solution $\phi$ 
on the time interval $[\, \tau, \, \tau+t_{0}(M)\, ]$ satisfying 
\begin{eqnarray*}
\left\{\begin{array}{ll}
\phi \in C^0\big( \, [\, \tau, \, \tau+t_{0}\, ] \, ; H^2 \, \big) \cap
L^2( \, \tau, \, \tau+t_{0} \, ; H^3 \, ),\\[5pt]
\partial _t\phi \in C^0\big( \, [\, \tau, \, \tau+t_{0}\, ] \, ; L^2 \, \big) \cap
L^2( \, \tau, \, \tau+t_{0} \, ; H^1 \, ).
\end{array} 
\right.\,
\end{eqnarray*}
}

\medskip

\noindent
Because the proof of Theorem 3.2 is similar to the one in Yoshida \cite{yoshida7}, 
we omit the details here 
(cf. \cite{kato1}, \cite{kato2}, \cite{lad-sol-ura}, \cite{lions}, \cite{yoshida3}, \cite{yoshida4}). 
The {\it a priori estimates} we establish in Sections 4 and 5 
are the following. 

\medskip

\noindent
{\bf Theorem 3.3} ({\it a priori} estimates){\bf .}\quad{\it 
Under the same assumptions as in Theorem 3.1, 
for any initial data 
$\phi_0 \in H^2$, 
there exists a positive constant $C_{\phi_0}$
such that 
if the Cauchy problem {\rm (3.4)}
has a solution 
$\phi$ 
on a time interval $[\, 0, \, T\, ]$ satisfying 
\begin{eqnarray*}
\left\{\begin{array}{ll}
\phi \in C^0\big( \, [\, 0, \, T\, ] \, ; H^2 \, \big) \cap
L^2( \, 0, \, T \, ; H^3\, ),\\[5pt]
\partial _t\phi \in C^0\big( \, [\, 0, \, T\, ] \, ; L^2 \, \big) \cap
L^2( \, 0, \, T \, ; H^1\, ),
\end{array} 
\right.\,
\end{eqnarray*}
for some constant $T>0$, 
then it holds that 
\begin{align}
\begin{aligned}
&\| \, \phi(t) \, \|_{{H}^{2}}^{2} 
+ \| \, \partial_{t}\phi(t) \, \|_{{L}^{2}}^{2}
+ \int^{t}_{0} 
   \big\| \, 
   (\, \sqrt{\, \partial_x U^r_{1}} \:  \phi \,)(\tau) 
   \, \big\|_{L^2}^{2} \, \mathrm{d}\tau 
\\
& \quad 
+ \int^{t}_{0}
   \big( 
 \left\| \, \partial_{x}\phi(\tau) \, \right\|_{{H}^2}^{2} 
 + \left\| \, \partial_{t}\phi(\tau) \, \right\|_{{H}^1}^{2} 
   \, \big) \, \mathrm{d}\tau 
\leq C_{\phi_{0}}
\qquad \big( \, t \in [\, 0, \, T\, ] \,\big).
\end{aligned}
\end{align}
}

\noindent
Once Theorem 3.3 is established, 
by combining the local existence Theorem 3.2
with  $M=M_0:=\sqrt{ \, C_{\phi_0}}$,
$\tau = n \, t_0(M_0)$, 
and $\phi_\tau=\phi\big( \, n \, t_0(M_0) \, \big)\ (n=0,1,2,\dots)$ together 
with the {\it a priori} estimates with $T=(n+1) \, t_0(M_0)$
inductively,
the
unique solution of (3.5) 
\begin{eqnarray*}
\left\{\begin{array}{ll}
\phi \in C^0\big( \, [\, 0,\, n\, t_0(M_0)\, ] \, ;H^2\,\big)
\cap L^2\big( \,0,\, n\, t_0(M_0)\, ;H^3\,\big),\\[5pt]
\partial _t\phi \in C^0\big( \, [\, 0,\, n\, t_0(M_0)\, ] \, ;L^2\,\big)
\cap L^2\big( \,0,\, n\, t_0(M_0)\, ;H^1\,\big)
\end{array} 
\right.\,
\end{eqnarray*}
for any $n \in \mathbb{N}$
is easily constructed,
that is, the global in time solution 
\begin{eqnarray*}
\left\{\begin{array}{ll}
\phi\in C^0\big( \, [\, 0,\, \infty\,) \, ;H^2\,\big)
\cap L^2_{\rm{loc}}( \, 0,\, \infty \, ; H^3\,),\\[5pt]
\partial _t\phi \in C^0\big( \, [\, 0,\, \infty\,) \, ;L^2\,\big)
\cap L^2( \, 0,\, \infty \, ; H^1\,). 
\end{array} 
\right.\,
\end{eqnarray*}
Then, the {\it a priori} estimates again
assert that
\begin{eqnarray}
\left\{\begin{array}{ll}
\displaystyle{
\sup_{t\ge 0}\, \big( \, 
                \| \, \phi (t)\,  \|_{H^2}
                + \| \, \partial _t\phi (t)\,  \|_{L^2}
                \, \big)}< \infty,\\[15pt]
\displaystyle{
\int _0^\infty \big(\,\| \, \partial_x\phi(t) \, \|_{H^2}^2
+\| \, \partial_t\phi(t) \, \|^2_{H^1} \, \big)
\, \mathrm{d}t }
< \infty,
\end{array} 
\right.\,
\end{eqnarray}
which easily gives
\begin{align}
\begin{aligned}
\int _0^{\infty }\bigg|\,\frac{\mathrm{d}}{\mathrm{d}t} 
\|  \, \partial_x\phi(t)  \, \|_{L^2}^2 \,\bigg| \, \mathrm{d}t
< \infty.
\end{aligned}
\end{align}
Hence, it follows from (3.8) and (3.9) that
\begin{equation}
\|  \, \partial_x\phi(t)  \, \|_{L^2} \to 0\quad (t \to \infty).
\end{equation}
Due to the Sobolev inequality, the
desired asymptotic behavior in Theorem 3.1 is obtained as
\begin{equation}
\sup_{x \in \mathbb{R}}|\,\phi(t,x)\,|
\leq \sqrt{2}\, \| \,\phi(t)\, \|^{\frac{1}{2}}_{L^2} 
\| \,\partial_x\phi(t)\, \| ^{\frac{1}{2}}_{L^2}
\to 0 \quad (t \to \infty),
\end{equation} 
\begin{equation}
\sup_{x \in \mathbb{R}}|\,\partial_{x}\phi(t,x)\,|
\leq \sqrt{2}\, \| \,\partial_{x}\phi(t)\, \|^{\frac{1}{2}}_{L^2} 
\| \,\partial_x^2\phi(t)\, \| ^{\frac{1}{2}}_{L^2}
\to 0 \quad (t \to \infty).
\end{equation} 
Further, we can show (for the proof, see Lemma 5.5)
\begin{align}
\begin{aligned}
\int _0^{\infty }\bigg|\,\frac{\mathrm{d}}{\mathrm{d}t} 
\|  \, \partial_t\phi(t)  \, \|_{L^2}^2 \,\bigg| \, \mathrm{d}t
< \infty.
\end{aligned}
\end{align}
Hence, it also follows from (3.8) and (3.13) that
\begin{equation}
\|  \, \partial_t\phi(t)  \, \|_{L^2} \to 0\quad (t \to \infty), 
\end{equation}
and by using (3.4) that 
\begin{equation}
\|  \, \partial_x^2\phi(t)  \, \|_{L^2} \to 0\quad (t \to \infty). 
\end{equation}
Thus, Theorem 3.1 is shown by combining Theorem 3.2 together with
Theorem 3.3. 

\medskip 

Next, we consider the time-decay estimates. 
In order to state them precisely 
(in partuclar, the derivatives, Theorem 1.8, see Remark 3.9), 
we also let 
\begin{equation}
u(t, x) = U^r_{2}(t, x) + \psi(t, x), 
\end{equation}
we reformulate the problem (1.1) in terms of 
the deviation $\psi $ from $U^r_{2}$ as 
\begin{eqnarray}
 \left\{\begin{array}{ll}
  \partial _t\psi + \partial_x \big( \, f(\psi+U^r_{2}) - f(U^r_{2}) \, \big) \\[5pt]
- \partial_x 
    \Big( \, 
    \sigma ( \, \partial_x \psi + \partial_x U^r_{2} \, ) 
    - \sigma ( \, \partial_x U^r_{2}  \, )  \, 
    \Big)
= 
    \partial_x \Big( \, \sigma ( \, \partial_x U^r_{2}  \,)  \, \Big) 
    \quad \big( \, t>0,\: x\in \mathbb{R} \, \big),\\[5pt]
  \psi(0, \, x) = \psi_0(x) 
  := u_0(x)-U^r_{2}(0, \, x) \rightarrow 0 \quad (x \rightarrow \pm \infty).
 \end{array}
 \right.\,
\end{eqnarray}

\medskip

\noindent
We note that, 
if the above theorems for $\phi$, Theorems 3.1-3.3, are proved, 
then by Lemma 2.2, 
we can similarly prove that 
$\psi$ defined by (3.16) satisfies the same results as in Theorems 3.1-3.3 
(cf. \cite{matsumura-yoshida'}), 
in particular, 

\medskip

\noindent
{\bf Theorem 3.4.} \quad{\it
Assume the far field states satisfy $u_- < u_+$, 
the convective flux $f \in C^4(\mathbb{R})$, {\rm(1.5)}, 
the viscous flux $\sigma \in C^2(\mathbb{R})$, {\rm(1.2)}, {\rm(1.3)}, 
and $p>1/3$. 
Further assume the initial data satisfy
$\psi_0 \in H^2$. 
Then the Cauchy problem {\rm(3.4)} has a 
unique global in time 
solution $\phi$ 
satisfying 
\begin{eqnarray*}
\left\{\begin{array}{ll}
\psi \in C^0\cap L^{\infty}
\bigl(\,[\, 0, \, \infty \, ) \, ; H^2 \,\bigr),\\[5pt]
\partial _x \psi  \in L^2\bigl( \, 0,\, \infty \, ; H^2 \,\bigr),\\[5pt]
\partial _t \psi  \in C^0\bigl(\,[\, 0, \, \infty \, ) \, ; L^2 \,\bigr)
                      \cap L^2\bigl( \, 0,\, \infty \, ; H^1 \,\bigr),
\end{array} 
\right.\,
\end{eqnarray*}
with 
\begin{align*}
\begin{aligned}
&\sup_{t\ge 0}\, \big( \, 
                \| \, \psi (t)\,  \|_{H^2}
                + \| \, \partial _t\psi (t)\,  \|_{L^2}
                \, \big)\\
& \quad 
+
\int _0^\infty \big(\,\| \, \partial_x\psi(t) \, \|_{H^2}^2
+\| \, \partial_t\psi(t) \, \|^2_{H^1} \, \big)
\, \mathrm{d}t 
\leq C_{\psi_{0}},
\end{aligned}
\end{align*}
and the asymptotic behavior 
$$
\lim _{t \to \infty}\, \sup_{x\in \mathbb{R}}
\big( \, \left|\,  \psi(t, x) \, \right|+
\left|\,  \partial_{x}\psi(t, x) \, \right|
\, \big) = 0, 
$$
$$
\lim _{t \to \infty} \, 
\big( \, \|\,  \partial_{t}\psi(t) \, \|_{L^2}+
\|\,  \partial_{x}\psi(t) \, \|_{H^1}
\, \big) = 0. 
$$
}

\medskip

\medskip

Thus, by Lemma 2.2 and Remark 2.3, Theorems 3.1 and 3.4, 
the corresponding time-decay theorems for $\phi$ and $\psi$ to Theorems 1.6, 1.7 and 1.8 
we should prove are the following. 

\medskip

\noindent
{\bf Theorem 3.5.} \quad{\it
Under the same assumptions as in Theorem 3.1, 
the unique global in time 
solution $\phi$ 
of the Cauchy problem {\rm(3.4)} 
has the following time-decay estimates 
\begin{eqnarray*}
\left\{\begin{array}{ll}
\| \, \phi(t) \, \|_{L^q}
\le C_{\phi_{0}} \, 
(1+t)^{-\frac{1}{4}\left( 1-\frac{2}{q} \right)}
\quad \big( \, t\ge0 \, \big),\\[5pt]
\| \, \phi(t) \, \|_{L^{\infty}}
\le C_{\phi_{0}, \epsilon} \, (1+t)^{-\frac{1}{4}+\epsilon}
\quad \big( \, t\ge0 \, \big),
\end{array} 
\right.\,
\end{eqnarray*}
for $q \in [\, 2, \, \infty)$ and any $\epsilon>0$. 
}

\medskip

\noindent
{\bf Theorem 3.6.} \quad{\it
Under the same assumptions as in Theorem 3.1, 
if the initial data further satisfies $\phi_{0}\in L^1$, 
then the unique global in time 
solution $\phi$ 
of the Cauchy problem {\rm(3.4)} 
has the following time-growth and time-decay estimates 
\begin{eqnarray*}
\left\{\begin{array}{ll}
\| \, \phi(t) \, \|_{L^1}
\le C_{\phi_{0}} \, 
\max \big\{ \, 1, \, \ln (1+t) \, \big\}
\quad \big( \, t\ge0 \, \big),\\[5pt]
\| \, \phi(t) \, \|_{L^q}
\le C_{\phi_{0}} \, 
(1+t)^{-\frac{1}{2}\left( 1-\frac{1}{q} \right)} \, 
\max \big\{ \, 1, \, \ln (1+t) \, \big\}
\quad \big( \, t\ge0 \, \big),\\[5pt]
\| \, \phi(t) \, \|_{L^{\infty}}
\le C_{\phi_{0}, \epsilon} \, (1+t)^{-\frac{1}{2}+\epsilon}
\quad \big( \, t\ge0 \, \big),
\end{array} 
\right.\,
\end{eqnarray*}
for $q \in ( 1, \, \infty)$ and any $\epsilon>0$. 
}

\medskip

\noindent
{\bf Theorem 3.7.} \quad{\it
Assume the same assumptions as in Theorem 3.4, 
then $\phi$ in Theorems 3.5 and 3.6 should be 
replaced by $\psi$. 
}
%
%

\medskip

\noindent
{\bf Theorem 3.8.} \quad{\it
Under the same assumptions as in Theorem 3.1, 
the unique global in time 
solution $\phi$ 
of the Cauchy problem {\rm(3.4)} 
has the following time-decay estimates 
for the derivatives 
\begin{eqnarray*}
\left\{\begin{array}{ll}
\| \, \partial_{x}\phi(t) \, \|_{L^{q+1}}
\le C_{\phi_{0}, \epsilon} \, 
(1+t)^{-\frac{2q+1}{2q+2}+\epsilon}
\quad \big( \, t\ge0 \, \big),\\[5pt]
\| \, \partial_{x}\phi(t) \, \|_{L^{\infty}}
\le C_{\phi_{0}, \epsilon} \, 
(1+t)^{-1+\epsilon} 
\quad \big( \, t\ge0 \, \big),\\[5pt]
\| \, \partial_{t}\phi(t) \, \|_{L^{2}}
\le C_{\phi_{0}, \epsilon} \, (1+t)^{-\frac{3}{4}+\epsilon}
\quad \big( \, t\ge0 \, \big),\\[5pt]
\| \, \partial_{x}^2\phi(t) \, \|_{L^{2}}
\le C_{\phi_{0}, \epsilon} \, (1+t)^{-\frac{3}{4}+\epsilon}
\quad \big( \, t\ge0 \, \big),
\end{array} 
\right.\,
\end{eqnarray*}
for $q \in [\, 1, \, \infty)$ and any $\epsilon>0$. 
}
%
%
%
%


\medskip

\noindent
{\it Remark 3.9.}\quad 
We note that the time-decay estimates 
of $\| \, \partial_{x}\phi \, \|_{L^{q+1}}$ with $1 \le q \le \infty$, 
$\| \, \partial_{t}\phi \, \|_{L^{2}}$ 
and $\| \, \partial_{x}^2\phi \, \|_{L^{2}}$ 
in Theorem 3.8 decay faster than 
$\| \, \partial_{x}\psi \, \|_{L^{q+1}}$ with $1 \le q \le \infty$, 
$\| \, \partial_{t}\psi \, \|_{L^{2}}$ 
and $\| \, \partial_{x}^2\psi \, \|_{L^{2}}$, 
respectively.  
This is because we easily know from Lemma 2.2 that 
the $L^r$-decay estimates of $\partial_{x}^2U^r_{1}$ and $\partial_{x}^3U^r_{1}$
are sharper than the estimates of $\partial_{x}^2U^r_{2}$ and $\partial_{x}^3U^r_{2}$, 
respectively, although 
the $L^r$-decay estimate of the difference $| \, U^r_{2}-u^r \, |$ 
is sharper than the estimate of $| \, U^r_{1}-u^r \, |$ 
(we also note that there have not been known the pointwise estimates 
of $| \, U^r-u^r \, |$ with $1/2<q \le 3/2$ so far, 
see \cite{liu-matsumura-nishihara}, 
\cite{matsu-nishi2-0}, \cite{matsu-nishi2}, 
\cite{matsumura-yoshida'}, 
cf. \cite{hashimoto-matsumura}, \cite{hattori-nishihara}, 
\cite{matsu-nishi1}, 
\cite{matsumura-yoshida}, \cite{yoshida1}, and so on). 
In fact, by the similar arguments as in Section 8 
(in particular, we note the fifth term on the right-hand side of (8.13) 
which is led from the nonlinear viscous flux term 
$\sigma' \big( \, \partial_{x}U^r_{2} \, \big) \, \partial_{x}^2U^r_{2}$), 
we can obtain the time-decay estimate for 
$\| \, \partial_x\psi \, \|_{L^2} $ correspond to (8.14) in Lemma 8.4 as 
\begin{equation}
\| \, \partial_x\psi(t) \, \|_{L^2} 
                    \leq C_{\psi_{0}, r}\, 
                    (1+t)^{-\frac{1}{2} \left( 1+\frac{1}{2r} \right)}
                    \quad \bigl(\, t\ge 0 \, \bigr)
\end{equation}
with $r\gg1$, for some $C_{\psi_{0}, r}>0$. 

\bigskip

In the following sections, we first give the proof of 
the {\it a priori} estimates, Theorem 3.3, in Sections 4 and 5.
To do that, we assume 
\begin{eqnarray*}
\left\{\begin{array}{ll}
\phi \in C^0\big( \, [\, 0, \, T\, ] \, ; H^2 \big) \cap
L^2( \, 0, \, T \, ; H^3),\\[5pt]
\partial _t\phi \in C^0\big( \, [\, 0, \, T\, ] \, ; L^2 \big) \cap
L^2( \, 0, \, T \, ; H^1)
\end{array} 
\right.\,
\end{eqnarray*}
is a solution of (3.4) for
some $T>0$, and for simplicity, in Sections 4 and 5, 
we use the notation $C_{\phi_0}$ to denote 
positive 
constants which may depend on the initial data $\phi_0 \in H^2$,
and the shape of the equation but not depend on $T$. 
Next we give the proofs of the time-decay estimates in Sections 6-9. 
The time-decay results, Theorems 3.5-3.8 are obtained 
by the technical time-weighted energy method 
(cf. \cite{yoshida1}, \cite{yoshida2}, \cite{yoshida4}) 
with the help of the uniform boundedness of 
$\| \, \phi \, \|_{L^{\infty}}$, $\| \, \psi \, \|_{L^{\infty}}$, 
$\| \, \partial_{x}\phi \, \|_{L^{\infty}}$ (see Lemmas 4.1 and 5.3) 
and $\| \, \partial_{x}\psi \, \|_{L^{\infty}}$, 
and the uniform energy estimates for $\phi$ and $\psi$ 
from (3.7), (3.8) and Theorem 3.4. 
In particular, $L^q$-estimates for $\phi$ and $\psi$ with $1<q<2$ 
need subtle and delicate treatment 
(see Lemmas 7.3, 7.5 and 7.6). 
Moreover, in order to obtain the time-decay estimates 
for the derivatives, Theorem 3.8, 
we apply the bootstrap argument 
by Yoshida \cite{yoshida2}, \cite{yoshida4} 
(see the process (8.16)-(8.23)). 
Since the proof of Theorem 3.7 is similarly given as the proofs of 
Theorems 3.5 and 3.6, 
we only show Theorems 3.5, 3.6 and 3.8 in Sections 6-9. 
We finally discuss and compare these time-decay resuts in Section 10. 
In the following sections, 
we rewrite $U^r_{1}$ 
as $U^r$ for simplicity. 

\bigskip 

\noindent
\section{{\it A priori} estimates I}
We first recall that, 
by the arguments in Matsumura-Yoshida \cite{matsumura-yoshida'}, 
the uniform boundedness 
of $L^{\infty}$-norm of the deviation $\phi$ can be obtained 
without using the maximum principle 
(cf. \cite{ilin-kalashnikov-oleinik}, \cite{lad-sol-ura}, \cite{lions}, \cite{matsumura}). 
In fact, according to 
\cite{matsumura-yoshida'}, 
the uniform boundedness of $\| \, \phi \, \|_{L^{\infty}}$ 
and the basic $L^2$-energy estimate for $\phi$ were obtained 
only by using 
a standard energy method. 
Therefore, we similarly get the following lemma and proposition. 

\medskip

\noindent
{\bf Lemma 4.1.} 
\quad {\it
There exists a positive constant $C_{\phi_{0}}$ such that 
$$
\|\, \phi(t) \,\|_{L^{\infty}} \leq C_{\phi_{0}}
\qquad \big( \, t \in [\, 0, \, T\, ] \, \big).
$$
}

\medskip

\noindent
{\bf Proposition 4.2.}\quad {\it
For $0<p<1$, there exists a positive constant $C_{\phi_{0}}$
such that 
\begin{align*}
\begin{aligned}
\| \, \phi(t) \, \|_{{L}^{2}}^{2} 
 &+ \int^{t}_{0} 
    \big|\big| \, 
    \big(\, \sqrt{\partial_x U^r}\: \phi \, \big)(\tau) 
    \, \big|\big|_{L^2}^{2} \, \mathrm{d}\tau\\
&+ \int^{t}_{0} \int ^{\infty }_{-\infty } 
     \braket{ \, \partial_x \phi \, }^{p-1} 
|\, \partial_x \phi \, |^2 
\, \mathrm{d}x\mathrm{d}\tau \leq C_{\phi_{0}}
\qquad \big( \, t \in [\, 0, \, T\, ] \, \big),
\end{aligned}
\end{align*}
where $\braket{ \, s \, }:= (1+s^2)^{1/2}\ (s \in \mathbb{R})$. 
}

\medskip

With the help of Lemma 4.1 and Proposition 4.2, 
we now show the {\it a pripri} estimate of $\partial_{x}\phi$ 
as follows. 

\medskip

\noindent
{\bf Proposition 4.3.} \quad {\it
For $0<p<1$, there exists a positive constant $C_{\phi_{0}}$ 
such that 
\begin{align*}
\begin{aligned}
& \int ^{\infty }_{-\infty } 
     \braket{ \, \partial_x \phi \, }^{p-1} 
     \, |\, \partial_x \phi \, |^2
\, \mathrm{d}x
   + \int^{t}_{0} 
     \| \, \partial_t \phi(\tau) \, \|_{L^2}^2
\, \mathrm{d}\tau \\[5pt]
&\leq C_{\phi_{0}} \, 
     \braket{\braket{\, 
     \partial_x\phi 
     \, }}_{L^{\infty}_{t,x}}^{1-p}
\qquad \big( \, t \in [\, 0, \, T\, ] \,\big), 
\end{aligned}
\end{align*}
where 
$$
\displaystyle{
\braket{\braket{\, v  \, }}_{L^{\infty}_{t,x}}
:= \esssup\displaylimits_{t\in [0,T], \; x \in \mathbb{R}} 
\, \braket{\, v(t, x)  \, }
}. 
$$
}

\medskip

{\bf Proof of Proposition 4.3.}
Multiplying the equation in (3.4) by 
$\partial_t \phi$, and integrating the resultant formula 
with respect to $x$, we have, after integration by parts, 
\begin{align}
\begin{aligned}
&\frac{\mathrm{d}}{\mathrm{d}t} \, 
\int ^{\infty }_{-\infty } \int ^{\partial_x \phi}_{0} 
\Big( \, 
\sigma\bigl( \, 
   \eta+ \partial_x U^r \, \bigr) - 
    \sigma\bigl( \, 
    \partial_x U^r  \, 
    \bigr)
    \, \Big)
\, \mathrm{d}\eta \, \mathrm{d}x 
+ \Vert \, \partial_t \phi(t) \, \Vert_{L^{2}}^{2} \\
&
= - \int ^{\infty }_{-\infty } 
\partial_t \phi \, 
  \partial_x \bigl( \, f(\phi+U^r) - f(U^r) \, \bigr) \, \mathrm{d}x
+\int ^{\infty }_{-\infty } 
\partial_t \phi \, 
  \partial_x 
  \Big( \, 
    \sigma\bigl( \, 
    \partial_x U^r  \, 
    \bigr)
    \, \Big)
\, \mathrm{d}x. 
\end{aligned}
\end{align}
We first note that 
\begin{equation}
\int ^{\infty }_{-\infty } \int ^{\partial_x \phi}_{0} 
\Big( \, 
\sigma\bigl( \, 
   \eta+ \partial_x U^r \, \bigr) - 
    \sigma\bigl( \, 
    \partial_x U^r  \, 
    \bigr)
    \, \Big)
\, \mathrm{d}\eta \, \mathrm{d}x
\sim \int ^{\infty }_{-\infty } 
     \braket{ \, \partial_x \phi \, }^{p-1} 
     \, |\, \partial_x \phi \, |^2
  \, \mathrm{d}x. 
\end{equation} 
Noting the uniform boundedness of $\| \, \phi \, \|_{L^{\infty}}$, Lemma 4.1, 
we estimate each terms on the right-hand side of (4.1) 
by the Young inequality as 
\begin{align}
\begin{aligned}
&\left| \,  
\int ^{\infty }_{-\infty } 
\partial_t \phi \, 
  \partial_x \bigl( \, f(\phi+U^r) - f(U^r) \, \bigr) \, \mathrm{d}x
\, \right|\\
&\leq \epsilon \, 
     \Vert \, \partial_t \phi(t) \, \Vert_{L^{2}}^{2} 
+ C_{\epsilon} \, \int ^{\infty }_{-\infty } 
  \big| \, 
  f'(\phi+U^r) - f'(U^r)
  \, \big|^2 \, 
  | \, \partial_xU^r \, |^2 
  \, \mathrm{d}x\\
&\quad + C_{\epsilon} \, \int ^{\infty }_{-\infty } 
  \big| \, 
  f'(\phi+U^r) 
  \, \big|^2 \, 
  | \, \partial_x\phi \, |^2 
  \, \mathrm{d}x\\
&\leq \epsilon \, 
     \Vert \, \partial_t \phi(t) \, \Vert_{L^{2}}^{2} 
+ C_{\phi_{0}, \epsilon} \, 
  \int ^{\infty }_{-\infty } 
  \phi^2 \, 
  | \, \partial_xU^r \,|^2
  \, \mathrm{d}x\\
&\quad 
+ C_{\phi_{0}, \epsilon} \, 
     \Vert \, \partial_x \phi(t) \, \Vert_{L^{2}}^{2}
\quad (\epsilon>0), 
\end{aligned}
\end{align} 
\begin{align}
\begin{aligned}
&\left| \, \int ^{\infty }_{-\infty } 
\partial_t \phi \, 
  \partial_x 
  \Big( \, 
    \sigma\bigl( \, 
    \partial_x U^r  \, 
    \bigr)
    \, \Big)
\, \mathrm{d}x \, \right|\\
&\quad 
\leq \epsilon \, 
     \Vert \, \partial_t \phi(t) \, \Vert_{L^{2}}^{2} 
+ C_{\epsilon} \, 
  \Vert \, \partial_x^2 U^r(t) \, \Vert_{L^{2}}^{2}
\quad (\epsilon>0). 
\end{aligned}
\end{align} 
Noting from Lemmas 2.2 and 3.5 that
\begin{align}
\begin{aligned}
&\int ^{\infty }_{-\infty } 
  \phi^2 \, 
  | \, \partial_xU^r \,|^2
  \, \mathrm{d}x\\
& \quad 
\leq C \, \min \bigg\{ \, 
  (1+t)^{-1} \, 
  \int ^{\infty }_{-\infty } 
  \phi^2 \, 
  \partial_xU^r 
  \, \mathrm{d}x, \, 
  (1+t)^{-2} \, \Vert \, \phi(t) \, \Vert_{L^{2}}^{2} 
 \, \bigg \}, 
\end{aligned}
\end{align} 
substituting (4.2)-(4.4) into (4.1) 
and integrating the resultant formula with respect to $t$, 
we obtain 
\begin{align}
\begin{aligned}
& \int ^{\infty }_{-\infty } 
     \braket{ \, \partial_x \phi \, }^{p-1} 
     \, |\, \partial_x \phi \, |^2
\, \mathrm{d}x
   + \int^{t}_{0} 
     \| \, \partial_t \phi(\tau) \, \|_{L^2}^2
\, \mathrm{d}\tau \\[5pt]
&\leq C_{\phi_{0}} + C_{\phi_{0}} \, 
     \braket{\braket{\, 
     \partial_x\phi 
     \, }}_{L^{\infty}_{t,x}}^{1-p}\, 
     \int^{t}_{0} \int ^{\infty }_{-\infty } 
     \braket{ \, \partial_x \phi \, }^{p-1} 
     \, |\, \partial_x \phi \, |^2
\, \mathrm{d}x\mathrm{d}\tau
\end{aligned}
\end{align}
provided $\epsilon$ is suitably small. 
Then, noting 
Proposition 4.2, 
we obtain the desired {\it a priori} estimate for $\partial_x \phi$. 

Thus, the proof of Proposition 4.3 is completed. 

\bigskip 

\noindent
\section{{\it A priori} estimates I\hspace{-.1em}I}
In this section, 
we proceed to the {\it a priori} estimate 
for the derivatives $\partial_t \phi$ and $\partial_x^2 \phi$. 
We further establish the uniform 
boundedness of $\partial_x \phi$, 
and then accomplish the proof of Theorem 3.3. 

\medskip

\noindent
{\bf Proposition 5.1.}\quad {\it
For $0<p<1$, there exists a positive constant $C_{\phi_{0}}$
such that 
\begin{align*}
\begin{aligned}
& \| \, \partial_t \phi(t) \, \|_{L^2}^2
   + \int^{t}_{0} 
     \int ^{\infty }_{-\infty } 
     \braket{ \, \partial_x \phi \, }^{p-1} \, 
     | \, \partial_t\partial_x \phi \, |^2
\, \mathrm{d}x\mathrm{d}\tau \\[5pt]
&\leq C_{\phi_{0}} \, 
     \braket{\braket{\, 
     \partial_x\phi 
     \, }}_{L^{\infty}_{t,x}}^{2(1-p)}
\qquad \big( \, t \in [\, 0, \, T\, ] \,\big).
\end{aligned}
\end{align*}
}

\medskip

{\bf Proof of Proposition 5.1.}
Differentiating the equation in (3.4) with respect to $t$, we have 
\begin{align}
\begin{aligned}
&\partial _t^2\phi 
 + \partial _t\partial_x \big( \, f(\phi+U^r) - f(U^r) \, \big) \\[5pt]
&- \partial _t\partial_x 
    \Big( \, 
    \sigma \big( \, \partial_x \phi + \partial_x U^r \, \big) 
    - \sigma \big( \, \partial_x U^r  \, \big)  \, 
    \Big)
    = \partial _t\partial_x 
     \Big( \, \sigma \big( \, \partial_x U^r  \, \big)  \, \Big). 
\end{aligned}
\end{align}
Multiplying (5.1) by 
$\partial _t\phi$
and integrating the resultant formula with respect to $x$, we have,
after integration by parts,
\begin{align}
\begin{aligned}
&\frac{1}{2} \, \frac{\mathrm{d}}{\mathrm{d}t} \, 
\| \, \partial_t \phi(t) \, \|_{L^2}^2
-\int ^{\infty }_{-\infty } 
  \partial _t \partial _x \phi \, 
  \partial _t \big( \, f(\phi+U^r) - f(U^r) \, \big) 
    \, \mathrm{d}x\\
&+\int ^{\infty }_{-\infty } 
  \partial _t \partial _x \phi \, 
  \partial _t \Big( \, 
  \sigma \big( \, \partial_x\phi\ + \partial_x U^r \, \big) - 
  \sigma \big( \, \partial_x U^r \, \big) \, \Big)
  \, \mathrm{d}x\\
&= -\int ^{\infty }_{-\infty } 
  \partial _t \partial _x \phi \, 
  \partial _t \Big( \, 
  \sigma \big( \, \partial_x U^r \, \big) \, \Big)
  \, \mathrm{d}x. 
\end{aligned}
\end{align}
Noting Lemmas 2.2 and 4.1, 
and using (4.5) and the Young inequality, we estimate the second term 
on the left-hand side of (5.1) as 
\begin{align}
\begin{aligned}
& \left| \,  
  \int ^{\infty }_{-\infty } 
  \partial _t \partial _x \phi \, 
  \partial _t \big( \, f(\phi+U^r) - f(U^r) \, \big) 
    \, \mathrm{d}x
\, \right|\\
& \leq \epsilon \, 
       \int ^{\infty }_{-\infty } 
       \braket{ \, \partial_x \phi \, }^{p-1} \, 
       | \, \partial_t\partial_x \phi \, |^2
       \, \mathrm{d}x\\
& \quad 
    + C_{\epsilon} \, 
      \int ^{\infty }_{-\infty } 
      \braket{ \, \partial_x \phi \, }^{1-p} \,
       \big| \, f'(\phi+U^r) - f'(U^r) \, \big|^2 \, 
       | \, \partial_{t}U^r \, |^2 
    \, \mathrm{d}x \\
& \quad 
    + C_{\epsilon} \, \int ^{\infty }_{-\infty } 
      \braket{ \, \partial_x \phi \, }^{1-p} \,
       \big| \, f'(\phi+U^r) \, \big|^2 \, 
       | \, \partial_{t}\phi \, |^2 
    \, \mathrm{d}x \\
& \leq \epsilon \, 
       \int ^{\infty }_{-\infty } 
       \braket{ \, \partial_x \phi \, }^{p-1} \, 
       | \, \partial_t\partial_x \phi \, |^2
       \, \mathrm{d}x\\
& \quad 
    + C_{\phi_{0}, \epsilon} \, 
      (1+t)^{-1}\, 
      \braket{\braket{\, 
     \partial_x\phi 
     \, }}_{L^{\infty}_{t,x}}^{1-p}\, 
      \int ^{\infty }_{-\infty } 
      \phi^2 \, 
      \partial_{x}U^r 
    \, \mathrm{d}x \\
& \quad 
    + C_{\phi_{0}, \epsilon} \, 
      \braket{\braket{\, 
     \partial_x\phi 
     \, }}_{L^{\infty}_{t,x}}^{1-p}\, 
      \| \, \partial_t \phi(t) \, \|_{L^2}^2 
\quad (\epsilon>0). 
\end{aligned}
\end{align}
\begin{align}
\begin{aligned}
& \int ^{\infty }_{-\infty } 
  \partial _t \partial _x \phi \, 
  \partial _t \Big( \, 
  \sigma \big( \, \partial_x\phi\ + \partial_x U^r \, \big) - 
  \sigma \big( \, \partial_x U^r \, \big) \, \Big)
  \, \mathrm{d}x\\
& \geq \int ^{\infty }_{-\infty } 
  \sigma' \big( \, \partial_x\phi\ + \partial_x U^r \, \big) \, 
   | \, \partial_t\partial_x \phi \, |^2
       \, \mathrm{d}x\\
& \quad 
  - \int ^{\infty }_{-\infty } 
  | \, \partial_t\partial_x \phi \, | \, 
  \Big( \, 
  \sigma' \big( \, \partial_x\phi\ + \partial_x U^r \, \big) - 
  \sigma' \big( \, \partial_x U^r \, \big) \, \Big)
  | \, \partial_t\partial_x U^r \, |
  \, \mathrm{d}x\\
& \geq \big( \, C_{\phi_0}^{-1} - \epsilon \, \big) \, 
  \int ^{\infty }_{-\infty } 
  \braket{ \, \partial_x \phi \, }^{p-1} \, 
   | \, \partial_t\partial_x \phi \, |^2
       \, \mathrm{d}x\\
& \quad 
    - C_{\phi_0, \epsilon} \, 
      \| \, \partial_t\partial_x U^r(t) \, \|_{L^{\infty}}^2 \,
      \int ^{\infty }_{-\infty } 
       \braket{ \, \partial_x \phi \, }^{1-p} \, 
       | \, \partial_x \phi \, |^2
       \, \mathrm{d}x\\     
& \geq \big( \, C_{\phi_0}^{-1} - \epsilon \, \big) \, 
  \int ^{\infty }_{-\infty } 
  \braket{ \, \partial_x \phi \, }^{p-1} \, 
   | \, \partial_t\partial_x \phi \, |^2
       \, \mathrm{d}x\\
& \quad 
    - C_{\phi_0, \epsilon} \, 
    \big( \, 
      \| \, \partial_x U^r(t) \, \|_{L^{\infty}}^4 
      + \| \, \partial_x^2 U^r(t) \, \|_{L^{\infty}}^2
      \, \big)\\
& \qquad \times 
      \braket{\braket{\, 
     \partial_x\phi 
     \, }}_{L^{\infty}_{t,x}}^{2(1-p)}\, 
      \int ^{\infty }_{-\infty } 
       \braket{ \, \partial_x \phi \, }^{p-1} \, 
       | \, \partial_x \phi \, |^2
       \, \mathrm{d}x 
\quad (\epsilon >0). 
\end{aligned}
\end{align}
Similarly, the right-hand side of (5.1) is estimated as 
\begin{align}
\begin{aligned}
& \left| \,
\int ^{\infty }_{-\infty } 
  \partial _t \partial _x \phi \, 
  \partial _t \Big( \, 
  \sigma \big( \, \partial_x U^r \, \big) \, \Big)
  \, \mathrm{d}x
\, \right|\\
& \leq \epsilon \, 
       \int ^{\infty }_{-\infty } 
  \braket{ \, \partial_x \phi \, }^{p-1} \, 
   | \, \partial_t\partial_x \phi \, |^2
       \, \mathrm{d}x
    + C_{\epsilon} \, 
      \int ^{\infty }_{-\infty } 
  \braket{ \, \partial_x \phi \, }^{1-p} \,  
  | \, \partial _t\partial_x U^r  \,|^2  
    \, \mathrm{d}x\\ 
& \leq \epsilon \, 
       \int ^{\infty }_{-\infty } 
  \braket{ \, \partial_x \phi \, }^{p-1} \, 
   | \, \partial_t\partial_x \phi \, |^2
       \, \mathrm{d}x\\
& \quad 
    + C_{\phi_0, \epsilon} \, 
    \big( \, 
      \| \, \partial_x U^r(t) \, \|_{L^{2}}^4 
      + \| \, \partial_x^2 U^r(t) \, \|_{L^{2}}^2
      \, \big)
      \braket{\braket{\, 
     \partial_x\phi 
     \, }}_{L^{\infty}_{t,x}}^{1-p}
\quad (\epsilon >0). 
\end{aligned}
\end{align}
Therefore, 
using Lemma 2.2 and Proposition 4.2, 
substituting (5.3)-(5.5) into (5.1), 
and integrating the resultant formula with respect to $t$, 
we have 
\begin{align}
\begin{aligned}
& \| \, \partial_t \phi(t) \, \|_{L^2}^2
   + \int^{t}_{0} 
     \int ^{\infty }_{-\infty } 
     \braket{ \, \partial_x \phi \, }^{p-1} \, 
     | \, \partial_t\partial_x \phi \, |^2
\, \mathrm{d}x\mathrm{d}\tau \\[5pt]
&\leq C_{\phi_{0}} + 
     C_{\phi_{0}} \, 
     \braket{\braket{\, 
     \partial_x\phi 
     \, }}_{L^{\infty}_{t,x}}^{2(1-p)} 
    + C_{\phi_{0}} \, 
     \braket{\braket{\, 
     \partial_x\phi 
     \, }}_{L^{\infty}_{t,x}}^{1-p} \, 
     \int_{0}^t 
     \| \, \partial_t \phi(\tau) \, \|_{L^2}^2
     \mathrm{d}\tau %
\end{aligned}
\end{align}
provided $\epsilon$ is suitably small. 
Noting Proposition 4.3, 
we obtain the desired {\it a priori} estimate for $\partial_t \phi$. 

Thus, the proof of Proposition 5.1 is completed. 

\medskip
Next, we show the following {\it a priori} estimate for $\partial_x^2 \phi$. 
\medskip

\noindent
{\bf Proposition 5.2.}\quad {\it
For $0<p<1$, there exists a positive constant $C_{\phi_{0}}$
such that 
\begin{align*}
\begin{aligned}
& \int ^{\infty }_{-\infty } 
     \braket{ \, \partial_x \phi \, }^{2(p-1)} \, 
     | \, \partial_x^2 \phi \, |^2
\, \mathrm{d}x\mathrm{d}\tau 
\leq C_{\phi_{0}} \, 
     \braket{\braket{\, 
     \partial_x\phi 
     \, }}_{L^{\infty}_{t,x}}^{2(1-p)}
\qquad \big( \, t \in [\, 0, \, T\, ] \,\big).
\end{aligned}
\end{align*}
}

\medskip

{\bf Proof of Proposition 5.2.}
With the help of Propositions 4.3 and 5.1, 
we see from (3.4) that 
\begin{align}
\begin{aligned}
& \big\| \, 
\sigma' \big( \, \partial_x \phi + \partial_x U^r \, \big)  
\partial_x^2 \phi
\, \big\|_{L^2}^2\\
&\leq
 \| \, \partial _t\phi(t) \, \|_{L^2}^2
 +  \big\| \, 
\partial_x \big( \, f(\phi+U^r) - f(U^r) \, \big) (t)
\, \big\|_{L^2}^2\\
& \quad 
 +  \big\| \, 
\sigma' \big( \, \partial_x \phi + \partial_x U^r \, \big) \,   
\partial_x^2 U^r(t)
\, \big\|_{L^2}^2
\leq C_{\phi_{0}} + C_{\phi_{0}} \, 
     \braket{\braket{\, 
     \partial_x\phi 
     \, }}_{L^{\infty}_{t,x}}^{2(1-p)}
\end{aligned}
\end{align}
This completes the proof of Proposition 5.2. 

\medskip

Now, combining Lemma 4.1, Propositions 4.3 and 5.2, 
we show the following 
uniform boundedness of $\| \, \partial_x \phi \, \|_{L^{\infty}}$ 
which plays the essential role to 
control the nonlinearity of $\sigma$ 
(cf. \cite{kanel}, \cite{matsumura-yoshida'}).

\medskip

\noindent
{\bf Lemma 5.3.} 
\quad {\it For $1/3 <p <1$,
there exists a positive constant $C_{\phi_{0}}$ such that 
$$
\|\,\partial_x\phi (t)\,\|_{L^{\infty}} \le C_{\phi_{0}} 
\quad \big(\,t \in [\,0,\,T\,] \, \big).
$$
}

\medskip

{\bf Proof of Lemma 5.3}.
According to \cite{matsumura-yoshida'}, we have 
\begin{align}
\begin{aligned}
\braket{\braket{\, 
     \partial_x\phi 
     \, }}_{L^{\infty}_{t,x}}^{\frac{3p+1}{2}}
\leq 1 &+ C_{\phi_{0}} \, \Big(\,
      \int ^{\infty }_{-\infty } 
       \braket{\, \partial_x\phi \, }^{p-1} 
     |\, \partial_x\phi \, |^2 
     \, \mathrm{d}x
     \, \Big)^{\frac{1}{2}} 
\\
& 
      \times \Big( \,
      \int ^{\infty }_{-\infty } 
       \braket{\, \partial_x\phi\, }^{2(p-1)} 
     |\, \partial_x^2 \phi \, |^2 
     \, \mathrm{d}x
      \, \Big)^{\frac{1}{2}}, 
\end{aligned}
\end{align}
which deduces from Propositions 4.3 and 5.2 that
\begin{equation}
\braket{\braket{\, 
     \partial_x\phi 
     \, }}_{L^{\infty}_{t,x}}^{\frac{3p+1}{2}}\, \le \,
C_{\phi_{0}} \, \braket{\braket{\, 
     \partial_x\phi 
     \, }}_{L^{\infty}_{t,x}}^{\frac{3}{2}(1-p)}.
\end{equation}
Hence, if we assume
$$
\frac{3\, p+1}{2} > \frac{3}{2}\, (1-p)\qquad \bigg(\,\Leftrightarrow p > \frac{1}{3}\,\bigg),
$$
we obtain for $1/3<p<1$
\begin{equation}
\braket{\braket{\, 
     \partial_x\phi 
     \, }}_{L^{\infty}_{t,x}}
\leq \,C_{\phi_{0}}. 
\end{equation}
Thus, the proof of Lemma 5.3 is completed. 

\medskip

By Lemma 4.1, Propositions 4.2, 4.3, 5.1 and 5.2 with the aid of
Lemma 5.3, we obtain the energy estimate
\begin{align}
\begin{aligned}
& \| \, \phi(t) \, \|_{{H}^{2}}^{2} 
+ \| \, \partial_{t}\phi(t) \, \|_{{L}^{2}}^{2} 
+ \int^{t}_{0} 
   \big\| \, 
   (\, \sqrt{\, \partial_x U^r} \:  \phi \,)(\tau) 
   \, \big\|_{L^2}^{2} \, \mathrm{d}\tau 
\\
& 
+ \int^{t}_{0}
   \big( 
 \left\| \, \partial_{x}\phi(\tau) \, \right\|_{{H}^1}^{2} 
 + \left\| \, \partial_{t}\phi(\tau) \, \right\|_{{H}^1}^{2} 
   \, \big) \, \mathrm{d}\tau 
\leq C_0
\qquad \big( \, t \in [\, 0, \, T\, ] \,\big).
\end{aligned}
\end{align}
Therefore, in order to accomplish the proof of Theorem 3.3, 
it suffices to show the following {\it a priori} estimate 
(the proof is also given in \cite{matsumura-yoshida'}). 

\medskip

\noindent
{\bf Lemma 5.4.} 
\quad {\it For $1/3 <p <1$,
there exists a positive constant $C_{\phi_{0}}$ such that 
\begin{equation*}
\int_0^t \| \, \partial_{x}^3\phi(\tau) \,\|_{L^{2}}^2\, \mathrm{d}\tau 
\le C_{\phi_{0}} \qquad \big( \, t \in [\, 0, \, T\, ] \,\big).
\end{equation*}
}

\medskip
Thus, we complete the proof of Theorem 3.3. 
\medskip

In order to accomplish the proof of Theorem 3.1, 
in particular, the asymptotic behavior (3.14) and (3.15), 
we finally show (3.13), that is, 
\medskip

\noindent
{\bf Lemma 5.5.} 
\quad {\it For $1/3 <p <1$, 
there exists a positive constant $C_{\phi_{0}}$ such that 
\begin{equation*}
\int _0^{\infty }\bigg|\,\frac{\mathrm{d}}{\mathrm{d}t} 
\|  \, \partial_t\phi(t)  \, \|_{L^2}^2 \,\bigg| \, \mathrm{d}t
\le C_{\phi_{0}}.
\end{equation*}
}

\medskip

{\bf Proof of Lemma 5.5}. 
Direct calculation shows by using (3.4) that 
\begin{align}
\begin{aligned}
&\int _0^{\infty }\bigg|\,\frac{\mathrm{d}}{\mathrm{d}t} 
\|  \, \partial_t\phi(t)  \, \|_{L^2}^2 \,\bigg| \, \mathrm{d}t\\
&\leq 2 \, \int _0^{\infty }\bigg|\, 
      \int _{-\infty }^{\infty }
      \partial_t\partial_x\phi \, 
      \partial_t\big( \, f(\phi+U^r)-f(U^r) \, \big) 
      \, \mathrm{d}x 
      \,\bigg| \, \mathrm{d}t\\
& \quad + 2 \, \int _0^{\infty }\bigg|\, 
      \int _{-\infty }^{\infty }
      \partial_t\partial_x\phi \, 
      \partial_t\Big( \, \sigma\big( \, \partial_{x}\phi+ \partial_{x}U^r \, \big)
      \, \Big) 
      \, \mathrm{d}x 
      \,\bigg| \, \mathrm{d}t\\
&\leq C_{\phi_{0}} \, 
      \int _0^{\infty }
      \bigg( \, 
      \|  \, \partial_t\phi(t)  \, \|_{H^1}^2 
      +\int _{-\infty }^{\infty }
      \phi^2 \, \partial_xU^r \, \mathrm{d}x
      + \|  \, \partial_t\partial_xU^r(t)  \, \|_{L^2}^2
      \, \bigg) 
      \, \mathrm{d}t\\
&\leq C_{\phi_{0}}.
\end{aligned}
\end{align}
This completes the proof of Lemma 5.5. 

\medskip
Thus, we complete the proof of Theorem 3.1. 

%

\medskip

\noindent
{\it Remark 5.6.}\quad 
The estimate in Proposition 5.2 is sharper than 
that in \cite{matsumura-yoshida'}, that is, 
\begin{align}
\begin{aligned}
& \int ^{\infty }_{-\infty } 
     \braket{ \, \partial_x \phi \, }^{2(p-1)} \, 
     | \, \partial_x^2 \phi \, |^2
\, \mathrm{d}x\mathrm{d}\tau 
\leq C_{\phi_{0}} \, 
     \braket{\braket{\, 
     \partial_x\phi 
     \, }}_{L^{\infty}_{t,x}}^{3(1-p)}
\qquad \big( \, t \in [\, 0, \, T\, ] \,\big).
\end{aligned}
\end{align}
Therefore, if substituting the estimate in Proposition 4.3 and (5.13) 
instead of the estimate in Proposition 5.2, 
then we have obtained in \cite{matsumura-yoshida'} that 
\begin{equation}
\braket{\braket{\, 
     \partial_x\phi 
     \, }}_{L^{\infty}_{t,x}}
\leq \,C_{\phi_{0}} 
\end{equation}
for $3/7<p<1$. 

\bigskip 

\noindent
\section{Time-decay estimates I}
In this section, 
we obtain the time-decay estimate 
for $\phi$ 
(not assuming the $L^1$-integrability 
to the initial perturbation), Theorem 3.5. 
To do that, we establish 
the following time-weighted $L^q$-energy estimate to $\phi$ with $2\leq q\leq \infty$.

\medskip

\noindent
{\bf Proposition 6.1.}\quad {\it
Suppose the same assumptions as in Theorem 3.1. 
For any $q \in [\, 2, \, \infty )$, 
there exist positive constants $\alpha$ 
and $C_{\alpha, \phi_{0}, \epsilon}$ 
such that the unique global in time solution $\phi$ 
of the Cauchy problem {\rm (3.4)} has the following $L^q$-energy estimate 
\begin{align*}
\begin{aligned}
&(1+t)^{\alpha}\, \| \, \phi(t) \, \|_{L^q}^q 
+\int_0^t (1+\tau)^{\alpha} \, 
 \int_{-\infty}^{\infty} 
 | \, \phi \, |^{q} 
 \, \partial_{x}U^r 
 \, \mathrm{d}x\mathrm{d}\tau\\
&+\int_0^t (1+\tau)^{\alpha} \, \int_{-\infty}^{\infty} 
 | \, \phi\,|^{q-2} \, | \, \partial_{x}\phi\,|^2
 \, \mathrm{d}x \mathrm{d}\tau \\
&\le C_{\phi_{0}}\, \| \,\phi_0 \,\|_{L^q}^q 
+ C_{\alpha, \phi_{0}}\, (1+t)^{\alpha- \frac{q-2}{4}}
\quad \bigl( \, t \ge 0 \, \bigr). 
\end{aligned}
\end{align*}
}

\medskip

The proof of Proposition 6.1 is provided by the following two lemmas. 

\medskip

\noindent
{\bf Lemma 6.2.}\quad {\it
For any $q \in [\, 1, \, \infty )$, 
there exist positive constants $\alpha$ and $C_{\phi_{0}}$ such that 
\begin{align}
\begin{aligned}
&(1+t)^{\alpha}\, \| \, \phi(t) \, \|_{L^q}^q 
+\int_0^t (1+\tau)^{\alpha} \, 
 \int_{-\infty}^{\infty} 
 | \, \phi \, |^{q} 
 \, \partial_{x}U^r 
 \, \mathrm{d}x\mathrm{d}\tau\\
&+\int_0^t (1+\tau)^{\alpha} \, \int_{-\infty}^{\infty} 
 | \, \phi\,|^{q-2} \, | \, \partial_{x}\phi\,|^2
 \, \mathrm{d}x \mathrm{d}\tau \\
&\le C_{\phi_{0}}\, \| \,\phi_0 \,\|_{L^q}^q 
+ C_{\phi_{0}} \cdot \alpha \, \int_0^t (1+\tau)^{\alpha-1} \, 
  \| \, \phi(\tau) \, \|_{L^q}^q \, \mathrm{d}\tau\\
&+ C_{\phi_{0}} 
       \, \int_0^t (1+\tau)^{\alpha} \, 
       \| \, \phi(\tau) \, \|_{L^{\infty}}^{q-1} \, 
       \| \, \partial_{x}^2U^r(\tau) \, \|_{L^{1}}
       \, \mathrm{d}\tau.
\quad \bigl( \, t \ge 0 \, \bigr). 
\end{aligned}
\end{align}
}

\bigskip

\noindent
{\bf Lemma 6.3.}\quad {\it
Assume $2\leq q < \infty$ and $2\leq r < \infty$. 
We have the following interpolation inequalities.  

 \noindent
 {\rm (1)}\ \ It follows that 
 \begin{align*}
 \begin{aligned}
 &\Vert \, \phi  (t) \, \Vert _{L^{\infty} }
  \leq 
  \left( \, 
  \frac{q+2}{2}
   \, \right)^{\frac{2}{q+2}} \, 
   \left( \, \int _{-\infty}^{\infty} | \, \phi  \, |^2 \, \mathrm{d}x \, \right)
  ^{\frac{1}{q+2}} \\
 & \qquad \quad \quad \quad \quad \; \: \, 
         \times 
         \left( \, \int _{-\infty}^{\infty} | \, \phi  \, |^{q-1} 
         | \, \partial_{x}\phi \, |^{2} \, \mathrm{d}x \, \right)
         ^{\frac{1}{q+2}} \quad \big( \, t \ge 0 \, \big). 
 \end{aligned}
 \end{align*}

 \noindent
 {\rm (2)}\ \ 
 It follows that 
 \begin{align*}
 \begin{aligned}
 &\Vert \, \phi  (t) \, \Vert _{L^{r} }
  \leq 
  \left( \, 
  \frac{q+2}{2}
   \, \right)^{\frac{2(r-2)}{(q+2)r}} \, 
   \left( \, \int _{-\infty}^{\infty} | \, \phi  \, |^2 \, \mathrm{d}x \, \right)
  ^{\frac{r-2}{(q+2)r}} \\
 & \qquad \quad \quad \quad \quad \; \: \, 
         \times 
         \left( \, \int _{-\infty}^{\infty} | \, \phi  \, |^{q-1} 
         | \, \partial_{x}\phi \, |^{2} \, \mathrm{d}x \, \right)
         ^{\frac{q+r}{(q+2)r}} \quad \big( \, t \ge 0 \, \big). 
 \end{aligned}
 \end{align*}
}

\medskip

In what follows, we first prove Lemma 6.2 and Lemma 6.3, 
and finally give the proof of Proposition 6.1, 
and further the $L^{\infty}$-estimate of $\phi$. 

\medskip

{\bf Proof of Lemma 6.2}.
Performing the computation 
$$
\int_0^t (1+\tau)^{\alpha} \, \int_{-\infty}^{\infty}
(3.4)\times |\, \phi  \,|^{q-2} \, \phi 
\, \mathrm{d}x \mathrm{d}\tau 
\quad \big( \, \alpha>0, \, q\in [\,2, \, \infty)  \, \big)
$$
allows us to get 
\begin{align}
\begin{aligned}
&\frac{1}{q}\, (1+t)^{\alpha}\, \| \, \phi(t) \, \|_{L^q}^q \\
& +(q-1) \, \int_0^t (1+t)^{\alpha} \, 
 \int_{-\infty}^{\infty} \int_{0}^{\phi} 
 \big( \, f'(\eta+U^r)-f'(U^r) \, \big) \\
 & \qquad \qquad \qquad \qquad \qquad \times 
 | \, \eta \, |^{q-2} 
 \, \partial_{x}U^r \, 
 \mathrm{d}\eta \, \mathrm{d}x\mathrm{d}\tau\\
&+(q-1) \, \int_0^t (1+\tau)^{\alpha} \, \int_{-\infty}^{\infty} 
 | \, \phi\,|^{q-2} \, \partial_{x}\phi\\
 & \qquad \qquad \qquad \qquad \qquad \times 
 \Big( \, 
 \sigma\big( \, \partial_{x}\phi+\partial_{x}U^r \, \big)
 -  \sigma\big( \, \partial_{x}U^r \, \big)
 \, \Big)
 \, \mathrm{d}x \mathrm{d}\tau \\
&= \frac{1}{q}\, \| \,\phi_0 \,\|_{L^q}^q 
+ \alpha \, \int_0^t (1+\tau)^{\alpha-1} \, 
  \| \, \phi(\tau) \, \|_{L^q}^q \, \mathrm{d}\tau\\
& \quad 
  + \int_0^t (1+\tau)^{\alpha} \, \int_{-\infty}^{\infty} 
 | \, \phi\,|^{q-2} \, \partial_{x}\phi\,
 \Big( \, 
 \sigma\big( \, \partial_{x}U^r \, \big)
 \, \Big)
 \, \mathrm{d}x \mathrm{d}\tau
\quad \bigl( \, t \ge 0 \, \bigr). 
\end{aligned}
\end{align}
By using the  uniform boundedness in Lemmas 4.1 and 5.3, 
we estimate the second and third terms on the left-hand side of (6.2) 
as follows. 
\begin{align}
\begin{aligned}
&(q-1) \, \int_0^t (1+\tau)^{\alpha} \, 
 \int_{-\infty}^{\infty} \int_{0}^{\phi} 
 \big( \, f'(\eta+U^r)-f'(U^r) \, \big) \\
 & \qquad \qquad \qquad \qquad \qquad \times 
 | \, \eta \, |^{q-2} 
 \, \partial_{x}U^r \, 
 \mathrm{d}\eta \, \mathrm{d}x\mathrm{d}\tau\\
& \geq C_{\phi_{0}}^{-1} 
       \, \int_0^t (1+\tau)^{\alpha} \, 
       \int_{-\infty}^{\infty} 
       | \, \phi \, |^{q} 
       \, \partial_{x}U^r 
       \, \mathrm{d}x\mathrm{d}\tau, 
\end{aligned}
\end{align}
\begin{align}
\begin{aligned}
&(q-1) \, \int_0^t (1+\tau)^{\alpha} \, \int_{-\infty}^{\infty} 
 | \, \phi\,|^{q-2} \, \partial_{x}\phi\\
 & \qquad \qquad \qquad \qquad \qquad \times 
 \Big( \, 
 \sigma\big( \, \partial_{x}\phi+\partial_{x}U^r \, \big)
 -  \sigma\big( \, \partial_{x}U^r \, \big)
 \, \Big)
 \, \mathrm{d}x \mathrm{d}\tau\\
& \geq C_{\phi_{0}}^{-1} 
       \, \int_0^t (1+\tau)^{\alpha} \, 
       \int_{-\infty}^{\infty} 
       | \, \phi \, |^{q-2} 
       \, | \, \partial_{x}\phi \, |^{2} 
       \, \mathrm{d}x\mathrm{d}\tau. 
\end{aligned}
\end{align}
The third term on the right-hand side of (6.2) is also estimated as 
\begin{align}
\begin{aligned}
&\int_0^t (1+\tau)^{\alpha} \, \int_{-\infty}^{\infty} 
 | \, \phi\,|^{q-2} \, \partial_{x}\phi\,
 \Big( \, 
 \sigma\big( \, \partial_{x}U^r \, \big)
 \, \Big)
 \, \mathrm{d}x \mathrm{d}\tau\\
& \leq C 
       \, \int_0^t (1+\tau)^{\alpha} \, 
       \| \, \phi(\tau) \, \|_{L^{\infty}}^{q-1} \, 
       \| \, \partial_{x}^2U^r(\tau) \, \|_{L^{1}}
       \, \mathrm{d}\tau. 
\end{aligned}
\end{align}
Substituting (6.3)-(6.5) into (6.1), we immediately have Lemma 6.2. 

\medskip
Next, we show Lemma 6.3. 
\medskip

{\bf Proof of Lemma 6.3}.
Noting that $\phi (t, \cdot \: ) \in H^2$ 
imply 
$\displaystyle{\lim _{x\rightarrow \pm \infty}\phi(t,\, x)}=0$ 
for $t \ge 0$ and 
\begin{align}
\begin{aligned}
 | \, \phi \, |^{s} 
 &\leq s \int _{-\infty}^{\infty} 
      | \, \phi \, |^{s-1}  \, 
      | \, \partial_{x}\phi \, |
      \, \mathrm{d}x\quad \big(\, s\ge1 \,\big).\\ 
\end{aligned}
\end{align}
By the Cauchy-Schwarz inequality, we have 
\begin{equation}
 | \, \phi \, |^{s} 
 \leq s \, 
       \left( \,
       \int _{-\infty}^{\infty} 
       | \, \phi \, |^{2s-q} 
       \, \mathrm{d}x
       \, \right)^{\frac{1}{2}} \, 
       \left( \, 
       \int _{-\infty}^{\infty} 
       | \, \phi \, |^{q-2} \,
       | \, \partial_{x}\phi \,|^{2}
       \, \mathrm{d}x 
       \, \right)^{\frac{1}{2}}.  
\end{equation}
Taking $s=\frac{q+2}{2}$, we immediately have (1). 
Next, substituting (1) into 
\begin{equation}
 \Vert \, \partial_{x}\phi \, \Vert _{L^{r} }^{r}
 \le \Vert \, \partial_{x}\phi \, \Vert _{L^\infty }^{r-2}
      \Vert \, \partial_{x}\phi \, \Vert _{L^2 }^2, 
\end{equation}
we also have (2). 

Thus, we complete the proof of Lemma 6.3. 

\medskip
We now show Proposition 6.1 with the help of Lemmas 6.2 and 6.3. 
\medskip

{\bf Proof of Proposition 6.1}.
By using (1) in Lemma 6.3, we estimate 
the second term on the right-hand side of (6.1) by the Young inequality as 
\begin{align}
\begin{aligned}
&C_{\phi_{0}} \cdot \alpha \, \int_0^t (1+\tau)^{\alpha-1} \, 
  \| \, \phi(\tau) \, \|_{L^q}^q \, \mathrm{d}\tau\\
&\leq C_{\alpha, \phi_{0}} \, 
       \int_0^t (1+\tau)^{\alpha-1} \, 
   \left( \, \int _{-\infty}^{\infty} | \, \phi  \, |^2 \, \mathrm{d}x \, \right)
  ^{\frac{2q}{q+2}} \\
 & \qquad \quad \quad \quad \quad 
         \times 
         \left( \, \int _{-\infty}^{\infty} | \, \phi  \, |^{q-1} 
         | \, \partial_{x}\phi \, |^{2} \, \mathrm{d}x \, \right)
         ^{\frac{q-2}{q+2}}
       \, \mathrm{d}\tau \\
&\leq \epsilon \, 
      \int_0^t (1+\tau)^{\alpha} \, 
      \int _{-\infty}^{\infty} | \, \phi  \, |^{q-1} 
         | \, \partial_{x}\phi \, |^{2} \, \mathrm{d}x\mathrm{d}\tau \\
& \quad + C_{\alpha, \phi_{0}, \epsilon} \, 
          \int_0^t (1+\tau)^{\alpha-\frac{q+2}{4}} \, 
          \left( \, \int _{-\infty}^{\infty} | \, \phi  \, |^2 \, \mathrm{d}x \, \right)
          ^{\frac{q}{2}}
          \, \mathrm{d}\tau \quad (\epsilon>0). 
\end{aligned}
\end{align}
Considering $r=q$ in (2), in Lemma 6.3, 
the third term on the right-hand side of (6.1) are similarly estimated as 
\begin{align}
\begin{aligned}
& C_{\phi_{0}} 
       \, \int_0^t (1+\tau)^{\alpha} \, 
       \| \, \phi(\tau) \, \|_{L^{\infty}}^{q-1} \, 
       \| \, \partial_{x}^2U^r(\tau) \, \|_{L^{1}}
       \, \mathrm{d}\tau\\
&\leq C_{\phi_{0}} \, 
       \int_0^t (1+\tau)^{\alpha-1} \, 
   \left( \, \int _{-\infty}^{\infty} | \, \phi  \, |^2 \, \mathrm{d}x \, \right)
  ^{\frac{q-1}{q+2}} \\
 & \qquad \quad \quad \quad \quad 
         \times 
         \left( \, \int _{-\infty}^{\infty} | \, \phi  \, |^{q-1} 
         | \, \partial_{x}\phi \, |^{2} \, \mathrm{d}x \, \right)
         ^{\frac{q-1}{q+2}}
       \, \mathrm{d}\tau \\
&\leq \epsilon \, 
      \int_0^t (1+\tau)^{\alpha} \, 
      \int _{-\infty}^{\infty} | \, \phi  \, |^{q-1} 
         | \, \partial_{x}\phi \, |^{2} \, \mathrm{d}x\mathrm{d}\tau \\
& \quad + C_{\phi_{0}, \epsilon} \, 
          \int_0^t (1+\tau)^{\alpha-\frac{q+2}{3}} \, 
          \left( \, \int _{-\infty}^{\infty} | \, \phi  \, |^2 \, \mathrm{d}x \, \right)
          ^{\frac{q-1}{3}}
          \, \mathrm{d}\tau \quad (\epsilon>0).
\end{aligned}
\end{align}
Substituting (6.9)-(6.10) into (6.1) and choosing $\epsilon$ suitably small, 
we then obtain the time-weighted $L^q$-energy estimate of $\phi$, 
in particular, for $q \in [\, 2, \, \infty)$, 
\begin{equation}
\| \, \phi(t) \, \|_{L^q} 
\leq C_{\phi_{0}} \, 
(1+t)^{-\frac{1}{4}\left(1- \frac{2}{q} \right)} 
\quad \big( \, t\geq0 \, \big)
\end{equation}
by choosing $\alpha$ suitably large. 

Thus, we complete the proof of Proposition 6.1. 

\medskip

\noindent
We finally show the $L^{\infty}$-estimate of $\phi$. 
We use the following Gagliardo-Nirenberg inequality: 
\begin{align}
\| \,\phi(t) \,\|_{L^\infty }
\le C_{q,\theta} \, 
\| \,\phi(t) \,\|_{L^{q}}^{1-\theta}
\| \,\partial_{x}\phi(t) \,\|_{L^{2}}^{\theta}
\end{align}
for any 
$q\in [\, 2,\, \infty )$, $\theta \in (0,\,1\, ]$ satisfying 
$$
\frac{\theta}{2} =\frac{1-\theta}{q}. 
$$
By using (6.12), 
we immediately have 
\begin{align}
\begin{aligned}
\| \,\phi(t) \,\|_{L^\infty }
&\le C_{\theta, \phi_0, \epsilon } \, 
     (1+t)^{- \frac{1}{4}\left(1- \frac{2}{q} \right)(1-\theta)}\\
&\le C_{\phi_0, \epsilon } \, 
     (1+t)^{- \frac{1}{4} + \epsilon} \quad \big( \, t \ge 0 \, \big)
\end{aligned}
\end{align}
for $\epsilon>0$. 
Consequently, we do complete the proof of 
$L^{q}$-estimate of $\phi$ with $2\le q \le \infty$.

Thus, the proof of Theorem 3.5 is completed. 

\bigskip 

\noindent
\section{Time-decay estimates I\hspace{-.1em}I}
In this section, 
we obtain the 
time-growth and time-decay estimates 
for $\phi$ with $1\leq q\leq \infty$ 
in the case where $\phi_0 \in L^1 \cap H^2$, Theorem 3.6. 
In the following, we are going to show the estimates, step by step, 
in the cases where $q=1$, $1< q<2$ and $2\le q \leq \infty$ in this order. 
We first establish the $L^1$-estimate of $\phi$. 
To do that, we use the Friedrichs mollifier $\rho_\delta \ast $, 
where 
$\rho_\delta(\phi):=\frac{1}{\delta}\rho \left( \frac{\phi}{\delta}\right)$
with 
\begin{align*}
\begin{aligned}
&\rho \in C^{\infty}_0(\mathbb{R}), \quad
\rho (\phi)\geq 0 \quad (\phi \in \mathbb{R}),
\\
&\mathrm{supp} \{\rho \} \subset 
\left\{\phi \in \mathbb{R}\: \left|\:  |\,\phi \, |\le 1 \right. \right\},\quad  
\int ^{\infty}_{-\infty} \rho (\phi)\, \mathrm{d}\phi=1.
\end{aligned}
\end{align*}

We then recall the following properties of the mollifier 
(see \cite{hashimoto-kawashima-ueda}, 
\cite{yoshida1}, \cite{yoshida2}, \cite{yoshida5} and so on). 
\medskip

\noindent
{\bf Lemma 7.1.}\quad{\it

\noindent
{\rm (1)}\ 
$\displaystyle{\lim_{\delta \to 0}\, 
\left( \, \rho_\delta \ast \mathrm{sgn} \, \right)(\phi)}
= \mathrm{sgn} (\phi)
\qquad (\phi \in \mathbb{R}),$

\noindent
{\rm (2)}\ 
$\displaystyle{\lim_{\delta \to 0}\, 
\int^{\phi}_0\left( \, \rho_\delta \ast \mathrm{sgn} \, \right)(\eta)
\, \mathrm{d}\eta}
=|\, \phi \, |
\qquad (\phi \in \mathbb{R}),$

\medskip

\noindent
{\rm (3)}\ 
$\big. \left( \, \rho_\delta \ast \mathrm{sgn} \, \right) \bigr|_{\phi=0} =0,$

\medskip

\noindent
{\rm (4)}\ 
$\displaystyle{ \frac{\mathrm{d}}{\mathrm{d}\phi}\, 
\left( \, \rho_\delta \ast \mathrm{sgn} \, \right)(\phi)}
=2\, \rho_\delta(\phi)
\ge 0
\qquad (\phi \in \mathbb{R}),$

\medskip

\noindent
where 
$$
\left( \, \rho_\delta \ast \mathrm{sgn} \, \right)(\phi)
:=\int^{\infty}_{-\infty}
\rho_\delta(\phi-\psi)\, \mathrm{sgn}(\psi)\, \mathrm{d}\psi
\qquad (\phi \in \mathbb{R})
$$
and 
\begin{equation*}
\mathrm{sgn}(\phi):=\left\{
\begin{array}{ll}
-1& \quad \; \bigl(\, \phi < 0 \, \bigr),\\[3pt]
\, \, 0& \quad \; \bigl(\, \phi = 0 \, \bigr),\\[3pt]
\, \, 1& \quad \; \bigl(\, \phi > 0 \, \bigr).
\end{array}
\right.
\end{equation*} 
}

\noindent
Making use of Lemma 7.1, 
we obtain the following $L^1$-energy estimate. 

\medskip

\noindent
{\bf Proposition 7.2.}\quad {\it
Suppose the same assumptions as in Theorem 3.1. 
the unique global in time solution $\phi$ 
of the Cauchy problem {\rm (3.4)} has the following $L^1$-energy estimate
\begin{align*}
\begin{aligned}
\| \, \phi(t) \, \|_{L^1}
\le C_{\phi_{0}}\, \max \big\{ \, 1, \, \ln (1+t) \, \big\}
\quad \bigl( \, t \ge 0 \, \bigr). 
\end{aligned}
\end{align*}
}

\medskip

\noindent
{\bf Proof of Proposition 7.2.}\quad
Multiplying the equation in (3.4) by 
$\left( \rho_\delta \ast \mathrm{sgn} \right)(\phi)$ 
and integrating the resultant formula with respect to $x$ and $t$, 
we obtain 
\begin{align}
\begin{aligned}
&\int ^{\infty}_{-\infty} \int^{\phi }_0 
 \left( \, \rho_\delta \ast \mathrm{sgn} \, \right)(\eta)
 \, \mathrm{d}\eta \, \mathrm{d}x \\
&+ 2 \, \int _0^t \int ^{\infty}_{-\infty} \int^{\phi}_0
   \bigl( \, f'(\eta+U^r)-f'(U^r) \, \bigr) \, 
   \rho_\delta(\eta) \, 
   \partial _x U^r 
   \, \mathrm{d}\eta 
   \, \mathrm{d}x \mathrm{d}\tau \\
&+ 
   \int _0^t \int ^{\infty}_{-\infty} 
   \rho_\delta(\phi) \, 
   \partial _x \phi \, 
    \Big( \, \sigma\bigl(\, \partial _x \phi + \partial _x U^r \, \bigr)
    -\sigma\bigl(\, \partial _x U^r \, \bigr) 
    \, \Big)
   \, \mathrm{d}x \mathrm{d}\tau \\
&= \int ^{\infty}_{-\infty} \int^{\phi _0}_0 
   \left( \, \rho_\delta \ast \mathrm{sgn} \, \right)(\eta)
   \, \mathrm{d}\eta \, \mathrm{d}x \\
&\quad 
 + \int _0^t \int ^{\infty}_{-\infty} 
   \left( \, \rho_\delta \ast \mathrm{sgn} \, \right)(\phi) \, 
   \partial _x 
   \Bigl( \, \sigma\bigl(\, \partial _x U^r \, \bigr)  \, \Bigr)
   \, \mathrm{d}x \mathrm{d}\tau. 
\end{aligned}
\end{align}
By using Lemma 7.1, and uniform boundedness in Lemma 4.1 and 5.3, 
we estimate the second and third terms on the left-hand side 
of (7.1) as follows. 
\begin{align}
\begin{aligned}
& 2 \, \int _0^t \int ^{\infty}_{-\infty} \int^{\phi}_0
   \bigl( \, f'(\eta+U^r)-f'(U^r) \, \bigr) \, 
   \rho_\delta(\eta) \, 
   \partial _x U^r 
   \, \mathrm{d}\eta 
   \, \mathrm{d}x \mathrm{d}\tau \\
& \ge C_{\phi_{0}}^{-1} \, 
      \int _0^t \int ^{\infty }_{-\infty } 
      \bigg|\, 
      \int^{| \phi |}_0 
      \eta \, \rho_\delta (\eta ) \, 
      \mathrm{d}\eta 
      \, \bigg|
      \, \partial _x U^r 
      \, \mathrm{d}x \mathrm{d}\tau,%
\end{aligned}
\end{align}
\begin{align}
\begin{aligned}
& \int _0^t \int ^{\infty}_{-\infty} 
   \rho_\delta(\phi) \, 
   \partial _x \phi \, 
    \Big( \, \sigma\bigl(\, \partial _x \phi + \partial _x U^r \, \bigr)
    -\sigma\bigl(\, \partial _x U^r \, \bigr) 
    \, \Big)
   \, \mathrm{d}x \mathrm{d}\tau \\
& \ge C_{\phi_{0}}^{-1} \, 
      \int _0^t \int ^{\infty }_{-\infty } 
      \rho_\delta (\phi ) \, 
      | \, \partial _x \phi \, |^2
      \, \mathrm{d}x \mathrm{d}\tau.%
\end{aligned}
\end{align}
Also from Lemma 7.1, we have 
\begin{align}
\bigg| \, 
\int^{\phi}_0\left( \, \rho_\delta \ast \mathrm{sgn} \, \right)(\eta)
\, \mathrm{d}\eta 
\, \bigg|
\le \left( \, \rho_\delta \ast \mathrm{sgn} \, \right)
    \big( \, |\,  \phi \,| \, \big)\, |\,  \phi \,|
\le |\,  \phi \,|, 
\end{align}
\begin{align}
\displaystyle {\lim_{\delta\rightarrow 0}}
\int ^{\infty}_{-\infty} \int^{\phi (t) }_0
\left( \, \rho_\delta \ast \mathrm{sgn} \, \right)(\eta)
\, \mathrm{d}\eta 
= \| \,\phi(t) \,\|_{L^1} \quad \bigl( \, t \ge 0 \, \bigr). 
\end{align}
Therefore, 
with the help of (7.4)-(7.5), 
substituting (7.2)-(7.3) into (7.1) 
and taking the limit $\delta \rightarrow 0$ 
to the resultant formula, 
we conclude that 
\begin{align}
\begin{aligned}
\| \,\phi (t) \,\|_{L^1} 
&\le \| \,\phi_0 \,\|_{L^1} \\
&\quad + \displaystyle {\lim_{\delta\rightarrow 0}}
   \int^t_0 \left|\, \int ^{\infty}_{-\infty} 
   \left( \, \rho_\delta \ast \mathrm{sgn} \, \right)(\phi) \, 
   \partial _x 
   \Bigl( \, \sigma\bigl(\, \partial _x U^r \, \bigr)  \, \Bigr)
   \, \mathrm{d}x \, \right|\, \mathrm{d}\tau\\
&\le \| \,\phi_0 \,\|_{L^1} 
 + C \, 
   \int^t_0 
   \| \, \partial _x^2 U^r(\tau) \, \|_{L^1} 
   \, \mathrm{d}\tau\\
&\le \| \,\phi_0 \,\|_{L^1} 
 + C \, \ln (1+t). 
\end{aligned}
\end{align}

This completes the proof of Proposition 7.2. 
\medskip

Next, we establish the $L^q$-estimate of $\phi$ in the case where $1<q<2$. 
The estimate in this case needs most subtle and delicate treatment 
than the all other cases (that is, the cases where $q=1$, $2 \leq q \leq \infty$), 
in particular, 
to obtain the interpolation inequalities in Lemma 7.6 with the aid of 
the third term on the left-hand side 
of the time-weighted $L^q$-estimate in Proposition 7.4 
(for a simpler way to obtain the estimate, see Remark 7.8, in particular, (7.31)). 
Therefore, we also introduce  
$\rho_\delta \ast |\, \cdot  \,|^s\,\mathrm{sgn}$ and 
$\rho_\delta \ast |\, \cdot  \,|^{s+1}$, 
and prepare their useful properties as follows. 
\medskip

\noindent
{\bf Lemma 7.3.}\quad{\it
For any $s\geq0$, we have the following properties. 
\medskip

\noindent
{\rm (1)}\ 
$\displaystyle{\lim_{\delta \to 0}\, 
\big( \, \rho_\delta \ast |\, \cdot  \,|^s\,\mathrm{sgn} \, \big)(\phi)}
= |\, \phi  \,|^s\, \mathrm{sgn} (\phi)
\qquad (\phi \in \mathbb{R}),$

\noindent
{\rm (2)}\ 
$\displaystyle{\lim_{\delta \to 0}\, 
\int^{\phi}_0
\big( \, \rho_\delta \ast |\, \cdot  \,|^s\,\mathrm{sgn} \, \big)(\eta)
\, \mathrm{d}\eta}
=\frac{1}{s+1}\, |\, \phi \, |^{s+1}
\qquad (\phi \in \mathbb{R}),$

\medskip

\noindent
{\rm (3)}\ 
$\Big. 
\big( \, \rho_\delta \ast |\, \cdot  \,|^s\,\mathrm{sgn} 
\, \big) \Bigr|_{\phi=0} =0,$

\medskip

\noindent
{\rm (4)}\ 
$\big| \, 
\big( \, \rho_\delta \ast |\, \cdot  \,|^s\,\mathrm{sgn} 
\, \big)(\phi) \, \big|
\leq \big( \, \rho_\delta \ast |\, \cdot  \,|^s
\, \big)(\phi)
\qquad (\phi \in \mathbb{R}),$

\medskip

\noindent
{\rm (5)}\ 
$\displaystyle{ \frac{\mathrm{d}}{\mathrm{d}\phi}\, 
\big( \, \rho_\delta \ast |\, \cdot  \,|^s\,\mathrm{sgn} 
\, \big)(\phi)}
=s\, \big( \, \rho_\delta \ast |\, \cdot  \,|^{s-1}
\, \big)(\phi)
\ge 0
\qquad (\phi \in \mathbb{R}),$

\medskip

\noindent
{\rm (6)}\ 
$\displaystyle{ \frac{\mathrm{d}}{\mathrm{d}\phi}\, 
\big( \, \rho_\delta \ast |\, \cdot  \,|^{s+1}
\, \big)(\phi)}
=(s+1)\, \big( \, \rho_\delta \ast |\, \cdot  \,|^{s}\,\mathrm{sgn} 
\, \big)(\phi)
\qquad (\phi \in \mathbb{R}),$

\medskip

\noindent
where 
\begin{align*}
\begin{aligned}
\big( \, \rho_\delta \ast |\, \cdot  \,|^s\,\mathrm{sgn} 
\, \big)(\phi)
&=\big( \, \rho_\delta \ast |\, \cdot  \,|^{s-1} \, \cdot 
  \, \big)(\phi)\\
&:=\int^{\infty}_{-\infty}
\rho_\delta(\phi-\psi)\, |\, \psi \,|^{s-1} \, \psi \, \mathrm{d}\psi
\qquad (\phi \in \mathbb{R}),
\end{aligned}
\end{align*}
\begin{align*}
\begin{aligned}
\big( \, \rho_\delta \ast |\, \cdot  \,|^{s+1}
  \, \big)(\phi):=\int^{\infty}_{-\infty}
\rho_\delta(\phi-\psi)\, |\, \psi \,|^{s+1} \, \mathrm{d}\psi
\qquad (\phi \in \mathbb{R}).
\end{aligned}
\end{align*}
}

\noindent
{\bf Proof of Lemma 7.3.}\quad
Because the others are standard or easier than the proof of (5), 
we only show (5). 
Noting the definition of 
$\big( \, \rho_\delta \ast |\, \cdot  \,|^s\,\mathrm{sgn} 
\, \big)(\phi)$, we have 
\begin{align}
\begin{aligned}
&\displaystyle{ \frac{\mathrm{d}}{\mathrm{d}\phi}}\, 
 \big( \, \rho_\delta \ast |\, \cdot  \,|^s\,\mathrm{sgn} 
 \, \big)(\phi)\\
&
=\displaystyle{ \frac{\mathrm{d}}{\mathrm{d}\phi}}\, 
  \left(\, \int^{\infty}_{0}- \int^{0}_{-\infty} \, \right) \, 
\rho_\delta(\phi-\psi)\, |\, \psi \,|^{s}  \, \mathrm{d}\psi
\qquad (\phi \in \mathbb{R}).
\end{aligned}
\end{align}
The direct calculation shows 
\begin{align}
\begin{aligned}
\displaystyle{ \frac{\mathrm{d}}{\mathrm{d}\phi}}\, 
 \int^{\infty}_{0} 
\rho_\delta(\phi-\psi)\, |\, \psi \,|^{s} \, \mathrm{d}\psi
&=- \int^{\infty}_{0} 
  \phi^s \, 
  \partial_{\psi}\big( \, 
  \rho_\delta (\phi-\psi) 
  \, \big)
  \, \mathrm{d}\psi\\
&=s \, \int^{\infty}_{0} 
  \rho_\delta (\phi-\psi) \,|\, \psi \,|^{s-1}\, \mathrm{d}\psi, 
\end{aligned}
\end{align}
\begin{align}
\begin{aligned}
\displaystyle{ - \frac{\mathrm{d}}{\mathrm{d}\phi}}\, 
 \int^{0}_{-\infty} 
\rho_\delta(\phi-\psi)\, |\, \psi \,|^{s} \, \mathrm{d}\psi
&=\int^{0}_{\infty} 
  \varphi^s \, 
  \partial_{\varphi}\big( \, 
  \rho_\delta (\phi+\varphi) 
  \, \big)
  \, \mathrm{d}\varphi\\
&=s \, \int^{0}_{-\infty}
  \rho_\delta (\phi-\psi) \,|\, \psi \,|^{s-1}\, \mathrm{d}\psi. 
\end{aligned}
\end{align}
Substituting (7.8) and (7.9) into (7.7), 
we then have (5). 

This ends the proof of Lemma 7.3. 
\medskip

\noindent
Making use of Lemma 7.3, 
we obtain the following time weighted $L^q$-energy estimate to $\phi$ 
with $1<q<2$. 

\medskip

\noindent
{\bf Proposition 7.4.}\quad {\it
Suppose the same assumptions as in Theorem 3.1. 
For any $q \in ( 1, \, 2 )$, 
there exist positive constants $\alpha$ 
and $C_{\alpha, \phi_{0}, \epsilon}$ 
such that the unique global in time solution $\phi$ 
of the Cauchy problem {\rm (3.4)} has the following $L^q$-energy estimate 
\begin{align*}
\begin{aligned}
&(1+t)^{\alpha}\, \| \, \phi(t) \, \|_{L^q}^q \\
& \quad 
\le C_{\phi_{0}}\, \| \,\phi_0 \,\|_{L^q}^q 
+ C_{\alpha, \phi_{0}}\, (1+t)^{\alpha- \frac{q-1}{2}}
    \, \big( \, \max \big\{ \, 1, \, \ln (1+t) \, \big\} \, \big)^q
\quad \bigl( \, t \ge 0 \, \bigr). 
\end{aligned}
\end{align*}
}

\medskip

The proof of Proposition 7.4 is provided by 
Lemmas 7.1 and 7.3 and the following two Lemmas. 

\medskip

\noindent
{\bf Lemma 7.5.}\quad {\it
For any $q \in [\, 1, \, \infty )$, 
there exist positive constants $\alpha$ and $C_{\phi_{0}}$ such that 
\begin{align}
\begin{aligned}
&(1+t)^{\alpha} \, 
\int_{-\infty}^{\infty} \int_{0}^{\phi} 
\big( \, \rho_\delta \ast |\, \cdot  \,|^{q-1}\,\mathrm{sgn} 
\, \big)(\eta)
\, \mathrm{d}\eta \, \mathrm{d}x \\
&+\int_0^t (1+\tau)^{\alpha} \, \int ^{\infty }_{-\infty } 
      \bigg|\, 
      \int^{| \phi |}_0 
      \eta \, \big( \, \rho_\delta \ast |\, \cdot  \,|^{q-2}\,\mathrm{sgn} 
      \, \big) (\eta ) \, 
      \mathrm{d}\eta 
      \, \bigg|
      \, \partial _x U^r 
      \, \mathrm{d}x \mathrm{d}\tau \\
&+\int_0^t (1+\tau)^{\alpha} \, \int_{-\infty}^{\infty} 
 \big( \, \rho_\delta \ast |\, \cdot  \,|^{q-2}
\, \big)(\phi) \, |\, \partial_{x}\phi \,|^2
 \, \mathrm{d}x \mathrm{d}\tau \\
&\le C_{\phi_{0}} \, 
 \int_{-\infty}^{\infty} \int_{0}^{\phi_{0}} 
\big( \, \rho_\delta \ast |\, \cdot  \,|^{q-1}\,\mathrm{sgn} 
\, \big)(\eta)
\, \mathrm{d}\eta \, \mathrm{d}x\\
&\quad 
 + C_{\phi_{0}} \cdot \alpha \,  
      \int ^t_0 (1+\tau )^{\alpha -1} \, 
      \big\| \,\big( \, \rho_\delta \ast |\, \cdot  \,|^{q-1}\,\mathrm{sgn} 
      \, \big)(\phi)(\tau ) \,\big\|_{L^{q}}^{q} 
      \, \mathrm{d}\tau\\
&\quad + C_{\phi_{0}}\, 
      \int ^t_0 (1+\tau )^{\alpha} \,
      \big\| \,\big( \, \rho_\delta \ast |\, \cdot  \,|^{q-1}\,\mathrm{sgn} 
      \, \big)(\phi)(\tau ) \,\big\|_{L^{\infty}} \, 
      \| \, \partial_{x}^2U^r(\tau ) \,\|_{L^{1}}
      \, \mathrm{d}\tau 
   \quad \big( \, t \ge 0 \, \big). 
\end{aligned}
\end{align}
}

\bigskip

\noindent
{\bf Lemma 7.6.}\quad {\it
Assume $1< q < 2$ and $r>1$. 
We have the following interpolation inequalities.  

 \noindent
 {\rm (1)}\ \ It follows that 
 \begin{align*}
 \begin{aligned}
 &\big\| \,\big( \, \rho_\delta \ast |\, \cdot  \,|\,\mathrm{sgn} 
      \, \big)(\phi)(\tau ) \,\big\|_{L^{\infty}}\\
 &\leq 
  \left( \, 
  \frac{q+1}{2}
   \, \right)^{\frac{2}{q+1}} \, 
   \left( \, \int _{-\infty}^{\infty} 
   \big( \, \rho_\delta \ast |\, \cdot  \,|
      \, \big)(\phi)
   \, \mathrm{d}x \, \right)
  ^{\frac{1}{q+1}} \\
 & \quad 
         \times 
         \left( \, \int _{-\infty}^{\infty} 
         \big( \, \rho_\delta \ast |\, \cdot  \,|^{q-2}
         \, \big)(\phi)
         \, | \, \partial_{x}\phi \, |^{2} 
         \, \mathrm{d}x \, \right)
         ^{\frac{1}{q+1}} 
         +R_{1}^{\frac{2}{q+1}}
         \quad \big( \, t \ge 0 \, \big). 
 \end{aligned}
 \end{align*}

 \noindent
 {\rm (2)}\ \ 
 It follows that 
 \begin{align*}
 \begin{aligned}
 &\big\| \,\big( \, \rho_\delta \ast |\, \cdot  \,|\,\mathrm{sgn} 
      \, \big)(\phi)(\tau ) \,\big\|_{L^{r} }\\
 &\leq 
  \left( \, 
  \frac{q+1}{2}
   \, \right)^{\frac{2(r-1)}{(q+1)r}} \, 
   \left( \, \int _{-\infty}^{\infty} 
   \big( \, \rho_\delta \ast |\, \cdot  \,|
      \, \big)(\phi)
   \, \mathrm{d}x \, \right)
  ^{\frac{r-1}{(q+1)r}} 
  \left( \, \int _{-\infty}^{\infty} 
   \big| \, 
   \big( \, \rho_\delta \ast |\, \cdot  \,|\,\mathrm{sgn} 
      \, \big)(\phi) 
   \,\big|
   \, \mathrm{d}x \, \right)
  ^{\frac{1}{r}}\\
 & \quad 
         \times 
         \left( \, \int _{-\infty}^{\infty} 
         \big( \, \rho_\delta \ast |\, \cdot  \,|^{q-2}
         \, \big)(\phi)
         \, | \, \partial_{x}\phi \, |^{2} 
         \, \mathrm{d}x \, \right)
         ^{\frac{r-1}{(q+1)r}}
         +R_{1}^{\frac{2(r-1)}{(q+1)r}}
         \quad \big( \, t \ge 0 \, \big). 
 \end{aligned}
 \end{align*}
Here, 
\begin{align*}
\begin{aligned}
R_{1}:=\Big| \, 
       \big| \, 
       \big( \, \rho_\delta \ast |\, \cdot  \,|^{\frac{q+1}{2}}
       \,\mathrm{sgn} 
       \, \big)(\phi) 
       \, \big|
       -\big| \, 
       \big( \, \rho_\delta \ast |\, \cdot  \,|
       \,\mathrm{sgn} 
       \, \big)(\phi) 
       \, \big|^{\frac{q+1}{2}}
       \, \Big|. 
\end{aligned}
\end{align*}
}

\medskip

In what follows, we show Lemmas 7.5 and 7.6, 
and finally give the proof of Proposition 7.4 
with $L^{\infty}$-estimate of $\phi$. 

\medskip

{\bf Proof of Lemma 7.5}.
Noting Lemma 7.3, performing the computation 
$$
\int_0^t (1+\tau)^{\alpha} \, \int_{-\infty}^{\infty}
(3.4)\times \big( \, \rho_\delta \ast |\, \cdot  \,|^{q-1}
       \,\mathrm{sgn} 
       \, \big)(\phi) 
\, \mathrm{d}x \mathrm{d}\tau 
\quad \big( \, \alpha>0, \, q\in (1, \, 2)  \, \big)
$$
and using integration by parts to the resultant formula 
allow us to get 
\begin{align}
\begin{aligned}
&(1+t)^{\alpha} \, 
\int_{-\infty}^{\infty} \int_{0}^{\phi(t)} 
\big( \, \rho_\delta \ast |\, \cdot  \,|^{q-1}\,\mathrm{sgn} 
\, \big)(\eta)
\, \mathrm{d}\eta \, \mathrm{d}x \\
&+(q-1)\, \int_0^t (1+\tau)^{\alpha} \, \int ^{\infty }_{-\infty } 
      \int^{\phi}_0 
      \bigl( \, f'(\eta+U^r)-f'(U^r) \, \bigr) \, \\
& \qquad \qquad \qquad \qquad \qquad \times 
      \big( \, \rho_\delta \ast |\, \cdot  \,|^{q-2}
      \, \big) (\eta ) \, \partial _x U^r 
      \, \mathrm{d}\eta 
      \, \mathrm{d}x \mathrm{d}\tau \\
&+(q-1)\, \int_0^t (1+\tau)^{\alpha} \, \int_{-\infty}^{\infty} 
 \big( \, \rho_\delta \ast |\, \cdot  \,|^{q-2}
\, \big)(\phi) \, \partial_{x}\phi \, \\
& \qquad \qquad \qquad \qquad \qquad \times 
\Big( \, \sigma\bigl(\, \partial _x \phi + \partial _x U^r \, \bigr)
    -\sigma\bigl(\, \partial _x U^r \, \bigr) 
    \, \Big)
 \, \mathrm{d}x \mathrm{d}\tau \\
&= \int_{-\infty}^{\infty} \int_{0}^{\phi_{0}} 
\big( \, \rho_\delta \ast |\, \cdot  \,|^{q-1}\,\mathrm{sgn} 
\, \big)(\eta)
\, \mathrm{d}\eta \, \mathrm{d}x\\
&\quad 
 + \alpha \,  
      \int ^t_0 (1+\tau )^{\alpha -1} \, 
      \int ^{\infty }_{-\infty } 
      \int^{\phi}_0 
      \big( \, \rho_\delta \ast |\, \cdot  \,|^{q-1}\,\mathrm{sgn} 
      \, \big)(\eta) \, \mathrm{d}\eta
      \, \mathrm{d}x\mathrm{d}\tau\\
&\quad + 
      \int ^t_0 (1+\tau )^{\alpha} \,
      \int ^{\infty }_{-\infty } 
      \big( \, \rho_\delta \ast |\, \cdot  \,|^{q-1}\,\mathrm{sgn} 
      \, \big)(\eta)\, 
      \partial_{x}
      \Big( \, \sigma\bigl(\, \partial _x U^r \, \bigr) 
      \, \Big)
      \, \mathrm{d}x\mathrm{d}\tau 
   \quad \big( \, t \ge 0 \, \big). 
\end{aligned}
\end{align}
By using the uniform boundedness in Lemmas 4.1 and 5.3,
we estimate the second and third terms 
on the left-hand side of (7.11) as follows. 
\begin{align}
\begin{aligned}
&(q-1)\, \int_0^t (1+\tau)^{\alpha} \, \int ^{\infty }_{-\infty } 
      \int^{\phi}_0 
      \bigl( \, f'(\eta+U^r)-f'(U^r) \, \bigr) \, \\
& \qquad \qquad \qquad \qquad \quad \times 
      \big( \, \rho_\delta \ast |\, \cdot  \,|^{q-2}
      \, \big) (\eta ) \, \partial _x U^r 
      \, \mathrm{d}\eta 
      \, \mathrm{d}x \mathrm{d}\tau \\
&\geq C_{\phi_{0}}^{-1} \, 
      \int _0^t (1+\tau)^{\alpha} \, \int ^{\infty }_{-\infty } 
      \bigg|\, 
      \int^{| \phi |}_0 
      \eta \, \big( \, \rho_\delta \ast |\, \cdot  \,|^{q-2}
      \, \big) (\eta ) \, 
      \mathrm{d}\eta 
      \, \bigg|
      \, \partial _x U^r 
      \, \mathrm{d}x \mathrm{d}\tau, 
\end{aligned}
\end{align}
\begin{align}
\begin{aligned}
&(q-1)\, \int_0^t (1+\tau)^{\alpha} \, \int_{-\infty}^{\infty} 
 \big( \, \rho_\delta \ast |\, \cdot  \,|^{q-2}
\, \big)(\phi) \, \partial_{x}\phi \, \\
& \qquad \qquad \qquad \qquad \quad \times 
\Big( \, \sigma\bigl(\, \partial _x \phi + \partial _x U^r \, \bigr)
    -\sigma\bigl(\, \partial _x U^r \, \bigr) 
    \, \Big)
 \, \mathrm{d}x \mathrm{d}\tau \\
&\geq C_{\phi_{0}}^{-1} \, 
      \int _0^t (1+\tau)^{\alpha} \, \int ^{\infty }_{-\infty } 
      \big( \, \rho_\delta \ast |\, \cdot  \,|^{q-2}
\, \big)(\phi) \, | \, \partial_{x}\phi \, |^2 
      \, \mathrm{d}x \mathrm{d}\tau. 
\end{aligned}
\end{align}
Noting Lemma 7.3, 
the second and third terms 
on the right-hand side of (7.11) are also estimated as 
\begin{align}
\begin{aligned}
&\alpha \,  
      \int ^t_0 (1+\tau )^{\alpha -1} \, 
      \Big| \, 
      \int ^{\infty }_{-\infty } 
      \int^{\phi}_0 
      \big( \, \rho_\delta \ast |\, \cdot  \,|^{q-1}\,\mathrm{sgn} 
      \, \big)(\eta) \, \mathrm{d}\eta
      \, \mathrm{d}x \, \Big| 
      \, \mathrm{d}\tau\\
& \leq \frac{\alpha}{q}\, 
       \int ^t_0 (1+\tau )^{\alpha -1} \, 
       \big\| \,\big( \, \rho_\delta \ast |\, \cdot  \,|\,\mathrm{sgn} 
      \, \big)(\phi)(\tau ) \,\big\|_{L^{q}}^{q}
      \, \mathrm{d}\tau, 
\end{aligned}
\end{align}
\begin{align}
\begin{aligned}
&\int ^t_0 (1+\tau )^{\alpha} \,
      \Big| \, 
      \int ^{\infty }_{-\infty } 
      \big( \, \rho_\delta \ast |\, \cdot  \,|^{q-1}\,\mathrm{sgn} 
      \, \big)(\eta)\, 
      \partial_{x}
      \Big( \, \sigma\bigl(\, \partial _x U^r \, \bigr) 
      \, \Big)
      \, \mathrm{d}x 
      \, \Big| 
      \, \mathrm{d}\tau\\
& \leq C \, \int ^t_0 (1+\tau )^{\alpha} \,
      \big\| \,\big( \, \rho_\delta \ast |\, \cdot  \,|^{q-1}\,\mathrm{sgn} 
      \, \big)(\phi)(\tau ) \,\big\|_{L^{\infty}} \, 
      \| \, \partial _x^2 U^r(\tau) \, \|_{L^{1}}
      \, \mathrm{d}\tau.\
\end{aligned}
\end{align}
Substituting (7.12)-(7.15) into (7.11), 
we immediately have Lemma 7.5. 

\medskip
Next, we give the proof of Lemma 7.6. 
\medskip

{\bf Proof of Lemma 7.6.} 
Noting Lemma 7.3, we have for $s\ge1$, 
\begin{align}
\begin{aligned}
&\big| \,\big( \, \rho_\delta \ast |\, \cdot  \,|^s\,\mathrm{sgn} 
      \, \big)(\phi) \,\big|\\
&\leq s \, \int _{-\infty}^{\infty} \int _{-\infty}^{\infty} 
      \rho_{\delta}(\phi-\psi)  \, | \, \psi \, |^{s-1}
      \, \mathrm{d}\psi \, 
      | \, \partial _x \phi \, |
      \, \mathrm{d}x. 
\end{aligned}
\end{align}
By the Cauchy-Schwarz inequality, we have 
\begin{align}
\begin{aligned}
&\big| \,\big( \, \rho_\delta \ast |\, \cdot  \,|^s\,\mathrm{sgn} 
      \, \big)(\phi) \,\big|
\leq s \, 
       \left( \,
       \int _{-\infty}^{\infty} \int _{-\infty}^{\infty}
       \rho_{\delta}(\phi-\psi)  \, | \, \psi \, |^{2s-q}
       \, \mathrm{d}\psi
       \, \mathrm{d}x \, 
       \right)^{\frac{1}{2}} \\
& \qquad \qquad \qquad \quad \times 
       \left( \, 
       \int _{-\infty}^{\infty} \int _{-\infty}^{\infty} 
       \rho_{\delta}(\phi-\psi)  \, | \, \psi \, |^{q-2}\, 
       | \, \partial _x \phi \,|^{2}
       \, \mathrm{d}\psi
       \, \mathrm{d}x \, 
       \right)^{\frac{1}{2}}.  
\end{aligned}
\end{align}
We also note 
\begin{align}
\begin{aligned}
&\big| \,\big( \, \rho_\delta \ast |\, \cdot  \,|^s\,\mathrm{sgn} 
      \, \big)(\phi) \,\big|
\geq \big| \,\big( \, \rho_\delta \ast |\, \cdot  \,|\,\mathrm{sgn} 
      \, \big)(\phi) \,\big|^s\\
& \qquad \qquad \qquad \quad 
      - \Big| \, 
        \big| \,\big( \, \rho_\delta \ast |\, \cdot  \,|^s\,\mathrm{sgn} 
        \, \big)(\phi) \,\big|
        - \big| \,\big( \, \rho_\delta \ast |\, \cdot  \,|\,\mathrm{sgn} 
        \, \big)(\phi) \,\big|^s
        \, \Big|. 
\end{aligned}
\end{align}
Combining (7.17) and (7.18), we have 
 for $s\ge1$, 
\begin{align}
\begin{aligned}
&\big\| \,\big( \, \rho_\delta \ast |\, \cdot  \,|\,\mathrm{sgn} 
      \, \big)(\phi)(t) \,\big\|_{L^s}
\leq s^{\frac{1}{s}} \, 
       \left( \,
       \int _{-\infty}^{\infty} 
       \big( \, \rho_\delta \ast |\, \cdot  \,|^{2s-q}
      \, \big)(\phi)
       \, \mathrm{d}x \, 
       \right)^{\frac{1}{2s}} \\
& \qquad \qquad \qquad \quad \times 
       \left( \, 
       \int _{-\infty}^{\infty} 
       \big( \, \rho_\delta \ast |\, \cdot  \,|^{q-2}
      \, \big)(\phi) \, | \, \partial _x \phi \,|^{2}
       \, \mathrm{d}x \, 
       \right)^{\frac{1}{2s}}\\
& \qquad \qquad \qquad \quad 
      + \Big| \, 
        \big| \,\big( \, \rho_\delta \ast |\, \cdot  \,|^s\,\mathrm{sgn} 
        \, \big)(\phi) \,\big|
        - \big| \,\big( \, \rho_\delta \ast |\, \cdot  \,|\,\mathrm{sgn} 
        \, \big)(\phi) \,\big|^s
        \, \Big|^{\frac{1}{s}}.
\end{aligned}
\end{align}
Taking $s=\frac{q+1}{2}$, we immediately have (1). 
Next, substituting (1) into 
\begin{align}
\begin{aligned}
&\big\| \,\big( \, \rho_\delta \ast |\, \cdot  \,|\,\mathrm{sgn} 
      \, \big)(\phi)(t) \,\big\|_{L^{r} }^{r}
 \le \big\| \,\big( \, \rho_\delta \ast |\, \cdot  \,|\,\mathrm{sgn} 
      \, \big)(\phi)(t) \,\big\|_{L^\infty }^{r-1}\\
& \qquad \qquad \qquad \quad \times 
      \big\| \,\big( \, \rho_\delta \ast |\, \cdot  \,|\,\mathrm{sgn} 
      \, \big)(\phi)(t) \,\big\|_{L^1 }, 
\end{aligned}
\end{align}
we also have (2). 

Thus, we complete the proof of Lemma 7.6. 

\medskip
We now show Proposition 7.4 with the help of Lemmas 7.5 and 7.6. 
\medskip

{\bf Proof of Proposition 7.4}.\
By using Lemma 7.6, 
we estimate the second term 
on the right-hand side of (7.10) as 
\begin{align}
\begin{aligned}
&C_{\phi_{0}} \cdot \alpha \,  
      \int ^t_0 (1+\tau )^{\alpha -1} \, 
      \big\| \,\big( \, \rho_\delta \ast |\, \cdot  \,|^{q-1}\,\mathrm{sgn} 
      \, \big)(\phi)(\tau ) \,\big\|_{L^{q}}^{q} 
      \, \mathrm{d}\tau\\
&\leq C_{\alpha, \phi_{0}} \,  
      \int ^t_0 (1+\tau )^{\alpha -1} \, 
      \big\| \,\big( \, \rho_\delta \ast |\, \cdot  \,|\,\mathrm{sgn} 
      \, \big)(\phi)(\tau ) \,\big\|_{L^{\infty}}^{q-1}\\
& \qquad \qquad \qquad \qquad \quad \times 
      \big\| \,\big( \, \rho_\delta \ast |\, \cdot  \,|\,\mathrm{sgn} 
      \, \big)(\phi)(\tau ) \,\big\|_{L^{1}}
      \, \mathrm{d}\tau\\
& \quad + C_{\alpha, \phi_{0}} \,  
      \int ^t_0 (1+\tau )^{\alpha -1} \, R_{2} \, \mathrm{d}\tau\\
&\leq C_{\alpha, \phi_{0}} \,  
      \int ^t_0 (1+\tau )^{\alpha -1} \, 
      \left( \, \int _{-\infty}^{\infty} 
   \big( \, \rho_\delta \ast |\, \cdot  \,|
      \, \big)(\phi)
   \, \mathrm{d}x \, \right)
  ^{\frac{q-1}{q+1}}\\
& \quad \times 
  \left( \, \int _{-\infty}^{\infty} 
   \big| \, 
   \big( \, \rho_\delta \ast |\, \cdot  \,|\,\mathrm{sgn} 
      \, \big)(\phi) 
   \,\big|
   \, \mathrm{d}x \, \right)\\
 & \quad 
         \times 
         \left( \, \int _{-\infty}^{\infty} 
         \big( \, \rho_\delta \ast |\, \cdot  \,|^{q-2}
         \, \big)(\phi)
         \, | \, \partial_{x}\phi \, |^{2} 
         \, \mathrm{d}x \, \right)
         ^{\frac{q-1}{q+1}}
      \, \mathrm{d}\tau\\
& \quad + C_{\alpha, \phi_{0}} \,  
      \int ^t_0 (1+\tau )^{\alpha -1} \, 
      \Big( \, 
      R_{1}^{\frac{2(q-1)}{q+1}} \, 
      \big\| \,\big( \, \rho_\delta \ast |\, \cdot  \,|\,\mathrm{sgn} 
      \, \big)(\phi)(\tau ) \,\big\|_{L^{1}}+R_{2}
      \, \Big)
      \, \mathrm{d}\tau\\
&\leq \epsilon \, 
      \int ^t_0 (1+\tau )^{\alpha} \, 
      \int _{-\infty}^{\infty} 
         \big( \, \rho_\delta \ast |\, \cdot  \,|^{q-2}
         \, \big)(\phi)
         \, | \, \partial_{x}\phi \, |^{2} 
         \, \mathrm{d}x\mathrm{d}\tau\\
& \quad + C_{\alpha, \phi_{0}, \epsilon} \,  
       \int ^t_0 (1+\tau )^{\alpha -\frac{q+1}{2}} \, 
       \big\| \,\big( \, \rho_\delta \ast |\, \cdot  \,|
      \, \big)(\phi)(\tau ) \,\big\|_{L^{1}}^q
      \, \mathrm{d}\tau\\
& \quad + C_{\alpha, \phi_{0}} \,  
      \int ^t_0 (1+\tau )^{\alpha -1} \, 
      \Big( \, 
      R_{1}^{\frac{2(q-1)}{q+1}} \, 
      \big\| \,\big( \, \rho_\delta \ast |\, \cdot  \,|\,\mathrm{sgn} 
      \, \big)(\phi)(\tau ) \,\big\|_{L^{1}}+R_{2}
      \, \Big)
      \, \mathrm{d}\tau \quad (\epsilon>0), 
\end{aligned}
\end{align}
where 
\begin{equation*}
R_{2}:= \bigg| \, 
        \int _{-\infty}^{\infty} \int _{0}^{\phi} 
        \big( \, \rho_\delta \ast |\, \cdot  \,|^{q-1}\,\mathrm{sgn} 
      \, \big)(\eta) 
      \, \mathrm{d}\eta \, \mathrm{d}x
      - \frac{1}{q} \, 
        \big\| \,\big( \, \rho_\delta \ast |\, \cdot  \,|\,\mathrm{sgn} 
      \, \big)(\phi)\,\big\|_{L^{q}}^q
        \, \bigg|. 
\end{equation*}
Similarly, the third term 
on the right-hand side of (7.10) is estimated as 
\begin{align}
\begin{aligned}
&C_{\phi_{0}}\, 
      \int ^t_0 (1+\tau )^{\alpha} \,
      \big\| \,\big( \, \rho_\delta \ast |\, \cdot  \,|^{q-1}\,\mathrm{sgn} 
      \, \big)(\phi)(\tau ) \,\big\|_{L^{\infty}} \, 
      \| \, \partial_{x}^2U^r(\tau ) \,\|_{L^{1}} 
      \, \mathrm{d}\tau\\
&\leq C_{\phi_{0}}\, 
      \int ^t_0 (1+\tau )^{\alpha-1} \,
      \big\| \,\big( \, \rho_\delta \ast |\, \cdot  \,|\,\mathrm{sgn} 
      \, \big)(\phi)(\tau ) \,\big\|_{L^{\infty}}^{q-1} 
      \, \mathrm{d}\tau\\
& \quad + C_{\alpha, \phi_{0}} \,  
      \int ^t_0 (1+\tau )^{\alpha -1} \, R_{3} \, \mathrm{d}\tau\\
&\leq \epsilon \, 
      \int ^t_0 (1+\tau )^{\alpha} \, 
      \int _{-\infty}^{\infty} 
         \big( \, \rho_\delta \ast |\, \cdot  \,|^{q-2}
         \, \big)(\phi)
         \, | \, \partial_{x}\phi \, |^{2} 
         \, \mathrm{d}x\mathrm{d}\tau\\
& \quad + C_{\phi_{0}, \epsilon} \,  
       \int ^t_0 (1+\tau )^{\alpha -\frac{q+1}{2}} \, 
       \big\| \,\big( \, \rho_\delta \ast |\, \cdot  \,|
      \, \big)(\phi)(\tau ) \,\big\|_{L^{1}}^{\frac{q-1}{2}}
      \, \mathrm{d}\tau\\
& \quad + C_{\phi_{0}} \,  
      \int ^t_0 (1+\tau )^{\alpha -1} \, 
      \Big( \, R_{1}^{\frac{2(q-1)}{q+1}} + R_{3} \, \Big)
      \, \mathrm{d}\tau \quad (\epsilon>0),
\end{aligned}
\end{align}
where 
\begin{equation*}
R_{3}:= \Big| \, 
        \big\| \,\big( \, \rho_\delta \ast |\, \cdot  \,|^{q-1}\,\mathrm{sgn} 
      \, \big)(\phi)\,\big\|_{L^{\infty}}
      - \big\| \,\big( \, \rho_\delta \ast |\, \cdot  \,|\,\mathrm{sgn} 
      \, \big)(\phi)\,\big\|_{L^{\infty}}^{q-1}
        \, \Big|. 
\end{equation*}
Substituting (7.21)-(7.22) into (7.10) and choosing $\epsilon$ suitably small, 
we then obtain 
\begin{align}
\begin{aligned}
&(1+t)^{\alpha} \, 
\int_{-\infty}^{\infty} \int_{0}^{\phi(t)} 
\big( \, \rho_\delta \ast |\, \cdot  \,|^{q-1}\,\mathrm{sgn} 
\, \big)(\eta)
\, \mathrm{d}\eta \, \mathrm{d}x \\
&+\int_0^t (1+\tau)^{\alpha} \, \int ^{\infty }_{-\infty } 
      \bigg|\, 
      \int^{| \phi |}_0 
      \eta \, \big( \, \rho_\delta \ast |\, \cdot  \,|^{q-2}\,\mathrm{sgn} 
      \, \big) (\eta ) \, 
      \mathrm{d}\eta 
      \, \bigg|
      \, \partial _x U^r 
      \, \mathrm{d}x \mathrm{d}\tau \\
&+\int_0^t (1+\tau)^{\alpha} \, \int_{-\infty}^{\infty} 
 \big( \, \rho_\delta \ast |\, \cdot  \,|^{q-2}
\, \big)(\phi) \, |\, \partial_{x}\phi \,|^2
 \, \mathrm{d}x \mathrm{d}\tau \\
&\le C_{\phi_{0}} \, 
 \int_{-\infty}^{\infty} \int_{0}^{\phi_{0}} 
\big( \, \rho_\delta \ast |\, \cdot  \,|^{q-1}\,\mathrm{sgn} 
\, \big)(\eta)
\, \mathrm{d}\eta \, \mathrm{d}x\\
&\quad 
 + C_{\alpha, \phi_{0}} \,  
       \int ^t_0 (1+\tau )^{\alpha -\frac{q+1}{2}} \\
&\qquad \times
       \Big( \, 
       \big\| \,\big( \, \rho_\delta \ast |\, \cdot  \,|
      \, \big)(\phi)(\tau ) \,\big\|_{L^{1}}^q
      + \big\| \,\big( \, \rho_\delta \ast |\, \cdot  \,|
      \, \big)(\phi)(\tau ) \,\big\|_{L^{1}}^{\frac{q-1}{2}}
      \, \Big)
      \, \mathrm{d}\tau\\
&\quad + C_{\alpha, \phi_{0}} \,  
      \int ^t_0 (1+\tau )^{\alpha-1} \\
&\qquad \times
      \Big( \, 
      R_{1}^{\frac{2(q-1)}{q+1}} \, 
      \Big( \, 
      \big\| \,\big( \, \rho_\delta \ast |\, \cdot  \,|\,\mathrm{sgn} 
      \, \big)(\phi)(\tau ) \,\big\|_{L^{1}} +1 \, \Big)+R_{2}+R_{3}
      \, \Big)
      \, \mathrm{d}\tau. %
\end{aligned}
\end{align}
Noting Lemmas 7.1 and 7.3, in particular, for $t\geq0$, 
\begin{align}
\begin{aligned}
\big\| \,\big( \, \rho_\delta \ast |\, \cdot  \,|\,\mathrm{sgn} 
      \, \big)(\phi)(t ) \,\big\|_{L^{1}} 
&\leq \big\| \,\big( \, \rho_\delta \ast |\, \cdot  \,|
      \, \big)(\phi)(t ) \,\big\|_{L^{1}} \\
&= \bigg\| \,
   \bigg( \, 
   \int^{\phi}_0\left( \, \rho_\delta \ast \mathrm{sgn} \, \right)(\eta)
   \, \mathrm{d}\eta 
   \, \bigg) 
   (t ) \,\bigg\|_{L^{1}} \\
&\rightarrow \| \, \phi(t) \, \|_{L^{1}} 
\quad (\delta \rightarrow 0),
\end{aligned}
\end{align}
\begin{equation}
\int_{-\infty}^{\infty} \int_{0}^{\phi(t)} 
\big( \, \rho_\delta \ast |\, \cdot  \,|^{q-1}\,\mathrm{sgn} 
\, \big)(\eta)
\, \mathrm{d}\eta \, \mathrm{d}x
\rightarrow \frac{1}{q} \, \| \, \phi(t) \, \|_{L^{q}}^q
\quad (\delta \rightarrow 0),
\end{equation}
\begin{equation}
R_{1}^{\frac{2(q-1)}{q+1}} \, 
      \Big( \, 
      \big\| \,\big( \, \rho_\delta \ast |\, \cdot  \,|\,\mathrm{sgn} 
      \, \big)(\phi)(t ) \,\big\|_{L^{1}} +1 \, \Big)+R_{2}+R_{3}
\rightarrow 0 
\quad (\delta \rightarrow 0), 
\end{equation}
and taking the limit $\delta \rightarrow 0$ 
to (7.23), 
we get 
\begin{align}
\begin{aligned}
&(1+t)^{\alpha}\, \| \, \phi(t) \, \|_{L^q}^q \\
& \quad 
\le C_{\phi_{0}}\, \| \,\phi_0 \,\|_{L^q}^q 
+ C_{\alpha, \phi_{0}}\, (1+t)^{\alpha- \frac{q+1}{2}}
    \, \Big( \, 
    \| \, \phi(t) \, \|_{L^{1}}^q
    + \| \, \phi(t) \, \|_{L^{1}}^{\frac{q-1}{2}}
    \, \Big).
\end{aligned}
\end{align}
Then, by using Proposition 7.2, we obtain the desired 
time-weighted $L^q$-energy estimate of $\partial_{x}\phi$ with $1< q < 2$.

Thus, we complete the proof of Proposition 7.4. 

\medskip
We further establish the $L^q$-estimate of $\phi$ 
in the case where $2\leq q < \infty$ and $q= \infty$. 
With an argument similar to that used to prove Proposition 6.1, 
we can obtain the following 
time-weighted $L^q$-energy estimate with $2\leq q < \infty$
(therefore, the proof can be omitted). 

\medskip

\noindent
{\bf Proposition 7.7.}\quad {\it
Suppose the same assumptions as in Theorem 3.1. 
For any $q \in [\, 2, \, \infty )$, 
there exist positive constants $\alpha$ 
and $C_{\alpha, \phi_{0}, \epsilon}$ 
such that the unique global in time solution $\phi$ 
of the Cauchy problem {\rm (3.4)} has the following $L^q$-energy estimate 
\begin{align*}
\begin{aligned}
&(1+t)^{\alpha}\, \| \, \phi(t) \, \|_{L^q}^q 
+\int_0^t (1+\tau)^{\alpha} \, 
 \int_{-\infty}^{\infty} | \, \phi \, |^{q} 
 \, \partial_{x}U^r \, 
 \mathrm{d}x\mathrm{d}\tau\\
&+\int_0^t (1+\tau)^{\alpha} \, \int_{-\infty}^{\infty} 
 | \, \phi\,|^{q-2} \, | \, \partial_{x}\phi\,|^2
 \, \mathrm{d}x \mathrm{d}\tau \\
&\le C_{\phi_{0}}\, \| \,\phi_0 \,\|_{L^q}^q 
+ C_{\alpha, \phi_{0}}\, (1+t)^{\alpha- \frac{q-1}{2}}
    \, \big( \, \max \big\{ \, 1, \, \ln (1+t) \, \big\} \, \big)^q
\quad \bigl( \, t \ge 0 \, \bigr). 
\end{aligned}
\end{align*}
}

\medskip

\noindent
Therefore, we can conclude from Propositions 7.2, 7.4 and 7.7 
that for $q \in [\,1, \, \infty)$, 
\begin{equation}
\| \, \phi(t) \, \|_{L^q}
\le C_{\phi_{0}}\, (1+t)^{- \frac{1}{2}\left(1-\frac{1}{q} \right)}
    \, \max \big\{ \, 1, \, \ln (1+t) \, \big\}
\quad \bigl( \, t \ge 0 \, \bigr) 
\end{equation}
by taking $\alpha$ suitably large. 
Furthermore, by using that 
$\| \, \partial_{x}\phi(t) \, \|_{L^2} \leq C_{\phi_0}$ 
for some $C_{\phi_0}>0$ 
and the Gagliardo-Nirenberg inequality, we have 
for any $\epsilon>0$, 
\begin{equation}
\| \, \phi(t) \, \|_{L^{\infty}}
\le C_{\phi_{0}, \epsilon}\, (1+t)^{- \frac{1}{2}+\epsilon}
\quad \bigl( \, t \ge 0 \, \bigr). 
\end{equation}

Thus, the proof of Theorem 3.6 is completed. 

\medskip

\noindent
{\it Remark 7.8.}\quad 
If we want to get the time-decay estimates 
of $\| \, \phi(t) \, \|_{L^r}$ with $r \in (1, \, 2)$ 
by a simpler way than that in the proof of Proposition 7.4, 
we use the following way though the time-decay rate in Proposition 7.4 
is a little sharper. 
Substituting the estimate in Proposition 7.2 
and (7.29) into 
\begin{equation}
\| \, \phi(t) \, \|_{L^{r}}^r
\le \| \, \phi(t) \, \|_{L^{\infty}}^{r-1} \, 
    \| \, \phi(t) \, \|_{L^{1}}
\end{equation}
we immediately have for $r \in (1, \, 2)$ and any $\epsilon>0$,
\begin{equation}
\| \, \phi(t) \, \|_{L^r}
\le C_{\phi_{0}, \epsilon}\, 
   (1+t)^{- \frac{1}{2}\left(1-\frac{1}{r} \right)+\epsilon}
\quad \bigl( \, t \ge 0 \, \bigr).
\end{equation}

\medskip

\noindent
{\it Remark 7.9.}\quad 
Although we use the time-decay estimate 
\begin{equation}
\| \, \partial_{x}\phi(t) \, \|_{L^2}
\leq C_{\phi_{0}, \epsilon} \, 
     (1+t)^{-\frac{2}{3}+ \epsilon}
\end{equation}
in Theorem 3.8 
instead of the uniform boundedness of 
$\| \, \partial_{x}\phi(t) \, \|_{L^2}$
to get the better decay rate of $\| \, \phi(t) \, \|_{L^{\infty}}$ 
than in (7.28), 
the decay rate only becomes the same as in (7.28).

\bigskip 

\noindent
\section{Time-decay estimates I\hspace{-.1em}I\hspace{-.1em}I}
In this section (and Sections 8 and 9), 
we try to show the time-decay estimate 
for the derivative $\partial_x \phi$ in Theorem 3.8 
with the help of Theorems 3.4 and 3.7. 
To do that, we establish 
the following time-weighted $L^{q+1}$-energy estimate 
to $\partial_x \phi$ with $1\leq q\leq \infty$. 

\medskip

\noindent
{\bf Proposition 8.1.}\quad {\it
Suppose the same assumptions as in Theorem 3.1. 
For any $q \in [\, 1, \, \infty )$, 
there exist positive constants $\alpha$,  $C_{\alpha, \phi_{0}}$ 
and $C_{\alpha, \phi_{0}, \epsilon}$ 
such that the unique global in time solution $\phi$ 
of the Cauchy problem {\rm (3.4)} has the following $L^{q+1}$-energy estimate 
\begin{align*}
\begin{aligned}
&(1+t)^{\alpha} \, 
\| \, \partial_{x}\phi(t)\,\|_{L^{q+1}}^{q+1}
+\int_0^t (1+\tau)^{\alpha} \, \int_{-\infty}^{\infty} 
 | \, \partial_{x}\phi\,|^{q-1} \, | \, \partial_{x}^2\phi\,|^2
 \, \mathrm{d}x \mathrm{d}\tau \\
&+\int_0^t (1+\tau)^{\alpha} \, \int_{-\infty}^{\infty} 
 | \, \partial_{x}\phi\,|^{q+1} \, \partial_{x}U^r
 \, \mathrm{d}x \mathrm{d}\tau \\
&\le C_{\alpha, \phi_{0}} \, 
 \| \, \partial_{x}\phi_{0}\,\|_{L^{q+1}}^{q+1} 
 + C_{\alpha, \phi_{0}, \epsilon}\, 
   (1+t)^{\alpha-\frac{2q+1}{2}+\epsilon} 
\quad \big(t\ge0\big), 
\end{aligned}
\end{align*}
for any $\epsilon>0$. 
}

\medskip
The proof of Proposition 8.1 is provided by the following two Lemmas. 
\medskip

\noindent
{\bf Lemma 8.2.}\quad {\it
For any $q \in [\, 1, \, \infty )$, 
there exist positive constants $\alpha$ and $C_{\phi_{0}}$ such that 
\begin{align}
\begin{aligned}
&(1+t)^{\alpha} \, 
\| \, \partial_{x}\phi(t)\,\|_{L^{q+1}}^{q+1}
+\int_0^t (1+\tau)^{\alpha} \, \int_{-\infty}^{\infty} 
 | \, \partial_{x}\phi\,|^{q-1} \, | \, \partial_{x}^2\phi\,|^2
 \, \mathrm{d}x \mathrm{d}\tau \\
&+\int_0^t (1+\tau)^{\alpha} \, \int_{-\infty}^{\infty} 
 f''(\phi+U^r)\, | \, \partial_{x}\phi\,|^{q+1} \, \partial_{x}U^r
 \, \mathrm{d}x \mathrm{d}\tau \\
&\le C_{\phi_{0}} \, 
 \| \, \partial_{x}\phi_{0}\,\|_{L^{q+1}}^{q+1} 
 + C_{\phi_{0}} \cdot \alpha \,  
      \int ^t_0 (1+\tau )^{\alpha -1} \, 
      \| \,\partial_{x}\phi(\tau ) \,\|_{L^{q+1}}^{q+1} 
      \, \mathrm{d}\tau\\
&\quad + C_{\phi_{0}}\, 
      \int ^t_0 (1+\tau )^{\alpha} \,
      \| \,\partial_{x}\phi(\tau ) \,\|_{L^{q+2}}^{q+2} 
      \, \mathrm{d}\tau \\
&\quad +  C_{\phi_{0}}\, 
      \int ^t_0 (1+\tau )^{\alpha} \,
      \int_{-\infty}^{\infty} 
 | \, \partial_{x}\phi\,|^{q-1} \, \phi^2 \, | \, \partial_{x}U^r\,|^2
 \, \mathrm{d}x \mathrm{d}\tau \\
&\quad + C_{\phi_{0}}\, 
      \int ^t_0 (1+\tau )^{\alpha} \,
      \int_{-\infty}^{\infty} 
 | \, \partial_{x}\phi\,|^{q-1} \, | \, \partial_{x}^2U^r\,|^2
 \, \mathrm{d}x \mathrm{d}\tau 
   \quad \big( \, t \ge 0 \, \big). 
\end{aligned}
\end{align}
}

\bigskip

\noindent
{\bf Lemma 8.3.}\quad {\it
Assume $1\leq q < \infty$ and $1\leq r < \infty$. 
We have the following interpolation inequalities.  

 \noindent
 {\rm (1)}\ \ It follows that 
 \begin{align*}
 \begin{aligned}
 &\Vert \, \partial_{x}\phi (t) \, \Vert _{L^{\infty} }
  \leq 
  \left( \, 
  \frac{q+3}{2}
   \, \right)^{\frac{2}{q+3}} \, 
   \left( \, \int _{-\infty}^{\infty} | \, \partial_{x}\phi \, |^2 \, \mathrm{d}x \, \right)
  ^{\frac{1}{q+3}} \\
 & \qquad \quad \quad \quad \quad \; \: \, 
         \times 
         \left( \, \int _{-\infty}^{\infty} | \, \partial_{x}\phi \, |^{q-1} 
         | \, \partial_{x}^2\phi \, |^{2} \, \mathrm{d}x \, \right)
         ^{\frac{1}{q+3}} \quad \big( \, t \ge 0 \, \big). 
 \end{aligned}
 \end{align*}

 \noindent
 {\rm (2)}\ \ 
 It follows that 
 \begin{align*}
 \begin{aligned}
 &\Vert \, \partial_{x}\phi (t) \, \Vert _{L^{r+1} }
  \leq 
  \left( \, 
  \frac{q+3}{2}
   \, \right)^{\frac{2(r-1)}{(q+3)(r+1)}} \, 
   \left( \, \int _{-\infty}^{\infty} | \, \partial_{x}\phi \, |^2 \, \mathrm{d}x \, \right)
  ^{\frac{r-1}{(q+3)(r+1)}} \\
 & \qquad \quad \quad \quad \quad \; \: \, 
         \times 
         \left( \, \int _{-\infty}^{\infty} | \, \partial_{x}\phi \, |^{q-1} 
         | \, \partial_{x}^2\phi \, |^{2} \, \mathrm{d}x \, \right)
         ^{\frac{q+r+2}{(q+3)(r+1)}} \quad \big( \, t \ge 0 \, \big). 
 \end{aligned}
 \end{align*}
}

\medskip

In what follows, we show Lemmas 8.2, 8.3, 
and finally give the proof of Proposition 8.1 
with $L^{\infty}$-estimate of $\partial_{x}\phi$.

\medskip

{\bf Proof of Lemma 8.2}.
Performing the computation 
$$
\int_0^t (1+\tau)^{\alpha} \, \int_{-\infty}^{\infty}
(3.4)\times \Big( \, 
- \partial_{x}\big( \,
| \, \partial_{x}\phi\,|^{q-1} \, \partial_{x}\phi \, \big) \, \Big) 
\, \mathrm{d}x \mathrm{d}\tau 
\quad \big( \, \alpha>0, \, q\in [\, 1, \, \infty)  \, \big)
$$
and using integration by parts to the resultant formula, 
we arrive at 
\begin{align}
\begin{aligned}
&\frac{1}{q+1} \, (1+t)^{\alpha} \, 
\| \, \partial_{x}\phi(t)\,\|_{L^{q+1}}^{q+1}\\
&- q \, \int_0^t (1+\tau)^{\alpha} \, \int_{-\infty}^{\infty} 
 | \, \partial_{x}\phi\,|^{q-1} \, \partial_{x}^2\phi \, 
 \partial_{x}\big( \,
f(\phi+U^r)-f(U^r) \, \big)
 \, \mathrm{d}x \mathrm{d}\tau \\
&+q \, \int_0^t (1+\tau)^{\alpha} \, \int_{-\infty}^{\infty} 
 | \, \partial_{x}\phi\,|^{q-1} \, \partial_{x}^2\phi \, 
 \partial_{x}\Big( \,
\sigma \big(\, \partial_{x}\phi+\partial_{x}U^r\, \big)
-\sigma \big(\, \partial_{x}U^r\, \big) \, \Big)
 \, \mathrm{d}x \mathrm{d}\tau \\
&=\frac{1}{q+1} \, 
 \| \, \partial_{x}\phi_{0}\,\|_{L^{q+1}}^{q+1} 
 + \frac{\alpha}{q+1} \, 
      \int ^t_0 (1+\tau )^{\alpha -1} \, 
      \| \,\partial_{x}\phi(\tau ) \,\|_{L^{q+1}}^{q+1} 
      \, \mathrm{d}\tau\\
&\quad + q\, 
      \int ^t_0 (1+\tau )^{\alpha} \, \int_{-\infty}^{\infty} 
 | \, \partial_{x}\phi\,|^{q-1} \, \partial_{x}^2\phi \, 
 \partial_{x}\Big( \,
\sigma \big(\, \partial_{x}U^r\, \big) \, \Big)
 \, \mathrm{d}x \mathrm{d}\tau
   \quad \big( \, t \ge 0 \, \big). 
\end{aligned}
\end{align}
By using integration by parts, the Young inequality, 
and Lemmas 4.1 and 5.3, 
we estimate the second and third terms 
on the left-hand side of (8.2) as follows. 
\begin{align}
\begin{aligned}
&- q \, \int_0^t (1+\tau)^{\alpha} \, \int_{-\infty}^{\infty} 
 | \, \partial_{x}\phi\,|^{q-1} \, \partial_{x}^2\phi \, 
 \partial_{x}\big( \,
f(\phi+U^r)-f(U^r) \, \big)
 \, \mathrm{d}x \mathrm{d}\tau \\
&\geq - q \, \int_0^t (1+\tau)^{\alpha} \, \int_{-\infty}^{\infty} 
 | \, \partial_{x}\phi\,|^{q-1} \, \partial_{x}\phi \, 
 \partial_{x}^2\phi \, f''(\phi + U^r)
 \, \mathrm{d}x \mathrm{d}\tau\\
& \quad 
  - C_{\phi_{0}}\, \int_0^t (1+\tau)^{\alpha} \, \int_{-\infty}^{\infty} 
 | \, \phi\,|\, | \, \partial_{x}\phi\,|^{q-1} \, 
 | \, \partial_{x}^2\phi \,|\,| \, \partial_{x}U^r \,|
 \, \mathrm{d}x \mathrm{d}\tau\\    
&\geq \frac{q}{q+1} \, \, \int_0^t (1+\tau)^{\alpha} \, \int_{-\infty}^{\infty} 
 f''(\phi + U^r) \, | \, \partial_{x}\phi\,|^{q+1} \, 
 \partial_{x}U^r
 \, \mathrm{d}x \mathrm{d}\tau\\
& \quad 
  - \frac{q}{q+1} \, \int_0^t (1+\tau)^{\alpha} \, \int_{-\infty}^{\infty} 
 f''(\phi + U^r) \, | \, \partial_{x}\phi\,|^{q+2} 
 \, \mathrm{d}x \mathrm{d}\tau\\
& \quad 
  - \epsilon \, \int_0^t (1+\tau)^{\alpha} \, \int_{-\infty}^{\infty} 
    | \, \partial_{x}\phi\,|^{q-1} \, 
 | \, \partial_{x}^2\psi \,|^2
 \, \mathrm{d}x \mathrm{d}\tau\\
& \quad 
  - C_{\phi_{0}, \epsilon} \, \int_0^t (1+\tau)^{\alpha} \, \int_{-\infty}^{\infty} 
    | \, \partial_{x}\phi\,|^{q-1} \, 
 \phi^2 \, | \, \partial_{x}U^r \,|^2
 \, \mathrm{d}x \mathrm{d}\tau \quad (\epsilon>0), 
\end{aligned}
\end{align}
\begin{align}
\begin{aligned}
&q \, \int_0^t (1+\tau)^{\alpha} \, \int_{-\infty}^{\infty} 
 | \, \partial_{x}\phi\,|^{q-1} \, \partial_{x}^2\phi \, 
 \partial_{x}\Big( \,
\sigma \big(\, \partial_{x}\phi+\partial_{x}U^r\, \big)
-\sigma \big(\, \partial_{x}U^r\, \big) \, \Big)
 \, \mathrm{d}x \mathrm{d}\tau \\
&\geq C^{-1}\, \int_0^t (1+\tau)^{\alpha} \, \int_{-\infty}^{\infty} 
 | \, \partial_{x}\phi\,|^{q-1} \, \braket{ \, \partial_{x}\phi \, }^{p-1} \, 
 | \, \partial_{x}^2\phi \,|^2
 \, \mathrm{d}x \mathrm{d}\tau\\
& \quad 
  - C_{\phi_{0}}\, \int_0^t (1+\tau)^{\alpha} \, \int_{-\infty}^{\infty} 
 | \, \partial_{x}\phi\,|^{q-1} \, 
 | \, \partial_{x}^2\phi \,|\,| \, \partial_{x}^2U^r \,|
 \, \mathrm{d}x \mathrm{d}\tau\\ 
&\geq (\,C_{\phi_{0}}^{-1}-\epsilon \,)\, \int_0^t (1+\tau)^{\alpha} \, \int_{-\infty}^{\infty} 
 | \, \partial_{x}\phi\,|^{q-1} \, 
 | \, \partial_{x}^2\phi \,|^2
 \, \mathrm{d}x \mathrm{d}\tau\\
& \quad 
  - C_{\phi_{0}, \epsilon}\, \int_0^t (1+\tau)^{\alpha} \, \int_{-\infty}^{\infty} 
 | \, \partial_{x}\phi\,|^{q-1} \, 
 | \, \partial_{x}^2U^r \,|^2
 \, \mathrm{d}x \mathrm{d}\tau \quad (\epsilon>0).
\end{aligned}
\end{align}
Similarly, third term on the right-hand side of (8.2) is estimated as 
\begin{align}
\begin{aligned}
&q\, \int ^t_0 (1+\tau )^{\alpha} \, \bigg| \, \int_{-\infty}^{\infty} 
 | \, \partial_{x}\phi\,|^{q-1} \, \partial_{x}^2\phi \, 
 \partial_{x}\Big( \,
\sigma \big(\, \partial_{x}U^r\, \big) \, \Big)
 \, \mathrm{d}x \, \bigg| \, \mathrm{d}\tau \\
&\leq \epsilon \, \int_0^t (1+\tau)^{\alpha} \, \int_{-\infty}^{\infty} 
    | \, \partial_{x}\phi\,|^{q-1} \, 
 | \, \partial_{x}^2\phi \,|^2
 \, \mathrm{d}x \mathrm{d}\tau\\
& \quad 
  + C_{\epsilon} \, \int_0^t (1+\tau)^{\alpha} \, \int_{-\infty}^{\infty} 
    | \, \partial_{x}\phi\,|^{q-1} \, 
 \phi^2 \, | \, \partial_{x}U^r \,|^2
 \, \mathrm{d}x \mathrm{d}\tau \quad (\epsilon>0).
\end{aligned}
\end{align}
Substituting (8.3)-(8.5) into (8.2) and choosing $\epsilon$ suitably small, 
we immediately have Lemma 8.2. 

\medskip

Next, we show Lemma 8.3. 

\medskip

{\bf Proof of Lemma 8.3}.
Noting that $\phi (t, \cdot \: ) \in H^2$ 
imply 
$\displaystyle{\lim _{x\rightarrow \pm \infty}\partial_{x}\phi(t,x)}=0$ 
for $t \ge 0$ and 
\begin{align}
\begin{aligned}
 | \, \partial_{x}\phi \, |^{s} 
 &\leq s \int _{-\infty}^{\infty} 
      | \, \partial_{x}\phi \, |^{s-1}  \, 
      | \, \partial _x^2 \phi \, |
      \, \mathrm{d}x\quad \big(\, s\ge1 \,\big).\\ 
\end{aligned}
\end{align}
By the Cauchy-Schwarz inequality, we have 
\begin{equation}
 | \, \partial_{x}\phi \, |^{s} 
 \leq s \, 
       \left( \,
       \int _{-\infty}^{\infty} 
       | \, \partial_{x}\phi \, |^{2s-q-1} 
       \, \mathrm{d}x
       \, \right)^{\frac{1}{2}} \, 
       \left( \, 
       \int _{-\infty}^{\infty} 
       | \, \partial_{x}\phi \, |^{q-1} \,
       | \, \partial _x^2 \phi \,|^{2}
       \, \mathrm{d}x 
       \, \right)^{\frac{1}{2}}.  
\end{equation}
Taking $s=\frac{q+3}{2}$, we immediately have (1). 
Next, substituting (1) into 
\begin{equation}
 \Vert \, \partial_{x}\phi \, \Vert _{L^{r+1} }^{r+1}
 \le \Vert \, \partial_{x}\phi \, \Vert _{L^\infty }^{r-1}
      \Vert \, \partial_{x}\phi \, \Vert _{L^2 }^2, 
\end{equation}
we also have (2). 

Thus, we complete the proof of Lemma 8.3. 

\medskip

{\bf Proof of Proposition 8.1}.\
By using Lemma 8.3, 
we estimate the second, third, fourth and fifth terms 
on the right-hand side of (8.1) as follows. 
\begin{align}
\begin{aligned}
&C_{\phi_{0}} \cdot \alpha \,  
      \int ^t_0 (1+\tau )^{\alpha -1} \, 
      \| \,\partial_{x}\phi(\tau ) \,\|_{L^{q+1}}^{q+1} 
      \, \mathrm{d}\tau\\
&\leq C_{\phi_{0}} \cdot \alpha \,  
      \int ^t_0 (1+\tau )^{\alpha -1} \, 
      \left( \, \int _{-\infty}^{\infty} | \, \partial_{x}\phi \, |^2 
      \, \mathrm{d}x \, \right)
      ^{\frac{2q+2}{q+3}}\\
& \quad \times 
         \left( \, \int _{-\infty}^{\infty} | \, \partial_{x}\phi \, |^{q-1} 
         | \, \partial_{x}^2\phi \, |^{2} \, \mathrm{d}x \, \right)
         ^{\frac{q-1}{q+3}}
      \, \mathrm{d}\tau\\
&\leq \epsilon \, \int_0^t (1+\tau)^{\alpha} \, \int_{-\infty}^{\infty} 
    | \, \partial_{x}\phi\,|^{q-1} \, 
 | \, \partial_{x}^2\phi \,|^2
 \, \mathrm{d}x \mathrm{d}\tau\\
& \quad + C_{\alpha, \phi_{0}, \epsilon}\, 
          \, \int_0^t (1+\tau)^{\alpha-\frac{q+3}{4}}\, 
          \left( \, \int _{-\infty}^{\infty} | \, \partial_{x}\phi \, |^{2} 
          \, \mathrm{d}x \, \right)
         ^{\frac{q+1}{2}}
         \, \mathrm{d}\tau \quad (\epsilon>0), 
\end{aligned}
\end{align}
\begin{align}
\begin{aligned}
&C_{\phi_{0}} \,  
      \int ^t_0 (1+\tau )^{\alpha} \, 
      \| \,\partial_{x}\phi(\tau ) \,\|_{L^{q+2}}^{q+2} 
      \, \mathrm{d}\tau\\
&\leq C_{\phi_{0}} \,  
      \int ^t_0 (1+\tau )^{\alpha} \, 
      \left( \, \int _{-\infty}^{\infty} | \, \partial_{x}\phi \, |^2 
      \, \mathrm{d}x \, \right)
      ^{\frac{2q+3}{q+3}}\\
& \quad \times 
         \left( \, \int _{-\infty}^{\infty} | \, \partial_{x}\phi \, |^{q-1} 
         | \, \partial_{x}^2\phi \, |^{2} \, \mathrm{d}x \, \right)
         ^{\frac{q}{q+3}}
      \, \mathrm{d}\tau\\
&\leq \epsilon \, \int_0^t (1+\tau)^{\alpha} \, \int_{-\infty}^{\infty} 
    | \, \partial_{x}\phi\,|^{q-1} \, 
 | \, \partial_{x}^2\phi \,|^2
 \, \mathrm{d}x \mathrm{d}\tau\\
& \quad + C_{\phi_{0}, \epsilon}\, 
          \, \int_0^t (1+\tau)^{\alpha}\, 
          \left( \, \int _{-\infty}^{\infty} | \, \partial_{x}\phi \, |^{2} 
          \, \mathrm{d}x \, \right)
         ^{\frac{2q+3}{3}}
         \, \mathrm{d}\tau \quad (\epsilon>0). 
\end{aligned}
\end{align}
Noting Lemma 2.2 and 
\begin{equation*}
\| \, \phi(t) \, \|_{L^{2}}^2 
\leq C_{\phi_{0}}\, (1+t)^{-\frac{1}{2}}\, 
\big( \, \max \big\{ \, 1, \, \ln (1+t) \, \big\} \, \big)^2
\end{equation*}
from Theorem 3.7, 
fourth and fifth terms 
on the right-hand side of (8.1) are estimated as follows. 
\begin{align}
\begin{aligned}
&C_{\phi_{0}}\, 
      \int ^t_0 (1+\tau )^{\alpha} \,
      \int_{-\infty}^{\infty} 
 | \, \partial_{x}\phi\,|^{q-1} \, \phi^2 \, | \, \partial_{x}U^r\,|^2
 \, \mathrm{d}x \mathrm{d}\tau\\
&\leq C_{\phi_{0}} \,  
      \int ^t_0 (1+\tau )^{\alpha -\frac{5}{2}} 
      \, \big( \, \max \big\{ \, 1, \, \ln (1+\tau) \, \big\} \, \big)^{2}
      \left( \, \int _{-\infty}^{\infty} | \, \partial_{x}\phi \, |^2 
      \,\mathrm{d}x \, \right)
      ^{\frac{q-1}{q+3}}\\
& \quad \times 
         \left( \, \int _{-\infty}^{\infty} | \, \partial_{x}\phi \, |^{q-1} 
         | \, \partial_{x}^2\phi \, |^{2} \, \mathrm{d}x \, \right)
         ^{\frac{q-1}{q+3}}
      \, \mathrm{d}\tau\\
&\leq \epsilon \, \int_0^t (1+\tau)^{\alpha} \, \int_{-\infty}^{\infty} 
    | \, \partial_{x}\phi\,|^{q-1} \, 
 | \, \partial_{x}^2\phi \,|^2
 \, \mathrm{d}x \mathrm{d}\tau\\
& \quad + C_{\phi_{0}, \epsilon}\, 
          \, \int_0^t (1+\tau)^{\alpha-\frac{5(q+3)}{8}}
          \, \big( \, \max \big\{ \, 1, \, \ln (1+\tau) \, \big\} \, \big)^{\frac{q+3}{2}}\\
& \qquad \qquad \qquad \qquad \qquad \qquad \times
          \left( \, \int _{-\infty}^{\infty} | \, \partial_{x}\phi \, |^{2} 
          \, \mathrm{d}x \, \right)
         ^{\frac{q-1}{4}}
         \, \mathrm{d}\tau \quad (\epsilon>0), 
\end{aligned}
\end{align}
\begin{align}
\begin{aligned}
&C_{\phi_{0}}\, 
      \int ^t_0 (1+\tau )^{\alpha} \,
      \int_{-\infty}^{\infty} 
 | \, \partial_{x}\phi\,|^{q-1} \, | \, \partial_{x}^2U^r\,|^2
 \, \mathrm{d}x \mathrm{d}\tau\\
&\leq \epsilon \, \int_0^t (1+\tau)^{\alpha} \, \int_{-\infty}^{\infty} 
    | \, \partial_{x}\phi\,|^{q-1} \, 
 | \, \partial_{x}^2\phi \,|^2
 \, \mathrm{d}x \mathrm{d}\tau\\
& \quad + C_{\phi_{0}, \epsilon}\, 
          \, \int_0^t (1+\tau)^{\alpha-\frac{5(q+3)}{8}}
          \, 
          \left( \, \int _{-\infty}^{\infty} | \, \partial_{x}\phi \, |^{2} 
          \, \mathrm{d}x \, \right)
         ^{\frac{q-1}{4}}
         \, \mathrm{d}\tau \quad (\epsilon>0). 
\end{aligned}
\end{align}
Noting $f''(\phi+U^r)\geq C_{\phi_{0}}^{-1}$ from 
the uniform boundedness of 
$\| \, \phi \, \|_{L^{\infty}}$, 
substituting (8.9)-(8.12) into (8.1) and choosing $\epsilon$ suitably small, 
we obtain 
\begin{align}
\begin{aligned}
&(1+t)^{\alpha} \, 
\| \, \partial_{x}\phi(t)\,\|_{L^{q+1}}^{q+1}
+\int_0^t (1+\tau)^{\alpha} \, \int_{-\infty}^{\infty} 
 | \, \partial_{x}\phi\,|^{q-1} \, | \, \partial_{x}^2\phi\,|^2
 \, \mathrm{d}x \mathrm{d}\tau \\
&+\int_0^t (1+\tau)^{\alpha} \, \int_{-\infty}^{\infty} 
 | \, \partial_{x}\phi\,|^{q+1} \, \partial_{x}U^r
 \, \mathrm{d}x \mathrm{d}\tau \\
&\le C_{\phi_{0}} \, 
 \| \, \partial_{x}\phi_{0}\,\|_{L^{q+1}}^{q+1} 
 + C_{\alpha, \phi_{0}}\, 
          \, \int_0^t (1+\tau)^{\alpha-\frac{q+3}{4}}\, 
          \| \, \partial_{x}\phi(\tau)\,\|_{L^{2}}
         ^{q+1}
         \, \mathrm{d}\tau\\
&\quad +  C_{\phi_{0}}\, 
          \, \int_0^t (1+\tau)^{\alpha}\, 
          \| \, \partial_{x}\phi(\tau)\,\|_{L^{2}}
         ^{\frac{4q+6}{3}}
         \, \mathrm{d}\tau \\
&
+ C_{\phi_{0}}\, 
          \, \int_0^t (1+\tau)^{\alpha-\frac{5(q+3)}{8}}
          \, \big( \, \max \big\{ \, 1, \, \ln (1+\tau) \, \big\} \, \big)^{\frac{q+3}{2}}\, 
          \| \, \partial_{x}\phi(\tau)\,\|_{L^{2}}
         ^{\frac{q-1}{2}}
         \, \mathrm{d}\tau\\
&
+ C_{\phi_{0}, \epsilon}\, 
          \, \int_0^t 
          (1+\tau)^{\alpha}\, 
          \| \, \partial_{x}^2U^r(\tau)\,\|_{L^{2}}
         ^{\frac{q+3}{2}}
          \| \, \partial_{x}\phi(\tau)\,\|_{L^{2}}
         ^{\frac{q-1}{2}}
         \, \mathrm{d}\tau. 
\end{aligned}
\end{align}
To complete the proof of Proposition 8.1, 
we claim 

\medskip

\noindent
{\bf Lemma 8.4.} \quad {\it
For any $\epsilon>0$, there exists a positive constant $C_{\phi_{0}, \epsilon}$ 
such that 
\begin{equation}
\|\, \partial_{x}\phi(t) \,\|_{L^{2}} \leq C_{\phi_{0}, \epsilon}\, 
(1+t)^{-\frac{3}{4}+\epsilon} \quad \big( \, t \ge 0 \, \big).
\end{equation}
}

\medskip

\noindent
Indeed, once Lemma 8.4 holds true, 
by substituting (8.14) into (8.13), Proposition 8.1 immediately follows. 
In particular, we have for $q \in [\, 1, \, \infty)$
\begin{equation}
\|\, \partial_{x}\phi(t) \,\|_{L^{q+1}} \leq C_{\phi_{0}, \epsilon}\, 
(1+t)^{-\frac{2q+1}{2q+2}+\epsilon} \quad \big( \, t \ge 0 \, \big). 
\end{equation}
Therefore, we finally show Lemma 8.4 
by applying the arguments in Yoshida \cite{yoshida2}, \cite{yoshida5}. 
Taking $q=1$ to (8.13), we get 
\begin{align}
\begin{aligned}
&(1+t)^{\alpha} \, 
\| \, \partial_{x}\phi(t)\,\|_{L^{2}}^{2}
+\int_0^t (1+\tau)^{\alpha} \, 
 \| \, \partial_{x}^2\phi(\tau)\,\|_{L^{2}}^{2}
\, \mathrm{d}\tau \\
&+\int_0^t (1+\tau)^{\alpha} \, 
 \big\| \, 
 \big( \, 
 \sqrt{\partial_{x}U^r}\:\partial_{x}\phi
 \, \big)
 (\tau)
 \, \big\|_{L^{2}}^{2}
 \, \mathrm{d}\tau \\
&\le C_{\phi_{0}} \, 
 \| \, \partial_{x}\phi_{0}\,\|_{L^{2}}^{2} 
 + C_{\alpha, \phi_{0}}\, 
          \, \int_0^t (1+\tau)^{\alpha-1}\, 
          \| \, \partial_{x}\phi(\tau)\,\|_{L^{2}}^{2}
         \, \mathrm{d}\tau\\
&\quad +  C_{\phi_{0}}\, 
          \, \int_0^t (1+\tau)^{\alpha}\, 
          \| \, \partial_{x}\phi(\tau)\,\|_{L^{2}}^{\frac{10}{3}}
         \, \mathrm{d}\tau \\
&\quad + C_{\alpha, \phi_{0}}\, 
          (1+t)^{\alpha-\frac{3}{2}}
          \, \big( \, \max \big\{ \, 1, \, \ln (1+t) \, \big\} \, \big)^{2}.
\end{aligned}
\end{align}
To estimate the right-hand side of (8.16), 
we prepare a time-weighted estimate for $\phi$ from Proposition 7.7 as 

\medskip

\noindent
{\bf Lemma 8.5.} \quad {\it
There exists a positive constant $\tilde{\beta}$ 
such that for any $\beta \geq \tilde{\beta}$, 
\begin{equation*}
\, \int_0^t (1+\tau)^{\beta} \, 
\|\, \partial_{x}\phi(\tau) \,\|_{L^{2}}^2 
\, \mathrm{d}\tau \leq C_{\beta, \phi_{0}}\, 
(1+t)^{\beta-\frac{1}{2}}
\, \big( \, 
\max \big\{ \, 1, \, \ln (1+t) \, \big\} 
\, \big)^{2}  \quad \big( \, t \ge 0 \, \big).
\end{equation*}
}
By using Lemma 8.5, the second and third terms on the right-hand side of (8.16) 
are estimated as for $\alpha \geq \tilde{\beta}$, 
\begin{align}
\begin{aligned}
&C_{\alpha, \phi_{0}}\, 
          \, \int_0^t (1+\tau)^{\alpha-1}\, 
          \| \, \partial_{x}\phi(\tau)\,\|_{L^{2}}^{2}
         \, \mathrm{d}\tau \\
& \leq C_{\alpha, \phi_{0}}\, 
(1+t)^{\alpha-\frac{3}{2}}
\, \big( \, \max \big\{ \, 1, \, \ln (1+t) \, \big\} \, \big)^{2}, 
\end{aligned}
\end{align}
\begin{align}
\begin{aligned}
&C_{\phi_{0}}\, 
          \, \int_0^t (1+\tau)^{\alpha}\, 
          \| \, \partial_{x}\phi(\tau)\,\|_{L^{2}}^{\frac{10}{3}}
         \, \mathrm{d}\tau \\
         &\leq C_{\phi_{0}}\, 
              \, \int_0^t (1+\tau)^{\alpha}\, 
              \| \, \partial_{x}\phi(\tau)\,\|_{L^{2}}^{2}
              \, \mathrm{d}\tau\\
         &\leq C_{\alpha, \phi_{0}}\, 
(1+t)^{\alpha-\frac{1}{2}}
\, \big( \, \max \big\{ \, 1, \, \ln (1+t) \, \big\} \, \big)^{2}. 
\end{aligned}
\end{align}
Substituting (8.17) and (8.18) into (8.16), 
and choosing $\alpha$ suitably large, we have 
\begin{equation}
\|\, \partial_{x}\phi(t) \,\|_{L^{2}}^2 \leq C_{\phi_{0}}\, 
(1+t)^{-\frac{1}{2}} 
\, \big( \, \max \big\{ \, 1, \, \ln (1+t) \, \big\} \, \big)^{2} \quad \big( \, t \ge 0 \, \big).
\end{equation}
Substituting (8.19) into the third term on the right-hand side of (8.16), 
this becomes 
\begin{align}
\begin{aligned}
&C_{\phi_{0}}\, 
          \int_0^t (1+\tau)^{\alpha}\, 
          \| \, \partial_{x}\phi(\tau)\,\|_{L^{2}}^{\frac{10}{3}}
         \, \mathrm{d}\tau \\
         &\leq C_{\phi_{0}}\, 
              \int_0^t (1+\tau)^{\alpha-\frac{1}{2}\cdot \frac{2}{3}}\, 
              \, \big( \, \max \big\{ \, 1, \, \ln (1+\tau) \, \big\} \, \big)^{2\cdot \frac{2}{3}}
              \| \, \partial_{x}\phi(\tau)\,\|_{L^{2}}^{2}
              \, \mathrm{d}\tau\\
         &\leq C_{\alpha, \phi_{0}}\, 
(1+t)^{\alpha-\frac{1}{2}\left( \frac{2}{3} +1 \right)}
\, \big( \, \max \big\{ \, 1, \, \ln (1+t) \, \big\} \, \big)^{2\left( \frac{2}{3} +1 \right)}. 
\end{aligned}
\end{align}
Also substituting (8.17) and (8.20) into (8.16), 
and choosing $\alpha$ suitably large again, we have 
\begin{equation}
\|\, \partial_{x}\phi(t) \,\|_{L^{2}}^2 \leq C_{\phi_{0}}\, 
(1+t)^{-\frac{1}{2}\left( \frac{2}{3} +1 \right)} 
\, \big( \, \max \big\{ \, 1, \, \ln (1+t) \, \big\} \, \big)^{2\left( \frac{2}{3} +1 \right)} \quad (t \ge 0).
\end{equation}
Noting 
$$\displaystyle{ \frac{1}{2} \, \left( \, 
1+ \frac{2}{3} + \left( \, \frac{2}{3} \, \right)^2 + \left( \, \frac{2}{3} \, \right)^3 
+ \cdots  
\, \right) \leq \frac{3}{2}}
$$ 
and 
iterating ``$\infty$''-times the above process, we will get 
\begin{align}
\begin{aligned}
&C_{\phi_{0}}\, 
          \int_0^t (1+\tau)^{\alpha}\, 
          \| \, \partial_{x}\phi(\tau)\,\|_{L^{2}}^{\frac{10}{3}}
         \, \mathrm{d}\tau \\
&\leq C_{\alpha, \phi_{0}}\, 
      (1+t)^{\alpha- \mathlarger{ \frac{1}{2}} 
         \substack{{\infty }\\{\substack{{\lsum }\\{n=0}}}}
         \mathlarger{\left( \frac{2}{3} \right)^{n}} }
         \, \big( \, \max \big\{ \, 1, \, \ln (1+t) \, \big\} \, \big)
         ^{2\substack{{\infty }\\{\substack{{\lsum }\\{n=0}}}}
         \mathlarger{\left( \frac{2}{3} \right)^{n}} }\\
 &\leq C_{\alpha, \phi_{0}, \epsilon}\, 
       (1+t)^{\alpha- \frac{3}{2} + \epsilon } \quad (\epsilon>0). 
\end{aligned}
\end{align}
Therefore, substituting (8.17) and (8.22) into (8.16), 
we obtain 
\begin{align}
\begin{aligned}
&(1+t)^{\alpha} \, 
\| \, \partial_{x}\phi(t)\,\|_{L^{2}}^{2}
+\int_0^t (1+\tau)^{\alpha} \, 
 \| \, \partial_{x}^2\phi(\tau)\,\|_{L^{2}}^{2}
\, \mathrm{d}\tau \\
&+\int_0^t (1+\tau)^{\alpha} \, 
 \big\| \, 
 \big( \, 
 \sqrt{\partial_{x}U^r}\:\partial_{x}\phi
 \, \big)
 (\tau)
 \, \big\|_{L^{2}}^{2}
 \, \mathrm{d}\tau \\
&
\le C_{\phi_{0}} \, 
 \| \, \partial_{x}\phi_{0}\,\|_{L^{2}}^{2} 
+ C_{\alpha, \phi_{0}, \epsilon}\, 
          (1+t)^{\alpha-\frac{3}{2}+\epsilon} \quad (\epsilon>0).
\end{aligned}
\end{align}
Further choosing $\alpha$ suitably large to (8.23), 
we immediately have Lemma 8.4, (8.14). 
Therefore, we obtain the desired 
time-weighted $L^q$-energy estimate of $\partial_{x}\phi$ with $1\le q < \infty$.

Thus, we complete the proof of Proposition 8.1. 

\medskip

It remains to show the time-decay estimate of $\partial_{x}\phi$ 
for the case where $q=\infty$.  
To show the $L^{\infty}$-estimate of $\partial_{x}\phi$, 
Noting that 
$\| \, \partial_{x}^2\phi(t) \, \|_{L^2} \leq C_{\phi_0}$ 
for some $C_{\phi_0}>0$, 
we use the following Gagliardo-Nirenberg inequality: 
\begin{align}
\| \,\partial_{x}\phi(t) \,\|_{L^\infty }
\le C_{q,\theta}\, 
\| \,\partial_{x}\phi(t) \,\|_{L^{q+1}}^{1-\theta}
\| \,\partial _x^2 \phi(t) \,\|_{L^{2}}^{\theta}
\end{align}
for any 
$q\in [\, 1,\, \infty )$, $\theta \in (0,\,1\, ]$ satisfying 
$$
\frac{\theta}{2} =\frac{1-\theta}{q+1}. 
$$
By using (8.24), 
we immediately have 
\begin{align}
\begin{aligned}
\| \,\phi(t) \,\|_{L^\infty }
&\le C_{\theta, \phi_0, \epsilon } \, 
     (1+t)^{- \left( {\frac{2q+1}{2q+2}} - \epsilon \right)(1-\theta)}\\
&\le C_{\phi_0, \epsilon } \, 
     (1+t)^{-1 + \epsilon} \quad \big( \, t \ge 0 \, \big)
\end{aligned}
\end{align}
for $\epsilon>0$. 
Consequently, we do complete the proof of 
$L^{q+1}$-estimate of $\partial_{x}\phi$ with $1\le q \le \infty$, 
that is, 
\begin{equation}
\|\, \phi(t) \,\|_{L^{q+1}} \leq C_{\phi_{0}, \epsilon}\, 
(1+t)^{-\frac{2q+1}{2q+2}+\epsilon} \quad \big( \, t \ge 0 \, \big) 
\end{equation}
for any $\epsilon>0$.

\medskip

\noindent
{\it Remark 8.6.}\quad 
Similarly in Remark 7.9, 
although we use (9.13) instead of the uniform boundedness of 
$\| \, \partial_{x}^2\phi \, \|_{L^2}$
to get the better decay rate of $\| \, \partial_{x}\phi \, \|_{L^{\infty}}$ 
than in (8.26), 
the decay rate only becomes the same as in (8.26).

\bigskip 

\noindent
\section{Time-decay estimates I\hspace{-.1em}V}
This section is a continuation of Section 8. 
In order to accomplish the proof of Theorem 3.6, 
it suffice to show the time-decay estimates 
for the derivatives $\partial_t \phi$ and $\partial_x^2 \phi$. 
At first, we establish 
the following time-weighted $L^{2}$-energy estimate to $\partial_t \phi$. 

\medskip

\noindent
{\bf Proposition 9.1.}\quad {\it
Suppose the same assumptions as in Theorem 3.1. 
For any $q \in [\, 1, \, \infty )$, 
there exist positive constants $\alpha$,  $C_{\alpha, \phi_{0}}$
and $C_{\alpha, \phi_{0}, \epsilon}$ 
such that the unique global in time solution $\phi$ 
of the Cauchy problem {\rm (3.4)} has the following $L^{q+1}$-energy estimate 
\begin{align*}
\begin{aligned}
&(1+t)^{\alpha} \, 
\| \, \partial_{t}\phi(t)\,\|_{L^{2}}^{2}
+\int_0^t (1+t)^{\alpha} \, 
\| \, \partial_{t}\partial_{x}\phi(t)\,\|_{L^{2}}^{2} \, \mathrm{d}\tau \\
&+\int_0^t (1+\tau)^{\alpha} \, \big\| \, 
 \big( \, 
 \sqrt{\partial_{x}U^r}\:\partial_{t}\phi 
 \, \big)
 (\tau)
 \, \big\|_{L^{2}}^{2} \, \mathrm{d}\tau \\
&\le C_{\alpha, \phi_{0}} \, 
 \| \, \partial_{t}\phi_{0}\,\|_{L^{2}}^{2} 
 + C_{\alpha, \phi_{0}, \epsilon}\, 
   (1+t)^{\alpha-\frac{3}{2}+\epsilon} 
\quad \big( \, t \ge 0 \, \big), 
\end{aligned}
\end{align*}
for any $\epsilon>0$. 
}

\medskip
The proof of Proposition 9.1 is provided by the following Lemma. 

\medskip

\noindent
{\bf Lemma 9.2.}\quad {\it
For any $q \in [\, 1, \, \infty )$, 
there exist positive constants $\alpha$ and $C_{\phi_{0}}$ such that 
\begin{align}
\begin{aligned}
&(1+t)^{\alpha} \, 
\| \, \partial_{t}\phi(t)\,\|_{L^{2}}^{2}
+\int_0^t (1+\tau)^{\alpha} \, 
\| \, \partial_{t}\partial_{x}\phi(\tau)\,\|_{L^{2}}^{2} 
\, \mathrm{d}\tau \\
&+\int_0^t (1+t)^{\alpha} \, \int_{-\infty}^{\infty} 
 f''(\phi+U^r) \, 
 | \, \partial_{t}\phi\,|^{2} \, \partial_{x}U^r
 \, \mathrm{d}x \mathrm{d}\tau \\
&\le C_{\phi_{0}} \, 
 \| \, \partial_{t}\phi_{0}\,\|_{L^{2}}^{2} 
 + C_{\phi_{0}} \cdot \alpha \,  
      \int ^t_0 (1+\tau )^{\alpha -1} \, 
      \| \,\partial_{t}\phi(\tau ) \,\|_{L^{2}}^{2} 
      \, \mathrm{d}\tau\\
&\quad + C_{\phi_{0}}\, 
      \int ^t_0 (1+\tau )^{\alpha} \,
      \int_{-\infty}^{\infty} 
 | \, \partial_{x}\phi\,| \, | \, \partial_{t}\phi\,|^2
 \, \mathrm{d}x \mathrm{d}\tau \\
&\quad +  C_{\phi_{0}}\, 
      \int ^t_0 (1+\tau )^{\alpha} \,
      \big\| \, 
 \big( \, 
 \phi \, \partial_{t}U^r
 \, \big)
 (\tau)
 \, \big\|_{L^{2}}^{2} \, \mathrm{d}\tau \\
&\quad + C_{\phi_{0}}\, 
      \int ^t_0 (1+\tau )^{\alpha} \,
      \| \, \partial_{t}\partial_{x}U^r(\tau)\,\|_{L^{2}}^{2} 
      \, \mathrm{d}\tau 
   \quad \big( \, t \ge 0 \, \big). 
\end{aligned}
\end{align}
}

\medskip

In what follows, we first prove Lemma 9.2 
and finally give the proof of Proposition 9.1. 

\medskip

{\bf Proof of Lemma 9.2}. 
Differentiating the equation in (3.4) with respect to $t$, we have 
\begin{align}
\begin{aligned}
&\partial _t^2\phi 
 + \partial _t\partial_x \big( \, f(\phi+U^r) - f(U^r) \, \big) \\[5pt]
&- \partial _t\partial_x 
    \Big( \, 
    \sigma \big( \, \partial_x \phi + \partial_x U^r \, \big) 
    - \sigma \big( \, \partial_x U^r  \, \big)  \, 
    \Big)
    = \partial _t\partial_x 
     \Big( \, \sigma \big( \, \partial_x U^r  \, \big)  \, \Big). 
\end{aligned}
\end{align}
If we perform the computation 
$$
\int_0^t (1+\tau)^{\alpha} \, \int_{-\infty}^{\infty}
(9.2)\times \partial_{t}\phi
\, \mathrm{d}x \mathrm{d}\tau 
\quad ( \, \alpha>0 \, )
$$
and after integrate the resultant formula by parts, 
then we arrive at 
\begin{align}
\begin{aligned}
&\frac{1}{2} \, (1+t)^{\alpha} \, 
\| \, \partial_{t}\phi(t)\,\|_{L^{2}}^{2}\\
&
+ \int_0^t (1+\tau)^{\alpha} \, \int_{-\infty}^{\infty} 
 \partial_{t}\phi \, 
 \partial_{t}\partial_{x}\big( \,
f(\phi+U^r)-f(U^r) \, \big)
 \, \mathrm{d}x \mathrm{d}\tau \\
&- \int_0^t (1+\tau)^{\alpha} \, \int_{-\infty}^{\infty} 
 \partial_{t}\phi \, 
 \partial_{t}\partial_{x}\Big( \,
\sigma \big(\, \partial_{x}\phi+\partial_{x}U^r\, \big)
-\sigma \big(\, \partial_{x}U^r\, \big) \, \Big)
 \, \mathrm{d}x \mathrm{d}\tau \\
&=\frac{1}{2} \, 
 \| \, \partial_{x}\phi_{0}\,\|_{L^{2}}^{2} 
 + \frac{\alpha}{2} \, 
      \int ^t_0 (1+\tau )^{\alpha -1} \, 
      \| \,\partial_{t}\phi(\tau ) \,\|_{L^{2}}^{2} 
      \, \mathrm{d}\tau\\
&\quad + 
      \int ^t_0 (1+\tau )^{\alpha} \, \int_{-\infty}^{\infty} 
 \partial_{t}\phi \, 
 \partial_{t}\partial_{x}\Big( \,
\sigma \big(\, \partial_{x}U^r\, \big) \, \Big)
 \, \mathrm{d}x \mathrm{d}\tau
   \quad \big( \, t \ge 0 \, \big). 
\end{aligned}
\end{align}
By using integration by parts, the Young inequality, 
and 
the uniform boundedness of $\| \, \phi \, \|_{L^{\infty}}$ 
and $\| \, \partial_{x}\phi \, \|_{L^{\infty}}$, 
we estimate the second and third terms 
on the left-hand side of (9.3) as follows. 
\begin{align}
\begin{aligned}
&\int_0^t (1+\tau)^{\alpha} \, \int_{-\infty}^{\infty} 
 \partial_{t}\phi \, 
 \partial_{t}\partial_{x}\big( \,
f(\phi+U^r)-f(U^r) \, \big)
 \, \mathrm{d}x \mathrm{d}\tau\\
&\geq - \int_0^t (1+\tau)^{\alpha} \, \int_{-\infty}^{\infty} 
 \partial_{t}\phi \, 
 \partial_{t}\partial_{x}\phi \, f''(\phi + U^r)
 \, \mathrm{d}x \mathrm{d}\tau\\
& \quad 
  - C_{\phi_{0}}\, \int_0^t (1+\tau)^{\alpha} \, \int_{-\infty}^{\infty} 
 | \, \phi \,|\, 
 | \, \partial_{t}\partial_{x}\phi \,|\,| \, \partial_{t}U^r \,|
 \, \mathrm{d}x \mathrm{d}\tau\\    
&\geq \frac{1}{2} \, \int_0^t (1+\tau)^{\alpha} \, \int_{-\infty}^{\infty} 
 f''(\phi + U^r) \, | \, \partial_{t}\phi \,|^{2} \, 
 \partial_{x}U^r
 \, \mathrm{d}x \mathrm{d}\tau\\
& \quad 
  - \frac{1}{2} \, \int_0^t (1+\tau)^{\alpha} \, \int_{-\infty}^{\infty} 
 f''(\phi + U^r) \, | \, \partial_{x}\phi\,|  \, | \, \partial_{t}\phi \,|^{2}
 \, \mathrm{d}x \mathrm{d}\tau\\
& \quad 
  - \epsilon \, \int_0^t (1+\tau)^{\alpha} \, 
  \| \, \partial_{t}\partial_{x}\phi(\tau)\,\|_{L^{2}}^{2} 
  \,\mathrm{d}\tau\\
& \quad 
  - C_{\phi_{0}, \epsilon} \, \int_0^t (1+\tau)^{\alpha} \, 
  \big\| \, 
 \big( \, 
 \phi \, \partial_{t}U^r
 \, \big)
 (\tau)
 \, \big\|_{L^{2}}^{2} \, \mathrm{d}\tau \quad (\epsilon>0), 
\end{aligned}
\end{align}
\begin{align}
\begin{aligned}
&- \, \int_0^t (1+\tau)^{\alpha} \, \int_{-\infty}^{\infty} 
 \partial_{t}\phi \, 
 \partial_{t}\partial_{x}\Big( \,
\sigma \big(\, \partial_{x}\phi+\partial_{x}U^r\, \big)
-\sigma \big(\, \partial_{x}U^r\, \big) \, \Big)
 \, \mathrm{d}x \mathrm{d}\tau \\
&\geq C^{-1}\, \int_0^t (1+\tau)^{\alpha} \, \int_{-\infty}^{\infty} 
 \braket{ \, \partial_{x}\phi \, }^{p-1} \, 
 | \, \partial_{t}\partial_{x}\phi \,|^2
 \, \mathrm{d}x \mathrm{d}\tau\\
& \quad 
  - C_{\phi_{0}}\, \int_0^t (1+\tau)^{\alpha} \, \int_{-\infty}^{\infty} 
 | \, \partial_{t}\partial_{x}\phi \,|\,| \, \partial_{t}\partial_{x}U^r \,|
 \, \mathrm{d}x \mathrm{d}\tau\\ 
&\geq (\,C_{\phi_{0}}^{-1}-\epsilon \,)\, \int_0^t (1+\tau)^{\alpha} \, 
 \| \, \partial_{t}\partial_{x}\phi(\tau)\,\|_{L^{2}}^{2} 
  \,\mathrm{d}\tau\\
& \quad 
  - C_{\phi_{0}, \epsilon}\, \int_0^t (1+\tau)^{\alpha} \, 
\big\| \, 
 \big( \, 
 \phi \, \partial_{t}U^r
 \, \big)
 (\tau)
 \, \big\|_{L^{2}}^{2} \, \mathrm{d}\tau \quad (\epsilon>0).
\end{aligned}
\end{align}
Similarly, third term on the right-hand side of (9.3) is estimated as 
\begin{align}
\begin{aligned}
&\int ^t_0 (1+\tau )^{\alpha} \, \bigg| \, \int_{-\infty}^{\infty} 
 \partial_{t}\phi \, 
 \partial_{t}\partial_{x}\Big( \,
\sigma \big(\, \partial_{x}U^r\, \big) \, \Big)
 \, \mathrm{d}x \, \bigg| \, \mathrm{d}\tau \\
&\leq \epsilon \, \int_0^t (1+\tau)^{\alpha} \, 
\| \, \partial_{t}\partial_{x}\phi(\tau)\,\|_{L^{2}}^{2}  
\, \mathrm{d}\tau\\
& \qquad 
  + C_{\epsilon} \, \int_0^t (1+\tau)^{\alpha} \, 
  \| \, \partial_{t}\partial_{x}U^r(\tau)\,\|_{L^{2}}^{2}  
\, \mathrm{d}\tau \quad (\epsilon>0).
\end{aligned}
\end{align}
Substituting (9.4)-(9.6) into (9.3) and choosing $\epsilon$ suitably small, 
we immediately have Lemma 9.2. 

\medskip
By using Lemma 9.2, we now show Proposition 9.1. 
\medskip

{\bf Proof of Proposition 9.1}.\
To estimate each terms on the right-hand side of (9.1), 
we prepare from (8.23) that 

\medskip

\noindent
{\bf Lemma 9.3.} \quad {\it
There exists a positive constant $\tilde{\beta}$ 
such that for any $\beta \geq \tilde{\beta}$ and $\epsilon>0$, 
\begin{equation*}
\, \int_0^t (1+\tau)^{\beta} \, 
\|\, \partial_{x}^2\phi(\tau) \,\|_{L^{2}}^2 
\, \mathrm{d}\tau \leq C_{\beta, \phi_{0}, \epsilon}\, 
(1+t)^{\beta-\frac{3}{2}+\epsilon} \quad \big( \, t \ge 0 \, \big).
\end{equation*}
}

\medskip

\noindent
Noting 
$$
|\, \partial_{t}\phi \,| 
\leq C_{\phi_{0}} \, 
     \big( \, 
     |\, \phi \,| \, \partial_{x}U^r + |\, \partial_{x}\phi \,| 
     + |\, \partial_{x}^2U^r \,| + |\, \partial_{x}^2\phi \,|
     \, \big)
$$
and using Lemmas 8.4 and 9.3, 
we estimate the second term 
on the right-hand side of (9.1) as follows: 
for any $\epsilon>0$ and $\alpha \gg 1$, 
\begin{align}
\begin{aligned}
&C_{\phi_{0}} \cdot \alpha \,  
      \int ^t_0 (1+\tau )^{\alpha -1} \, 
      \| \,\partial_{t}\phi(\tau ) \,\|_{L^{2}}^{2} 
      \, \mathrm{d}\tau\\
&\leq C_{\alpha, \phi_{0}} \, 
      \int ^t_0 (1+\tau )^{\alpha -1} \, 
      \big\| \, 
 \big( \, 
 \phi \, \partial_{x}U^r
 \, \big)
 (\tau)
 \, \big\|_{L^{2}}^{2}
      \, \mathrm{d}\tau\\
& \quad + C_{\alpha, \phi_{0}} \, 
          \int ^t_0 (1+\tau )^{\alpha -1} \, 
          \| \,\partial_{x}\phi(\tau ) \,\|_{H^{1}}^{2} 
          \, \mathrm{d}\tau\\
& \quad + C_{\alpha, \phi_{0}} \, 
          \int ^t_0 (1+\tau )^{\alpha -1} \, 
          \| \,\partial_{x}^2U^r(\tau ) \,\|_{L^{2}}^{2} 
          \, \mathrm{d}\tau\\
&\leq C_{\alpha, \phi_{0}, \epsilon}\, 
          (1+t)^{\alpha-\frac{3}{2}+\epsilon}. 
\end{aligned}
\end{align}
From (9.7), we immediately have 

\medskip

\noindent
{\bf Lemma 9.4.} \quad {\it
There exists a positive constant $\tilde{\beta}$ 
such that for any $\beta \geq \tilde{\beta}$ and $\epsilon>0$, 
\begin{equation*}
\, \int_0^t (1+\tau)^{\beta} \, 
\|\, \partial_{t}\phi(\tau) \,\|_{L^{2}}^2 
\, \mathrm{d}\tau \leq C_{\beta, \phi_{0}, \epsilon}\, 
(1+t)^{\beta-\frac{1}{2}+\epsilon} \quad \big( \, t \ge 0 \, \big).
\end{equation*}
}

\medskip

By using the Lemma 9.4, 
the third term on the right-hand side of (9.1) is estimated as follows: 
for any $\epsilon>0$ and $\alpha \gg 1$, 
\begin{align}
\begin{aligned}
&C_{\phi_{0}} \,  
      \int ^t_0 (1+\tau )^{\alpha} \, 
      \int_{-\infty}^{\infty} 
 | \, \partial_{x}\phi\,| \, | \, \partial_{t}\phi\,|^2
 \, \mathrm{d}x
      \mathrm{d}\tau\\
&\leq C_{\phi_{0}} \,  
      \int ^t_0 (1+\tau )^{\alpha} \, 
 \| \, \partial_{x}\phi(\tau)\,\|_{L^{\infty}} \, 
 \| \, \partial_{t}\phi(\tau)\,\|_{L^{2}}^2
 \, \mathrm{d}\tau\\
&\leq C_{\phi_{0}, \epsilon} \,  
      \int ^t_0 (1+\tau )^{\alpha-1+\epsilon} \, 
 \| \, \partial_{t}\phi(\tau)\,\|_{L^{2}}^2
 \, \mathrm{d}\tau\\
&\leq C_{\alpha, \phi_{0}, \epsilon} \,
          (1+t )^{\alpha-\frac{3}{2}+2\epsilon}. 
\end{aligned}
\end{align}
Also the fourth and fifth terms 
on the right-hand side of (9.1) is estimated as follows: 
for any $\epsilon>0$ and $\alpha \gg 1$, 
\begin{align}
\begin{aligned}
&C_{\phi_{0}} \, \int_0^t (1+\tau)^{\alpha} \, 
  \big\| \, 
 \big( \, 
 \phi \, \partial_{t}U^r
 \, \big)
 (\tau)
 \, \big\|_{L^{2}}^{2} \, \mathrm{d}\tau\\
&\leq C_{\phi_{0}} \, \int_0^t (1+\tau)^{\alpha} \, 
  \big\| \, 
 \big( \, 
 \phi \, \partial_{x}U^r
 \, \big)
 (\tau)
 \, \big\|_{L^{2}}^{2} \, \mathrm{d}\tau\\
&\leq C_{\alpha, \phi_{0}} \,
          (1+t )^{\alpha-\frac{3}{2}}
          \, \big( \, \max \big\{ \, 1, \, \ln (1+t) \, \big\} \, \big)^{2}, 
\end{aligned}
\end{align}
\begin{align}
\begin{aligned}
&C_{\phi_{0}} \, \int_0^t (1+\tau)^{\alpha} \, 
  \| \, \partial_{t}\partial_{x}U^r(\tau)\,\|_{L^{2}}^{2}
   \, \mathrm{d}\tau\\
&\leq C_{\phi_{0}} \, \int_0^t (1+\tau)^{\alpha} \, 
  \Big( \, 
  \| \, \partial_{x}U^r(\tau)\,\|_{L^{4}}^{4} 
  + \| \, \partial_{x}^2U^r(\tau)\,\|_{L^{2}}^{2} 
  \, \Big)
  \, \mathrm{d}\tau\\
&\leq C_{\alpha, \phi_{0}} \,
          (1+t )^{\alpha-\frac{3}{2}}. 
\end{aligned}
\end{align}
Substituting (9.7)-(9.10) into (9.1), 
we obtain the desired 
time-weighted $L^2$-energy estimate of $\partial_{t}\phi$. 

Thus, the proof of Proposition 9.1 is completed. 

\medskip

To complete the proof of Theorem 3.8, 
we finally obtain the time-decay estimate of $\partial_{x}^2\phi$. 
By using that for any $\epsilon>0$, 
\begin{equation}
\|\, \partial_{t}\phi(t) \,\|_{L^{2}} 
\leq C_{\phi_{0}, \epsilon}\, 
(1+t)^{-\frac{3}{4}+\epsilon} \quad \big( \, t \ge 0 \, \big)
\end{equation}
from Proposition 9.1, 
we obtain the time-decay estimate of $\partial_{x}^2\phi$ as follows. 
\begin{align}
\begin{aligned}
&\|\, \partial_{x}^2\phi(t) \,\|_{L^{2}}^2 \\
& \leq C_{\phi_{0}}\, 
       \Big( \, 
  \|\, \partial_{t}\phi(t) \,\|_{L^{2}}^2 
  +\big\| \, 
 \big( \, 
 \phi \, \partial_{x}U^r
 \, \big)
 (t)
 \, \big\|_{L^{2}}^{2}
  + \|\, \partial_{x}\phi(t) \,\|_{L^{2}} 
  + \| \, \partial_{x}^2U^r(t)\,\|_{L^{2}}^{2} 
  \, \Big)\\
&\leq C_{\phi_{0}, \epsilon}\, 
(1+t)^{-\frac{3}{2}+\epsilon} \quad (\epsilon>0), 
\end{aligned}
\end{align}
that is, 
\begin{equation}
\|\, \partial_{x}^2\phi(t) \,\|_{L^{2}} 
\leq C_{\phi_{0}, \epsilon}\, 
(1+t)^{-\frac{3}{4}+\epsilon} \quad \big( \, t \ge 0 \, \big)
\end{equation}
for any $\epsilon>0$. 

Thus, the proof of Theorem 3.8 is completed.

\bigskip 

\noindent

\section{Conservation law wothout convective flux} 
In this section, we finally consider 
the global structure of 
the conservation law without convective flux, 
and compare the structure with theorems in Section 1. 

At first, we consider 
the following Cauchy problem 
for the non-convective viscous conservation law 
\begin{eqnarray}
 \left\{\begin{array}{ll}
  \partial_tu -\partial_x \big( \, \sigma(\partial_xu )\, \big) =0
  \qquad &\big( \, t>0, \: x\in \mathbb{R} \, \big), \\[5pt]
  u(0, x) = u_0(x) \rightarrow  \tilde{u}  \qquad &( x \rightarrow  \pm \infty ).   
 \end{array}
 \right.\,
\end{eqnarray}
By the similar arguments as on previous sections, 
because we need controll the terms from convective flux no longer, 
such as 
the first term on the right-hand side of (4.1) 
and the second term on the left-hand side of (5.2), 
we particularly note that, in any cases for $p>0$, we can get
the uniform boundedness of 
$\| \, \partial_{x}u \, \|_{L^{\infty}}$ (cf. Lemma 5.3). 
Namely, for any $p>0$, 
$$
\|\,\partial_xu (t)\,\|_{L^{\infty}} \le C_{u_{0}} 
\quad \big(\,t \in [\,0,\,T\,] \, \big). 
$$
Therefore, we can obtein the global structure for any $p>0$ of the solution 
to (10.1) as follows. 

\medskip

\noindent
{\bf Theorem 10.1.}\quad{\it
Assume the viscous flux 
$\sigma\in C^2(\mathbb{R})$, {\rm(1.2)}, {\rm(1.3)} 
and $p>0$. 
Further assume the initial data satisfy
$u_0-\tilde{u} \in H^2$.
Then the Cauchy problem {\rm(1.1)} has a 
unique global in time 
solution $u$ 
satisfying 
\begin{eqnarray*}
\left\{\begin{array}{ll}
u-\tilde{u} \in C^0\cap L^{\infty}
( [\, 0, \, \infty \, ) \, ; H^2 \,),\\[5pt]
\partial _x u \in L^2\bigl( \, 0,\, \infty \, ; H^2 \,\bigr),\\[5pt]
\partial _t u \in C^0
                  \cap L^{\infty}\bigl( \, [\, 0,\, \infty \, ) \, ; L^2 \,\bigr)
                  \cap L^2\bigl( \, 0,\, \infty \, ; H^1 \,\bigr),
\end{array} 
\right.\,
\end{eqnarray*}
and the asymptotic behavior 
\begin{eqnarray*}
\left\{\begin{array}{ll}
\displaystyle{
\lim _{t \to \infty}\sup_{x\in \mathbb{R}} \, 
\left| \, u(t,x) - \tilde{u} \, \right| = 0, 
\quad 
\lim _{t \to \infty}\sup_{x\in \mathbb{R}} \, 
\left| \, \partial_{x}u(t,x) \, \right| = 0
},\\[10pt]
\displaystyle{
\lim _{t \to \infty} \, 
\| \, \partial_{x}u(t) \, \|_{H^1} = 0, 
\quad 
\lim _{t \to \infty} \, 
\| \, \partial_{t}u(t) \, \|_{L^2} = 0. 
}
\end{array} 
\right.\,
\end{eqnarray*}
}

\medskip 

The decay results correspond to Theorem 10.1 are the following. 

\medskip

\noindent
{\bf Theorem 10.2.} \quad{\it
Under the same assumptions as in Theorem 10.1, 
the unique global in time 
solution $u$ 
of the Cauchy problem {\rm(10.1)} 
has the following time-decay estimates 
\begin{eqnarray*}
\left\{\begin{array}{ll}
\| \, u(t)-\tilde{u} \, \|_{L^q}
\le C_{u_{0}} \, 
(1+t)^{-\frac{1}{4}\left( 1-\frac{2}{q} \right)}
\quad \big( \, t\ge0 \, \big),\\[5pt]
\| \, u(t)-\tilde{u} \, \|_{L^{\infty}}
\le C_{u_{0}, \epsilon} \, (1+t)^{-\frac{1}{4}+\epsilon}
\quad \big( \, t\ge0 \, \big),
\end{array} 
\right.\,
\end{eqnarray*}
for $q \in \ [\, 2, \, \infty)$ and any $\epsilon>0$. 
}

\medskip

\noindent
{\bf Theorem 10.3.} \quad{\it
Under the same assumptions as in Theorem 10.1, 
if the initial data further satisfies $u_0-\tilde{u}\in L^1$, 
then the unique global in time 
solution $u$ 
of the Cauchy problem {\rm(10.1)} 
has the following time-growth and time-decay estimates 
\begin{eqnarray*}
\left\{\begin{array}{ll}
\| \, u(t)-\tilde{u} \, \|_{L^1}
\le \| \, u_{0}-\tilde{u} \, \|_{L^1}
\quad \big( \, t\ge0 \, \big),\\[5pt]
\| \, u(t)-\tilde{u} \, \|_{L^q}
\le C_{u_{0}} \, 
(1+t)^{-\frac{1}{2}\left( 1-\frac{1}{q} \right)} 
\quad \big( \, t\ge0 \, \big),\\[5pt]
\| \, u(t)-\tilde{u} \, \|_{L^{\infty}}
\le C_{u_{0}, \epsilon} \, (1+t)^{-\frac{1}{2}+\epsilon}
\quad \big( \, t\ge0 \, \big),
\end{array} 
\right.\,
\end{eqnarray*}
for $q \in \ ( 1, \, \infty)$ and any $\epsilon>0$. 
}

\medskip

\noindent
{\bf Theorem 10.4.} \quad{\it
Under the same assumptions as in Theorem 10.1, 
the unique global in time 
solution $u$ 
of the Cauchy problem {\rm(10.1)} 
has the following time-decay estimates 
for the derivatives 
\begin{eqnarray*}
\left\{\begin{array}{ll}
\| \, \partial_{x}u(t) \, \|_{L^{q+1}}
\le C_{u_{0}} \, 
(1+t)^{-\frac{2q+1}{2q+2}}
\quad \big( \, t\ge0 \, \big),\\[5pt]
\| \, \partial_{x}u(t) \, \|_{L^{\infty}}
\le C_{u_{0}, \epsilon} \, 
(1+t)^{-1+\epsilon} 
\quad \big( \, t\ge0 \, \big),\\[5pt]
\| \, \partial_{t}u(t) \, \|_{L^{2}}
\le C_{u_{0}} \, (1+t)^{-\frac{3}{4}}
\quad \big( \, t\ge0 \, \big),\\[5pt]
\| \, \partial_{x}^2u(t) \, \|_{L^{2}}
\le C_{u_{0}} \, (1+t)^{-\frac{3}{4}}
\quad \big( \, t\ge0 \, \big),
\end{array} 
\right.\,
\end{eqnarray*}
for $q \in \ ( 1, \, \infty)$ and any $\epsilon>0$. 
}

\medskip
From above Theorems 10.1-10.4, 
we easily see that the global structures in Theorems 1.1-1.8 
are almost the same as in the case where $f(u)\equiv 0$, 
that is Theorems 10.1-10.4 (see also Remarks 1.9 and 1.10). 

\medskip
Next, we compare the decay results in theorems in both Sections 1 and 10 with 
the viscous contact wave (1.11) 
considered in \cite{matsumura-yoshida}, \cite{yoshida1}, \cite{yoshida7}, 
as a typical example of the exact solution 
to the non-convective conservation law, 
that is 
\begin{align}
\begin{aligned}
&U\left(\frac{x}{\sqrt{t}}\: ;\: u_- ,\: u_+ \right)
 :=u_- +\frac{u_+ - u_-}{\sqrt{\pi}}
   \lint ^{\frac{\mathlarger{x}}{\mathlarger{\sqrt{4\mu t}}}}_{-\infty} 
 \mathrm{e}^{-\xi^2}\, \mathrm{d}\xi, 
\end{aligned}
\end{align}
where $\mu>0$ and $u_{-}<u_{+}$. 
The above viscous contact wave is the unique $C^{\infty}$-solution to the 
following Cauchy problem 
\begin{equation}
           \left\{
              \begin{array}{l} 
              \partial _tU - \mu \, \partial _x^2U = 0 
              \qquad \qquad \qquad \qquad  
              \bigl(\, t>0,\:  x\in \mathbb{R}\, \bigr),\\[5pt]
              U(0,x) = u_0 ^{\rm{R}} (\, x\: ;\: u_- ,\: u_+)
               = 
               \left\{\begin{array} {ll}
               u_-  & \, \; \qquad (x < 0),\\[5pt]
               u_+  & \, \; \qquad (x > 0),
               \end{array}\right.\\[15pt]
              \displaystyle{\lim_{x\to \pm \infty}} U(t,x) =u_{\pm} 
              \, \, \: \: \: \; \quad \qquad \qquad \qquad \qquad 
              \bigl(\, t\ge 0 \, \bigr),
              \end{array}
            \right.\,  
\end{equation}
and has the following optimal time-decay properties. 
             \begin{eqnarray}
                 \begin{array}{l}
                    \| \, \partial_x U(t) \, \|_{L^r} 
                    \sim (1+t)^{-\frac{1}{2}\left( 1- \frac{1}{r} \right)} 
                    \quad \bigl(\, t\ge 0 \, \bigr),\\[5pt]
                    \| \, \partial_t U(t) \, \|_{L^2} 
                    \sim (1+t)^{-\frac{3}{4}}
                    \quad \bigl(\, t\ge 0 \, \bigr),\\[5pt]
                    \| \, \partial_x^2 U(t) \, \|_{L^2} 
                    \sim (1+t)^{-\frac{3}{4}}
                    \quad \bigl(\, t\ge 0 \, \bigr),
                    \end{array}       
              \end{eqnarray}
for any $r \in [\,1,\, \infty \,]$, 
where we note that these time-decay rates 
fully coinside with the time-decay rates of the viscous contact wave (1.12) 
of the $p$-Laplacian evolution equation 
considered in \cite{yoshida3}, \cite{yoshida4}, \cite{yoshida5}, when $p=1$ 
(for $p>1$, the precise time-decay properties are given in \cite{yoshida5}). 

\medskip
From (10.4), we also easily see that the time-decay rates 
in $L^2$-norm of $\partial_t U$ and $\partial_x^2 U$ are almost the same as 
in $L^2$-norm of $\partial_t u$ and $\partial_x^2 u$ in Theorems 1.4 and 10.4, 
and in $L^2$-norm of $\partial_x^2 u-\partial_x^2 u^r$, $\partial_t u-\partial_t U^r$ 
and $\partial_x^2 u-\partial_x^2 U^r$ in Theorem 1.8. 
Besides, it is interesting that the time-decay rates 
in $L^r$-norm with $1\le r \le \infty$ of $\partial_x U$ is almost the same as 
in $L^r$-norm of $u$ and $u-U^r$ 
in Theorems 1.3, 1.7 and 10.3 (see also the arguments 
in \cite{vaz1}, \cite{vaz2}, \cite{yoshida2}).

\bigskip









\bibliographystyle{model6-num-names}
\bibliography{<your-bib-database>}

\begin{thebibliography}{99}

\bibitem{barenblatt}
G. I. Barenblatt, 
{\it On the motion of suspended particles in a turbulent flow 
taking up a half-space or a plane open channel of finite depth},
Prikl. Mat. Meh., {\bf19} (1955), pp. 61-88 (in Russian).

\bibitem{carillo-toscani}
J. A. Carrillo and G. Toscani, 
{\it Asymptotic $L^1$-decay of solutions of the porous medium equation 
to self-similarity}, 
Indiana Univ. Math. J., {\bf 49} (2000), pp. 113-142.

\bibitem{chh1}
R. P. Chhabra, 
{\it Bubbles, drops and particles in non-Newtonian Fluids}, 
CRC, Boca Raton, FL, 2006.

\bibitem{chh2}
R. P. Chhabra, 
{\it Non-Newtonian Fluids: An Introduction}, 
URL http://www.physics.iitm.ac.in/$\tilde{\: }$compflu/Lect-notes/chhabra.pdf.

\bibitem{chh-ric}
R. P. Chhabra and J. F. Richardson,  
{\it Non-Newtonian flow and applied rheology}, 
2nd edn. Butterworth-Heinemann, Oxford, 2008.

\bibitem{cra-tar} 
M. Crandall and L. Tartar, 
{\it Some relations between nonexpansive 
and order preserving mappings}, 
Proc. Amer. Math. Soc., {\bf 78} (1980), pp. 385-390.

\bibitem{de waele} 
A. de Waele, 
{\it Viscometry and plastometry}, 
J. Oil Colour Chem. Assoc., {\bf 6} (1923), pp. 3369.

\bibitem{den-wan} 
S. Deng and W. Wang, 
{\it Pointwise decaying rate of large perturbation around 
viscous shock for scalar viscous conservation law}, 
Sci. Chin. Math., {\bf 56} (2013), pp. 729-736.

\bibitem{du-gu} 
Q. Du and M. D. Gunzburger, 
{\it Analysis of a Ladyzhenskaya model for incompressible viscous flow}, 
J. Math. Anal. Appl., {\bf 155} (1991), pp. 21-45.

\bibitem{fre-ser} 
H. Freist{\"{u}}hler and D. Serre 
{\it $L^1$-stability of shock waves in scalar viscous conservation laws}, 
Comm. Pure Appl. Math., {\bf 51} (1998), pp. 291-301. 

\bibitem{har} 
E. Harabetian, 
{\it Rarefactions and large time behavior for parabolic equation 
and monotone schemes}, 
Comm. Math. Phys., {\bf 114} (1988), pp. 527-536.

\bibitem{hashimoto-kawashima-ueda} 
I. Hashimoto, Y. Ueda and S. Kawashima, 
{\it Convergence rate to the nonlinear waves 
for viscous conservation laws on the half line},
Methods Appl. Anal., {\bf 16} (2009), pp. 389-402.

\bibitem{hashimoto-matsumura} 
I. Hashimoto and A. Matsumura, 
{\it Large time behavior 
of solutions to an initial boundary value problem 
on the half space for scalar viscous conservation law}, 
Methods Appl. Anal., {\bf 14} (2007), pp. 45-59.

\bibitem{hattori-nishihara} 
Y. Hattori and K. Nishihara, 
{\it A note on the stability of rarefaction wave of the 
Burgers equation}, 
Japan J. Indust. Appl. Math., {\bf 8} (1991), pp. 85-96.

\bibitem{huang-pan-wang} 
F. Huang, R. Pan and Z. Wang, 
{\it $L^1$Convergence to the Barenblatt solution for 
compressible Euler equations with damping}, 
Arch. Rational Mech. Anal., {\bf 200} (2011), pp. 665-689.

\bibitem{kur-lev-ros} 
A. Kurganov, D. Levy and P. Rosenau, 
{\it On Burgers-type equations with nonmonotonic dissipative fluxes}, 
Comm. Pure Appl. Math., {\bf 51} (1998), pp. 443-473. 

\bibitem{ilin-kalashnikov-oleinik} 
A. M. Il'in, A. S. Kala{\v{s}}nikov and O. A. Ole{\u\i}nik,
{\it Second-order linear equations of parabolic type}, 
Uspekhi Math. Nauk SSSR, {\bf 17} (1962), pp. 3-146 (in Russian); 
Russian Math. Surveys, {\bf 17} (1962), pp. 1-143 (in English). 

\bibitem{ilin-oleinik} 
A. M. Il'in and O. A. Ole{\u\i}nik, 
{\it Asymptotic behavior of the solutions of the Cauchy problem for 
some quasi-linear equations for large values of the time}, 
Mat. Sb., {\bf 51} (1960), pp. 191-216 (in Russian).

\bibitem{jah-str-mul} 
P. Jahangiri, R. Streblow and D. M\"{u}ller, 
{\it Simulation of Non-Newtonian Fluids using Modelica}, 
Proceedings of the 9th International Modelica Conference September 3-5, 
Munich, Germany, (2012), pp. 57-62. 

\bibitem{kamin} 
S. Kamin, 
{\it Source-type solutions for equations of nonstationary filtration}, 
J. Math. Anal. Appl., {\bf 64} (1978), pp. 263-276. 

\bibitem{kanel} 
Ya. I. Kanel', 
{\it A model system of equations 
for the one-dimensional motion of a gas}, 
Differencial'nya Uravnenija, {\bf 4} (1968), pp. 721-734 (in Russian).

\bibitem{kato1} 
T. Kato, 
{\it Linear evolution equations of ``hyperbolic type''}, 
J. Fac. Sci. Univ. Tokyo Sect. A, {\bf 17} (1970), pp. 241-258. 

\bibitem{kato2} 
T. Kato, 
{\it Linear evolution equations of ``hyperbolic type'', I\hspace{-.1em}I}, 
J. Math. Soc. Japan, {\bf 19} (1973), pp. 648-666.


\bibitem{lad} 
O. A. Lady{\v{z}}enskaja, 
{\it New equations for the description of the viscous incompressible
fluids and solvability in the large of the boundary value problems for them}, 
in Boundary Value Problems of Mathematical Physics V, 
AMS, 
Providence, RI, 1970.

\bibitem{lad-sol-ura} 
O. A. Lady{\v{z}}enskaja, V. A. Solonnikov and N. N. Ural'ceva, 
{\it Linear and quasilinear equations of parabolic type}, 
Transl. Math. Monogr. {\bf 23}, AMS, 
Providence, RI, 1968.

\bibitem{lax} 
P. D. Lax, 
{\it Hyperbolic systems of conservation laws  I\hspace{-.1em}I}, 
Comm. Pure Appl. Math., {\bf 10} (1957), pp. 537-566.

\bibitem{liep-rosh} 
H. W. Liepmann and A. Roshko, 
{\it Elements of Gas Dynamics}, 
John Wiley \& Sons, Inc., New York, 1957.

\bibitem{lions} 
J. L. Lions, 
{\it Quelques m{\`{e}}thodes de r{\`{e}}solution 
des probl{\`{e}}mes aux limites non lin{\`{e}}aires}, 
Dunod Gauthier-Villars, Paris, 1969 (in French).

\bibitem{liu-matsumura-nishihara} 
T.-P. Liu, A. Matsumura and K. Nishihara, 
{\it Behaviors of solutions for the Burgers equation 
with boundary corresponding to rarefaction waves}, 
SIAM J. Math. Anal., {\bf 29} (1998), pp. 293-308.

\bibitem{ma} 
J. M{\'{a}}lek, 
{\it Some frequently used models for non-Newtonian fluids}, 
URL http://www.karlin.mff.cuni.cz/$\tilde{\: }$malek/new/images/Lecture4.pdf.

\bibitem{ma-pr-st} 
J. M{\'{a}}lek, D. Pra{\v{z}}{\'{a}}k and M. Steinhauer,
{\it On the existence and regularity of solutions 
for degenerate power-law fluids}, 
Differential Integral Equations, {\bf 19} (2006), pp. 449-462.


\bibitem{matsumura} 
A. Matsumura, 
{\it Waves in compressible fuids: 
viscous shock, rarefaction and contact waves}, 
in Handbook of Mathematical Analysis in Mechanics of Viscous Fluids, 
Springer-Verlag, 
New York, 2018. 

\bibitem{matsu-nishi1} 
A. Matsumura and K. Nishihara,
{\it Asymptotic toward the rarefaction wave of solutions of a
one-dimensional model system for compressible viscous gas}, Japan
J. Appl. Math., {\bf 3} (1986), pp. 1-13.

\bibitem{matsu-nishi2-0} 
A. Matsumura and K. Nishihara,
{\it Global stability of the rarefaction wave of a 
one-dimensional model system for compressible viscous gas}, 
Comm. Math. Phys. {\bf 144} (1992), pp. 325-335.

\bibitem{matsu-nishi2} 
A. Matsumura and K. Nishihara, 
{\it Asymptotics toward the rarefaction wave 
of the solutions of Burgers' equation with 
nonlinear degenerate viscosity}
Nonlinear Anal., {\bf 23} (1994), pp. 605-614.

\bibitem{matsu-nishi3} 
A. Matsumura and K. Nishihara, 
{\it Asymptotic stability of traveling waves for scalar 
viscous conservation laws with non-convex nonlinearity}, 
Comm. Math. Phys., {\bf 165} (1994), pp. 83-96.

\bibitem{matsumura-yoshida} 
A. Matsumura and N. Yoshida, 
{\it Asymptotic behavior of solutions to the Cauchy problem 
for the scalar viscous conservation law 
with partially linearly degenerate flux}, 
SIAM J. Math. Anal., {\bf 44} (2012), pp. 2526-2544.

\bibitem{matsumura-yoshida'} 
A. Matsumura and N. Yoshida, 
{\it Global asymptotics toward 
the rarefaction waves for solutions to the Cauchy problem of
the scalar conservation law with nonlinear viscosity}
(to appear in Osaka Math. J.)
/{\it Global asymptotics toward rarefaction waves 
for solutions of the scalar conservation law with nonlinear viscosity}, 
Preprint, arXiv:1804.10841.



\bibitem{osh-ral} 
S. Osher and J. Ralston 
{\it $L^1$ stability of traveling waves with applications to 
convective porous media flow}, 
Comm. Pure Appl. Math., {\bf 35} (1982), pp. 737-751.

\bibitem{ost} 
W. Ostwald, 
{\it \"{U}ber die Geschwindigkeitsfunktion der Viskositat disperser Systeme}, 
I. Colloid Polym. Sci., {\bf 36} (1925), pp. 99-117 (in German).

\bibitem{pattle} 
R. E. Pattle, 
{\it Diffusion from an instantaneous point source with 
a concentration-dependent coefficient}, 
Quart. J. Mech. Appl. Math., {\bf 12} (1959), pp. 407-409.

\bibitem{smoller} 
J. Smoller, 
{\it Shock Waves and Reaction-diffusion Equations}, Springer-Verlag, 
New York, 1983. 

\bibitem{soc} 
T. Sochi, 
{\it Pore-Scale Modeling of Non-Newtonian Flow in Porous Media}, 
PhD thesis, Imperial College London, 2007. 

\bibitem{vaz1} 
J. L. V{\'{a}}zquez, 
{\it Smoothing and Decay Estimates for Nonlinear Diffusion Equations: 
Equations of Porous Medium Type}, Oxford Math. and Appl., 2006. 

\bibitem{vaz2} 
J. L. V{\'{a}}zquez, 
{\it The Porous Medium Equation: Mathematical Theory}, 
Oxford Math. Monogr., 2007. 

\bibitem{yoshida1} 
N. Yoshida, 
{\it Decay properties of solutions toward a multiwave pattern 
for the scalar viscous conservation law 
with partially linearly degenerate flux}, 
Nonlinear Anal., {\bf 96} (2014), pp. 189-210.

\bibitem{yoshida2} 
N. Yoshida, 
{\it Decay properties of solutions 
to the Cauchy problem for the scalar conservation law 
with nonlinearly degenerate viscosity}, 
Nonlinear Anal., {\bf 128} (2015), pp. 48-76.

\bibitem{yoshida3} 
N. Yoshida, 
{\it Large time behavior of solutions toward 
a multiwave pattern for the Cauchy problem 
of the scalar conservation law with degenerate flux and viscosity}, 
S\={u}rikaisekikenky\={u}sho K\={o}ky\={u}roku 
gMathematical Analysis in Fluid and Gas Dynamicsh, 
{\bf 1947} (2015), pp. 205-222.

\bibitem{yoshida4} 
N. Yoshida, 
{\it Asymptotic behavior of solutions toward a multiwave pattern 
for the scalar conservation law 
with the Ostwald-de Waele-type viscosity}, 
SIAM J. Math. Anal., {\bf 49} (2017), pp. 2009-2036.

\bibitem{yoshida5} 
N. Yoshida, 
{\it Decay properties of solutions toward a multiwave pattern 
to the Cauchy problem for the scalar conservation law 
with degenerate flux and viscosity}, 
J. Differential Equations, {\bf 263} (2017), pp. 7513-7558.

\bibitem{yoshida6} 
N. Yoshida, 
{\it Asymptotic behavior of solutions toward the viscous shock waves 
to the Cauchy problem 
for the scalar conservation law with nonlinear flux and viscosity}, 
SIAM J. Math. Anal., {\bf 50} (2018), pp. 891-932.

\bibitem{yoshida7} 
N. Yoshida, 
{\it Asymptotic behavior of solutions toward a multiwave pattern
to the Cauchy problem for the dissipative wave equation 
with partially linearly degenerate flux} 
(to appear in Funkcialaj Ekvacioj).

\bibitem{yoshida8} 
N. Yoshida, 
{\it Asymptotic behavior of solutions toward 
the rarefaction waves to the Cauchy problem for 
the scalar diffusive dispersive conservation laws}, 
Nonlinear Anal., {\bf 189} (2019), pp. 1-19.

\bibitem{zel-kom} 
Ya. B. Zel'dovi{\v{c}} and A. S. Kompaneec, 
{\it On the theory of propagation 
of heat with the heat conductivity depending upon the temperature}, 
in Collection in honor of the seventieth birthday of academician A. F. Ioffe, 
Izdat. Akad. Nauk SSSR, 1950, pp. 61-71 (in Russian). 

\end{thebibliography}







\end{document}